\RequirePackage{fix-cm}
\documentclass[smallextended]{svjour3}       % onecolumn (second format)
\smartqed  % flush right qed marks, e.g. at end of proof
\usepackage{graphicx}%
\usepackage{multirow}%
\usepackage{amsmath,amssymb,amsfonts}%
\usepackage{mathrsfs}%
\usepackage[title]{appendix}%
\usepackage{xcolor}%
\usepackage{textcomp}%
\usepackage{manyfoot}%
\usepackage{booktabs}%
\usepackage{algorithm}%
\usepackage{algorithmicx}%
\usepackage[noEnd]{algpseudocodex}
\usepackage{listings}%
\usepackage{comment}%
\usepackage{enumitem}% for costum enumerates
\usepackage{hyperref}

\usepackage{subfig}
\graphicspath{{figures}}
 \newtheorem{mydefinition}{Definition}[section]
 \numberwithin{mydefinition}{section}
 
 \newtheorem{mytheorem}{Theorem}[section]
 \numberwithin{mytheorem}{section}
 
 \newtheorem{mylemma}{Lemma}[section]
 \numberwithin{mylemma}{section}
 
 \newtheorem{myproposition}{Proposition}[section]
 \numberwithin{myproposition}{section}
 
\newtheorem{myassumption}{Assumption}[section]
\numberwithin{myassumption}{section}

\newtheorem{mycorollary}{Corollary}[section]
\numberwithin{mycorollary}{section}

\newtheorem{myremark}{Remark}[section]
\numberwithin{myremark}{section}

 \newtheorem{myexample}{Example}[section]
 \numberwithin{myexample}{section}

\newcommand{\R}{{\mathbb{R}}}
\newcommand{\T}{{\mathcal{T}}}
\newcommand{\I}{{\mathcal{I}}}
\newcommand{\D}{{\mathcal{D}}}
\newcommand{\A}{{\mathcal{A}}}

\newcommand{\BNLP}{{\mathrm{BNLP}}}

\newcommand{\LPEC}{{\mathrm{LPEC}}}
\newcommand{\LPCC}{{LPEC}}
\newcommand{\LPCCs}{{LPECs}}
\newcommand{\MPEC}{{MPEC}}
\newcommand{\MPECs}{{MPECs}}
\newcommand{\MPECopt}{{MPECopt}}

\newcommand{\ipopt}{\texttt{IPOPT}}
\newcommand{\casadi}{\texttt{CasADi}}
\newcommand{\gurobi}{\texttt{Gurobi}}
\newcommand{\highs}{\texttt{HiGHS}}

\newcommand{\MacMPEC}{\texttt{MacMPEC}}

\begin{document}
%\title[A Globally Convergent Method for Computing B-stationary Points of Mathematical Programs with Complementarity Constraints]
%{A Globally Convergent Method for Computing B-stationary Points of  Mathematical Programs with Complementarity Constraints}

\title{A Globally Convergent Method for Computing B-stationary Points of Mathematical Programs with Equilibrium Constraints}
\titlerunning{A Globally Convergent Method for Computing B-stationary Points of MPECs}% Part of RIGHT running header

\author{Armin Nurkanovi\'c        \and
	Sven Leyffer
}
\institute{Armin Nurkanovi\'c\at
	Systems Control and Optimization Laboratory, Department of Microsystems Engineering\\ 
	{University of Freiburg, Germany}\\
	\email{armin.nurkanovic@imtek.uni-freiburg.de}\\
	\and
	Sven Leyffer \at
	Mathematics and Computer Science Division\\
	Argonne National Laboratory, USA\\
	\email{leyffer@anl.gov}
}

\date{Received: 07.2026/ Accepted: xx}

\maketitle
\begin{abstract}
This paper introduces a computationally efficient method that converges globally to B-stationary points of mathematical programs with equilibrium constraints \sloppy(\MPECs).
B-stationarity is necessary for optimality and means that no feasible first-order direction can improve the objective. 
It can be certified by solving a linear program with equilibrium constraints (\LPCC) constructed at a given feasible point.
The proposed method solves a finite sequence of \LPCCs, which either certify B-stationarity or provide an active-set estimate for the complementarity constraints, along with branch nonlinear programs (BNLPs) obtained by fixing the complementarity active set in the \MPEC.
In particular, the method proceeds in two phases: the first identifies a feasible BNLP or a stationary point of a constraint infeasibility minimization problem, and the second solves a sequence of BNLPs until a B-stationary point of the {\MPEC} is found.
We prove that under the \MPEC-MFCQ, the method requires solving only a finite number of BNLPs and {\LPCCs} for convergence.
Moreover, we show that, unless the current iterate is B-stationary, the combinatorial {\LPCCs} need not be solved to optimality.
For convergence, it suffices to compute a nonzero feasible point, which in practice often requires solving a single linear program, yielding significant computational savings.
Numerical experiments show that the proposed method is more robust and faster than relaxation-based methods and mixed-integer NLP reformulations (which, in contrast to the proposed approach, do not provide a certificate of B-stationarity), even on medium- to large-scale instances.
		
\keywords{MPEC \and MPCC \and nonlinear programming \and complementarity constraints \and equilibrium constraints}
	\subclass{90C30 \and 90C33 \and 49M37 \and 65K10\and 90C11  }
\end{abstract}

\section{Introduction}\label{sec:introduction}
We study mathematical programs with equilibrium constraints (\MPECs) of the following form:
\begin{subequations}\label{eq:mpec}
\begin{align}
	\underset{x\in \R^{n}}{\mathrm{min}} \;  \quad &f(x)\\
	\textnormal{s.t.} \quad 
%	&g(x)=0, \label{eq:mpec_eq}\\
	&c(x)\geq0, \label{eq:mpec_ineq}\\
	&0 \leq x_1 \perp x_2 \geq 0 \label{eq:mpec_comp},
\end{align}
\end{subequations}
with the partition of variables $x=(x_0,x_1,x_2) \in \R^n$, where $x_0 \in \R^p$, $x_1,x_2\in \R^m$. 
The functions $f:\R^n \to \R$ and $c:\R^n \to \R^{n_c} $ are assumed to be at least twice continuously differentiable.
The notation $0 \leq x_1 \perp x_2 \geq 0$ means that $x_{1,i},x_{2,i} \geq 0$ and $x_{1,i}x_{2,i}=0$ for all $i$.
{There is no loss of generality in using \eqref{eq:mpec_comp}, since the more general form $0 \leq g(x) \perp h(x) \geq 0$, with $g,h:\R^{p}\to\R^{m}$ twice continuously differentiable, can be transformed into \eqref{eq:mpec_comp} by introducing slack variables $x_1 = g(x)$ and $x_2 = h(x)$.}
For lighter notation we omit equality constraints in~\eqref{eq:mpec}. 
If given by twice continuously differentiable functions, all results extend directly to this case.

Equilibrium constraints are more general than complementarity constraints \eqref{eq:mpec_comp}, but coincide under suitable conditions~\cite{Luo1996}. 
We focus on complementarity constraints, yet follow the literature in using the acronym MPEC.

Solving \MPECs~reliably and fast is of great practical interest, as they arise in numerous applications including process engineering~\cite{Baumrucker2008}, robotics~\cite{Wensing2023}, optimal control of hybrid dynamical systems~\cite{Nurkanovic2023f,Nurkanovic2022b}, bilevel optimization~\cite{Kim2020}, nonsmooth optimization~\cite{Hegerhorst2020}, inverse optimization~\cite{Albrecht2017,Hu2012}, to name a few examples.  
A related class of problems, equivalent to MPECs, is that of mathematical programs with vanishing constraints.

The complementarity constraints \eqref{eq:mpec_comp} complicate the theoretical and computational aspects of the optimization problem~\eqref{eq:mpec}. 
These constraints can be replaced by a set of inequality constraints, leading to a standard nonlinear program (NLP):
\begin{subequations}\label{eq:mpec_nlp}
	\begin{align}
		\underset{x\in \R^{n}}{\mathrm{min}} \;  \quad &f(x)\\
		\textnormal{s.t.} \quad 
		%	&g(x)=0, \label{eq:mpec_eq}\\
		&c(x)\geq0, \label{eq:mpec_nlp_ineq}\\
		&x_{1,i} \geq 0,\ x_{2,i} \geq 0,\ x_{1,i} x_{2,i} \leq 0,\ \forall i \in \{1,\ldots,m\} \label{eq:mpec_nlp_comp}.
	\end{align}
\end{subequations}
However, this NLP is numerically and theoretically degenerate.
Standard constraint qualifications such as the Mangasarian-Fromovitz constraint qualification (MFCQ) fail at all feasible points~\cite{Jane2005,Scheel2000}, making standard NLP methods inefficient~\cite{Kim2020,Nurkanovic2024b} due to unbounded Lagrange multiplier sets~\cite{Gauvin1977}.
Moreover, standard KKT conditions may not apply to~\eqref{eq:mpec_nlp}, complicating stationary point definition and certification~\cite{Scheel2000}.
Various reformulations of~\eqref{eq:mpec_comp} exist, including smooth and nonsmooth C-functions~\cite{Facchinei2003,Leyffer2006a}, but none resolve the inherent degeneracy of NLP reformulations of \MPECs.

A necessary optimality condition is the absence of feasible first-order descent directions.
In the {\MPEC} literature, such points are called B-stationary points.
For regular NLPs, B-stationary points simply satisfy KKT conditions, but directly solving KKT conditions for~\eqref{eq:mpec_nlp} may be difficult or impossible due to degeneracy.
Weaker stationarity notions exist, but most allow first-order descent directions, making them too weak as practical stopping criteria~\cite{Leyffer2007}.
%This is further discussed in Section~\ref{sec:preliminaries}.

To cope with the degeneracy of NLP reformulations of {\MPECs} such as~\eqref{eq:mpec_nlp}, the optimization community has also developed \MPEC-tailored algorithms in the attempt to characterize and compute local minimizers of~\eqref{eq:mpec}, cf.~\cite{Hoheisel2013,Kanzow2015,Kim2020,Luo1996,Nurkanovic2024b} for surveys.
We briefly survey some of these methods.
They can be divided into regularization and active set/combinatorial methods.

Regularization methods replace the difficult constraints ~\eqref{eq:mpec_nlp_comp} with a \sloppy smoothed or relaxed version of them. 
For example, in Scholtes' global relaxation~\cite{Scholtes2001} $x_{1,i} x_{2,i} \leq 0$ is replaced by $x_{1,i} x_{2,i} \leq \tau$ for all $i$ and $\tau >0$. 
This yields a regular NLP that usually satisfies MFCQ and can be solved by standard NLP algorithms.  
In practice, a sequence of NLPs with decreasing $\tau$ is solved to approximate the {\MPEC} solution.  
Many refinements of this idea exist~\cite{DeMiguel2005,Kadrani2009,Kanzow2013,Lin2003,Scholtes2001,Steffensen2010}. 
Penalty methods~\cite{Anitescu2005,Anitescu2005a,Ralph2004} handle $x_{1,i} x_{2,i} \leq 0$ by adding a penalty term to the objective, e.g.\ $\rho \sum_{i=1}^m x_{1,i} x_{2,i}$.  
A sequence of NLPs with increasing $\rho > 0$ is then solved to approximate a solution.  
Under suitable conditions, even finite $\rho$ may yield an {\MPEC} stationary point via a standard NLP~\cite{Anitescu2005,Ralph2004}.  
However, unless restrictive assumptions hold, all these methods may converge to points that are not B-stationary~\cite{Hoheisel2013,Kanzow2015}.

%Interestingly, under appropriate assumptions, there is a one-to-one correspondence between the iterates of some relaxation and penalty methods~\cite{Leyffer2006}.

Active set methods take a different approach to obtaining a sequence of regular NLPs, whose solutions eventually solve the {\MPEC}. 
The main idea is to fix the active set complementarity constraints, i.e., replace \eqref{eq:mpec_comp} by $x_{1,i} = 0, x_{2,i} \geq 0, \forall i \in \I_1$ and $x_{1,j} \geq 0, x_{2,j} = 0, \forall j \in \I_2$, where $\I_1$ and $\I_2$ are a partition of $\{1,\ldots,m\}$.
The resulting NLPs -- called branch NLPs (BNLPs) -- are now regular problems that can be solved with classical NLP methods.
Given the correct partition, the solution of such an NLP also solves the {\MPEC}~\eqref{eq:mpec}. 
However, finding such a partition and certifying optimality have combinatorial complexity~\cite{Scheel2000}.
Most active set methods assume that {\MPEC}-LICQ holds (cf. Def.~\ref{def:mpec_cqs}) and search for  partitions by looking at the Lagrange multipliers of the active complementarity constraints to either find a new active set or to certify B-stationarity.
Such methods are reported in~\cite{Fukushima2002,Giallombardo2008,Izmailov2008,Jiang1999,Lin2006,Liu2006,Luo1996,Scholtes1999}.
However, if {\MPEC}-LICQ does not hold, only weaker concepts of stationarity -- possibly allowing descent directions -- can be used as stopping criteria.

To avoid the restrictive assumption of {\MPEC}-LICQ, B-stationarity can be verified by solving a linear program with equilibrium constraints ({\LPCC})~\cite{Luo1996}, obtained by linearizing the functions $f$ and $c$ in \eqref{eq:mpec}, but retaining the complementarity constraints.
If the current point is not B-stationary then the {\LPCC} provides a descent direction and a corresponding active set guess.
Such an approach is taken in~\cite{Guo2022,Kirches2022,Leyffer2007}.
Leyffer and Munson~\cite{Leyffer2007} propose a sequential {\LPCC} QP method.  
Since solving an {\LPCC} can be expensive, Kirches et al.~\cite{Kirches2022} study simpler LPECs with only complementarity and box constraints, and suggest handling of other constraints via an augmented Lagrangian.  
Guo and Deng~\cite{Guo2022} analyze this approach and prove convergence to M-stationary points (Def.~\ref{def:mpec_stationarity}).  
However, M-stationarity may still admit descent directions, cf.~\cite[Sec. 2.1]{Leyffer2007}.
A main drawback of these methods is their implementation complexity, since all components of a robust solver must be provided.  

%, e.g. promotion of global convergence, treatment of indefinite Hessians, second-order corrections, and efficient subproblem solvers.

Finally, we mention the class of hybrid methods~\cite{Kazi2024,Lin2006}, which combine active set and regularization-based methods. 
The method presented here, which we call MPEC optimizer ({\MPECopt}), can be classified as such a method and consists of two phases. 
Phase I finds a feasible BNLP by solving relaxed {\MPECs} for solution approximations, then {\LPCCs} to identify valid active sets.
Phase II starts from a feasible BNLP and solves {\LPCCs} at each BNLP's stationary point until reaching B-stationarity.
{\LPCC} solutions either certify B-stationarity or provide active set guesses for BNLPs with better objective values.
Notably, the {\LPCCs} steps never update the iterates, but only provide  better active sets or certify B-stationarity.

We highlight key differences from other hybrid approaches.
Unlike Lin and Fukushima~\cite{Lin2006}, our method solves {\LPCCs} for active set guesses rather than using multipliers, eliminating {\MPEC}-LICQ requirements and allowing convergence to B-stationary points that are not S-stationary.
Unlike Kazi et al.~\cite{Kazi2024}, we never apply {\LPCC} steps directly, ensuring all Phase II iterates remain feasible.

It is important to note that {\LPCCs} are nonconvex combinatorial problems whose efficient solution is crucial for our method's effectiveness.
{\LPCCs} with only bound and complementarity constraints have linear complexity~\cite[Propositions 2.3 and 2.4]{Kirches2022}.
As special cases of {\MPECs}, {\LPCCs} can be solved with regularization and penalty methods.
Tailored {\LPCC} methods include simplex extensions~\cite{Fang2010}, mixed-integer linear programming (MILP) approaches with tailored cuts~\cite{Hu2012,JaraMoroni2020}, and specialized penalty methods~\cite{Hall2024,JaraMoroni2018}; see~\cite{Judice2012} for a survey.
However, regularization and penalty methods may converge to points with trivial descent directions and fail to verify B-stationarity, making combinatorial methods preferable.
Importantly, we show that, in our proposed method, {\LPCCs} need only be solved to optimality when certifying a B-stationary point at the last iteration of Phase II.
Otherwise, in Phase I, it is sufficient to find a feasible point of the {\LPCC}, and in Phase II, a feasible descent direction.
Our numerical experiments show that this is often equivalent to solving a single LP.

%, and the BNLPs are often easier to solve than relaxed {\MPEC} subproblems.

%In summary, having an NLP and {\LPCC} solver at hand, the proposed method is easy to implement and finds B-stationary points in a finite number of steps without the usual restrictive assumptions.
%{\LPCCs} have due to the complementarity constraints also a combinatorial difficulty.
%Fortunately, as we demonstrated in numerous experiments, {\LPCCs} can be solved efficiently in practice either with regularization-based or mixed-integer linear programming (MILP) methods~\cite{Gurobi,Huangfu2018}. 
%Furthermore, we discuss the theoretical reasons why this is the case.
%We mention that also dedicated {\LPCC} methods exist~\cite{Hu2012,Fang2010,JaraMoroni2018,JaraMoroni2020}, which could be used within the proposed algorithm.
%All algorithmic details are discussed in Section~\ref{sec:algorithm}.

\textit{Contributions:}
This paper presents a new method, called MPEC optimizer ({\MPECopt}), that computes B-stationary points by solving a finite number of {\LPCCs} and NLPs in two phases.
We prove that, by finding a feasible point of an LPEC constructed at a point sufficiently close to the feasible set of the MPEC, the active set obtained from the LPEC corresponds to a BNLP that is guaranteed to be feasible.
For our implementation, this result guarantees that Phase I, which combines well-established regularization methods and an LPEC solver, finds feasible points of a BNLP (in fact, already stationary points of the \MPEC).
Phase II either certifies B-stationarity of this point or iterates over a sequence of feasible BNLPs while strictly decreasing the objective until a B-stationary point is found.
We provide a detailed convergence analysis and prove that, under the MPEC-MFCQ, Phase II is guaranteed to find a B-stationary point in a finite number of \LPCC\ and BNLP solves.

%Furthermore, we discuss in detail two different {\LPCC} formulations, their relation, and theoretical properties.
We analyze two LPEC formulations and their relationship.
We show that, unless we need to certify B-stationarity, it is sufficient to find a feasible (improving) point of the {\LPCC} to ensure progress, which often reduces the computational effort significantly.
In practice, solving a single linear program is often sufficient.
We also discuss in detail why {\LPCCs} can be solved efficiently via their MILP reformulation.

We compare the proposed method with several relaxation-based methods and a mixed-integer nonlinear programming (MINLP) reformulation on two benchmarks: the {\MacMPEC} set and a synthetic large-scale nonlinear {\MPEC} benchmark.
Higher success rates and generally faster computation times are reported.
Interestingly, solving {\LPCCs} turns out not to be a computational bottleneck in practice, even for large problem instances.
We highlight that most other methods with an available implementation -- including those we compare to -- do not, in general, provide a certificate of B-stationarity for generic \MPECs\ and they may converge to a point with a first-order descent direction.
Additionally, an open-source implementation of the method is provided at \url{github.com/nosnoc/mpecopt}.

\textit{Notation:} 
Subscripts denote vector/matrix components ($x_i$), superscripts denote iterates ($x^k$).
Functions evaluated at particular iterates are denoted as $f^k := f(x^k)$.
The concatenation of two vectors $a \in R^n$ and $b \in \R^m$ is written compactly as $(a,b) = \begin{bmatrix}
	a^\top & b^\top
\end{bmatrix}^\top$. 
The notation extends analogously to more than two vectors.
For scalars $a,b\in\R$, $\max(a,b)$ returns the larger value.
For vectors $a, b \in \R^n$: $\max(a)$ returns the largest component, $c=\max(a,b)$ gives $c_i = \max(a_i,b_i)$.
%Given two scalars $a,b\in\R$, $\max(a,b)$ returns the larger of them.
%Let $a, b \in \R^n$, then $\max(a)$ returns the largest component of $a$, and $c=\max(a,b)$ is a vector in $\R^n$ with the components $c_i = \max(a_i,b_i), i =1,\ldots,n$. 
The $\min$ function is defined in all cases accordingly.
Let $x\in\R^{n}$, then $\mathrm{diag}(x)\in\R^{n\times n}$ returns a diagonal matrix with $x_i$ as its diagonal elements.

\textit{Outline:}
Section \ref{sec:preliminaries} recalls some preliminaries.
Section presents \ref{sec:algorithm} {\MPECopt} and discusses the formulation and solution of {\LPCCs}.
Section \ref{sec:convergence_theory} gives a detailed convergence analysis of the new method.
Section \ref{sec:numerical_results} {provides all implementation details and numerical results on two benchmark test sets}. 
Section \ref{sec:conclusion} provides concluding remarks.
\section{Preliminaries}\label{sec:preliminaries}
%{
%	This section reviews first-order optimality conditions, stationarity concepts, which are needed for formulating the proposed algorithm and analyzing its convergence.
%}
%\subsection{Optimality conditions for {\MPECs}}
We start by defining some useful index sets for the complementarity constraints of the {\MPEC}~\eqref{eq:mpec}.
Denote the feasible set of the {\MPEC}~\eqref{eq:mpec} by $\Omega = \{ x\in \R^n \mid c(x)\geq0, 0 \leq x_1 \perp x_2 \geq 0\}$.
The following index sets, which depend on a feasible point $x$, define a partition of $\{1,\ldots,m\}$:
%\allowdisplaybreaks
\begin{align*}
	\I_{0+}(x) &= \{i \in \{1,\ldots,m\}\mid x_{1,i}=0, x_{2,i}>0\},\\
	\I_{+0}(x) &=	\{i \in \{1,\ldots,m\} \mid x_{1,i}>0, x_{2,i}=0\},\\
	\I_{00}(x) &= \{i \in \{1,\ldots,m\}\mid x_{1,i}=0, x_{2,i}=0\}.
\end{align*}
Most of the computational and theoretical difficulties arise if the so-called biactive set $\I_{00}(x)$ is nonempty.
The active set of the standard inequality constraints \eqref{eq:mpec_ineq} is defined as 
${\A(x) = \{ i \in \{1,\ldots,n_c\} \mid c_i(x)=0\}}$.

Next, we state first-order necessary optimality conditions.
For this we need the notion of the tangent cone at $x\in\Omega$ to the set $\Omega$, which is defined as $\mathcal{T}_{\Omega}(x)  
= \{ d\in \R^n \mid \exists \{x^k\} \subset \Omega, \{t^k\}\subset \R_{\geq0}: \lim\limits_{k\to \infty} t^k = 0, \lim\limits_{k\to \infty} x^k = x, \lim\limits_{k\to \infty} \frac{x^k-x}{t^k} = d\}$. 
%Let $d = (d_0,d_1,d_2) \in \R^{p}\times \R^{m}\times \R^{m}$. 
%The first-order necessary optimality conditions for \eqref{eq:mpec} are given by the following theorem~\cite[Corollary 3.3.1]{Luo1996}. 
\begin{mytheorem}(\cite[Corollary 3.3.1]{Luo1996})\label{th:b_stationarity}
	If $x^*\in\Omega$ is a {local} minimizer of \eqref{eq:mpec}, then it holds that
	\begin{align}\label{eq:geometric_b_stationarity}
		\nabla f(x^*)^\top d \geq 0\; \textrm{ for all } d \in \mathcal{T}_{\Omega}(x^*),
	\end{align}
	or equivalently, $d = 0$ is {a global optimizer} of the following optimization problem:
	\begin{align}\label{eq:b_stationariry}
		\min_{d \in \R^{n}} \quad &   \nabla f(x^*)^\top d\quad 	\textnormal{s.t.} \quad   d \in \mathcal{T}_{\Omega}(x^*).
	\end{align}
\end{mytheorem}
If a point $x^*$ satisfies the condition above, it is said that geometric Bouligand stationarity holds, or for short, $x^*$ is geometric B-stationary~\cite{Flegel2005a,Luo1996,Scheel2000}.

Note that Theorem \ref{th:b_stationarity} is purely geometric and thus difficult to verify.
In standard optimization theory, to obtain an algebraic characterization of a stationary point, the tangent cone $\T_{\Omega}(x)$ is replaced by the linearized feasible cone~\cite[Def. 12.3]{Nocedal2006}.
However, the latter is always a convex polyhedral cone, and if $\mathcal{I}_{00}(x)$ is nonempty the cone $\mathcal{T}_{\Omega}(x)$ is always nonconvex~\cite{Luo1996}, making this transition impossible.
A more suitable notion is the so-called {\MPEC}-linearized feasible cone~\cite{Flegel2005a,Pang1999,Scheel2000}, which for a feasible point $x\in\Omega$ reads as
\begin{align}\label{eq:mpec_lin_cone}
	\begin{split}
		\mathcal{F}_{\Omega}(x) = \{ d \in \R^n \mid \
%		& g(x)+\nabla g(x)^\top d =0,\\
		& c_i(x)+ \nabla c_{i} (x)^\top d \geq 0,\ \forall  i \in \A(x),\\
		&  x_{1,i} + d_{1,i} = 0,\ \forall  i \in \mathcal{I}_{0+}(x),\\
		&  x_{2,i} + d_{2,i} = 0,\ \forall  i \in \mathcal{I}_{+0}(x),\\
		&0 \leq  x_{1,i} + d_{1,i}  \perp  x_{2,i}+ d_{2,i}  \geq0, \forall  i \in \mathcal{I}_{00}(x)
		\}.
	\end{split}
\end{align}
Similarly to $x = (x_0,x_1,x_2)$, the vector $d$ is partitioned as $d = (d_0,d_1,d_2) \in \mathbb{R}^{p} \times \mathbb{R}^m \times \mathbb{R}^m$.
Observe that at a feasible point $x\in\Omega$ all constraint residuals above are zero but they are kept for convenience.
In contrast to the usual linearized feasible cone, here the combinatorial structure is kept for the degenerate index set $\I_{00}(x)$ and $\mathcal{F}_{\Omega}(x)$ is nonconvex if $\mathcal{I}_{00}$ is nonempty.
%It is a union of a finite number of convex polyhedral cones.

In general, it holds that $\T_{\Omega}(x) \subseteq \mathcal{F}_{\Omega}(x)$~\cite{Flegel2005a}. 
{If the so-called {\MPEC}-Abadie constraint qualification (MPEC-ACQ) is satisfied, i.e., $\T_{\Omega}(x) = \mathcal{F}_{\Omega}(x)$~\cite{Flegel2005a}, we can use $\mathcal{F}_{\Omega}(x)$ in Theorem \ref{th:b_stationarity}.} 
In particular, one can define an {\LPCC} to characterize an {algebraic B-stationary point (sometimes also called linearized B-stationarity~\cite{Flegel2005a}), or B-stationary point for short.}

\begin{mydefinition}[B-stationarity]\label{def:b_stationarity}
A point $x^*$ is called (algebraic) B-stationary if $d=0$ solves the following {\LPCC}:
\begin{subequations}\label{eq:lpec_reduced_theory}
	\begin{align}
		\underset{d\in \R^{n}}{\mathrm{min}} \;  \quad & \nabla f(x^*)^\top d \\
		\textnormal{s.t.} \quad 
		& c_i(x^*)+ \nabla c_i(x^*)^\top d\geq0, & \forall i \in \A(x^*) \\
		&  x_{1,i}^* + d_{1,i} = 0,\ 
%		\color{red} x_{2,i}^* + d_{2,i} \geq 0 ,
		&\forall i \in  \I_{0+}(x^*), \label{eq:lpec_reduced_theory_branch1}\\
		& 
%		\color{red} x_{1,i}^* + d_{1,i} \geq 0 ,\  
		x_{2,i}^* + d_{2,i} = 0, &\forall i \in \I_{+0}(x^*), \label{eq:lpec_reduced_theory_branch2}\\
		&0 \leq  x_{1,i}^* + d_{1,i}  \perp x_{2,i}^* + d_{2,i}  \geq0, &\forall i \in \I_{00}(x^*).
	\end{align}
\end{subequations}
\end{mydefinition}
To verify that $x^*$ is B-stationary, one has to verify that $d =0$ is a minimizer of an {\LPCC}.
Note that if $| \I_{00}(x^*)| = 0$, then \eqref{eq:lpec_reduced_theory} is a linear program (LP).

Next, we review multiplier-based stationary concepts for \MPECs\ and some auxiliary NLPs associated with the \MPEC~\eqref{eq:mpec}. 
One of these problems will be solved by the algorithm introduced below, while others serve to define additional concepts relevant to {\MPECs}. 
To this end, we introduce further index sets.
Let $x^*\in \Omega$ and $\D_1(x^*), \D_2(x^*)$ be a partition of the biactive set $\I_{00}(x^*)$, with $\mathcal{D}_1(x^*) \cup \mathcal{D}_2(x^*)=\I_{00}(x^*)$ and $\mathcal{D}_1(x^*) \cap \mathcal{D}_2(x^*)=\emptyset$.
\begin{align}\label{eq:active_set_partition}
&\I_1(x^*) = \I_{0+}(x^*)\cup \mathcal{D}_1(x^*), \quad
\I_2(x^*) = \I_{+0}(x^*)\cup \mathcal{D}_2(x^*).
\end{align}
{Let us denote the set of all possible partitions by $\mathcal{P}(x^*) = \{ (\I_1(x^*), \I_2(x^*))\}$.
The number of possible partitions is  $|\mathcal{P}(x^*)|= 2^{|\I_{00}(x^*)|}$. 
}
For every partition, we can define a branch NLP, denoted by {$\mathrm{BNLP}{(\I_1(x^*),\I_2(x^*))}$}, as follows:
	\begin{subequations} \label{eq:bnlp}
	\begin{align}
		\min_{x \in \R^{n}} \quad &   f(x) \label{eq:bnlp_objective}\\
		&c(x)\geq 0, \label{eq:bnlp_ineq}\\
		&x_{1,i} = 0,\ x_{2,i} \geq 0,\; &\forall  i \in 	\I_{1}(x^*), \label{eq:bnlp_branch_h}\\
		&x_{1,i} \geq 0,\  x_{2,i}= 0,\; &\forall  i \in 	\I_{2}(x^*). \label{eq:bnlp_branch_g}
	\end{align}
\end{subequations}
%The difficult complementarity constraints are replaced by simple equality and non-negativity constraints in~\eqref{eq:bnlp_branch_h} and \eqref{eq:bnlp_branch_g}. 
%\red{If the standard constraints $c(x) \geq0$ do not lead to a violation of constraint qualifications, \eqref{eq:bnlp} are regular NLPs for which the standard optimization theory is applicable and existing numerical algorithms perform well.}
Additionally, the {\MPEC} literature frequently uses the so-called tight NLP ({TNLP$(x^*)$}):
	\begin{subequations} \label{eq:tnlp}
	\begin{align}
		\min_{x \in \R^{n}} \quad &   f(x) \label{eq:tnlp_objective}\\
		\textnormal{s.t.} \quad  
%		&g(x) = 0 \label{eq:tnlp_eq},\\
%		&h(x)\geq 0, \label{eq:rnlp_ineq}\\
		&c(x)\geq 0, \label{eq:rnlp_ineq}\\
		&x_{1,i} = 0,\ x_{2,i} \geq 0,\; &\forall i \in 	\mathcal{I}_{0+}(x^*), \label{eq:tnlp_branch_h}\\
		&x_{1,i} \geq 0,\ x_{2,i} = 0,\; &\forall  i \in 	\mathcal{I}_{+0}(x^*), \label{eq:tnlp_branch_g}\\
		&x_{1,i}= 0,\ x_{2,i}=0,\; &\forall  i \in 	\mathcal{I}_{00}(x^*) \label{eq:tnlp_biactive},
	\end{align}
\end{subequations}
{The so-called relaxed NLP (RNLP$(x^*)$) is defined by replacing \eqref{eq:tnlp_biactive} by
\[x_{1,i} \geq 0,\ x_{2,i} \geq 0,\; \forall  i \in 	\mathcal{I}_{00}(x^*).\]
We denote the feasible sets of the BNLP$(\mathcal{I}_1(x^*),\mathcal{I}_2(x^*))$ and TNLP$(x^*)$ by $\Omega_{\mathrm{BNLP}{(\mathcal{I}_1,\mathcal{I}_2)}}$ and $\Omega_{\mathrm{TNLP}}$, respectively, and if clear from the context, we omit the dependency on $x^*$.
}

It can be seen that the following holds locally around a point $x^* \in \Omega$~\cite{Scheel2000}:
\begin{align}\label{eq:mpec_feasible_sets}
	\Omega_{\mathrm{TNLP}}  = \bigcap_{(\mathcal{I}_1,\mathcal{I}_2){\in \mathcal{P}}} \Omega_{\mathrm{BNLP}{(\mathcal{I}_1,\mathcal{I}_2)}} \subset \Omega
	= \bigcup_{(\mathcal{I}_1,\mathcal{I}_2) {\in \mathcal{P}}} \Omega_{\mathrm{BNLP}{(\mathcal{I}_1,\mathcal{I}_2)}}. 
\end{align}
From \eqref{eq:mpec_feasible_sets} for a feasible point $x^* \in \Omega$ the following can be concluded \cite{Scheel2000}. 
If $x^*$ is a local minimizer of the {\MPEC}~\eqref{eq:mpec} then it is a local minimizer of the TNLP$(x^*)$. 
The point $x^*$ is a local minimizer of the {\MPEC} if and only if it is a local minimizer of every $\mathrm{BNLP}_{(\I_1(x^*),\I_2(x^*))}$. 
%The last assertion once again highlights the combinatorial nature of {\MPECs} as -- unless some strong assumptions specified below do not hold -- an exponential number of branch NLPs must be checked to make conclusions about stationarity.
The TNLP is used to define some {\MPEC}-specific concepts, and we recall those relevant for this paper.
Note that if $\mathcal{I}_{00}(x^*) = \emptyset$, there is only one BNLP, which is  equal to the TNLP.
%Moreover, there are no complementarity conditions in \eqref{eq:lpec_reduced_theory}, the {\LPCC} reduces to a linear program (LP), and for $d=0$ the Karush-Kuhn-Tucker (KKT) conditions of this LP coincide with the KKT conditions of the TNLP.

%, such as tailored constraint qualifications (CQs), an {\MPEC}-Lagrangian, and multiplier-based stationarity concepts for~\eqref{eq:mpec}.
%First, we define some MPEC CQs.
\begin{mydefinition}\label{def:mpec_cqs}
The {\MPEC} \eqref{eq:mpec} is said to satisfy the {\MPEC}-LICQ ({\MPEC}-MFCQ) at a feasible point $x^*\in \Omega$ if the corresponding $\mathrm{TNLP}$ \eqref{eq:tnlp} satisfies the LICQ (MFCQ) at the same point $x^*$.
\end{mydefinition}
Next, we define the {\MPEC} Lagrangian. 
This is simply the \textit{standard} Lagrangian for the BNLP/TNLP and reads as:
\begin{align}
	\mathcal{L}^{\mathrm{\MPEC}}(x,\mu,\nu,\xi) = f(x)  - \mu^\top c(x)- \nu^\top x_1- \xi^\top x_2,
\end{align}
with the {\MPEC} Lagrange multipliers $\mu \in \R^{n_c}$, $\nu \in \R^{m}$ and $\xi \in \R^{m}$.
%In contrast to the standard Lagrangian of the \MPEC \eqref{eq:mpec_nlp} it does not have the bilinear terms $x_{1,i}x_{2,i}\leq0$ and their multipliers.
{
By applying the KKT conditions to the TNLP~\eqref{eq:tnlp}, we obtain the concept of W-stationarity, whereas applying them to the RNLP yields S-stationarity~\cite{Scheel2000}. 
Further restrictions on the multipliers 
$ \nu^*_i, \xi^*_i, \ i \in \I_{00}(x^*)$, lead to M-, C-, and A-stationarity, which are not directly related to the KKT conditions of an auxiliary NLP, see e.g. ~\cite{Kim2020,Outrata1999,Scheel2000}.
These concepts are summarized in the next definition.
}

\begin{mydefinition}[Stationarity concepts for {\MPECs}]\label{def:mpec_stationarity}
	For a feasible point $x^* \in \Omega$, we distinguish the following stationarity concepts.
	\begin{itemize}
		\item Weak stationarity (W-stationarity) \cite{Scheel2000}: 
		A point $x^*$ is called W-stationary if the corresponding {TNLP$(x^*)$}~\eqref{eq:tnlp} admits the satisfaction of the KKT conditions, i.e., there exist Lagrange multipliers $\mu^*,\nu^*$ and $\xi^*$ such that:
		\begin{align*}
			&\nabla_x \mathcal{L}^{\mathrm{\MPEC}}(x^*,\mu^*,\nu^*,\xi^*) = 0,\\
			&0 \leq \mu^* \perp c(x^*) \geq 0,\\
			& x^{*}_{1,i} \geq0, \nu^*_i = 0, \; &\forall i \in 	\I_{+0}(x^*),\\
			& x^{*}_{2,i} \geq0, \xi^*_i = 0,\; &\forall i \in 	\I_{0+}(x^*),\\
			& x^{*}_{1,i} =0,\; \nu_i^* \in \R, \; &\forall i \in 	\I_{0+}(x^*)\cup\I_{00}(x^*),\\
			&x^{*}_{2,i} =0,\; \xi_i^* \in \R, \; &\forall i \in 	\I_{+0}(x^*)\cup\I_{00}(x^*).
		\end{align*}
		\item Strong stationarity (S-stationarity) \cite{Scheel2000}:  
		A point $x^*$ is called S-stationary if it is weakly stationary and $\nu^*_i \geq 0, \xi^*_i \geq0$ for all $i \in \I_{00}(x^*)$. 
		\item Clarke stationarity (C-stationarity) \cite{Scheel2000}:  
		A point $x^*$ is called C-stationary if it is weakly stationary and $\nu^*_i\xi^*_i \geq0$ for all $i \in \I_{00}(x^*)$.
		\item Mordukhovich stationarity (M-stationarity) \cite{Outrata1999}:  
		A point $x^*$ is called M-stationary if it is weakly stationary and if either $\nu^*_i >0$ and $\xi^*_i >0$ or $\nu^*_i\xi^*_i =0$ for all $i \in \I_{00}(x^*)$.
		\item Abadie stationarity (A-stationarity) \cite{Flegel2005a}: 
		A point $x^*$ is called A-stationary if it is weakly stationary and $\nu^*_i \geq 0$ or $\xi^*_i \geq0$ for all $i \in \I_{00}(x^*)$.
	\end{itemize}
\end{mydefinition}
Note that, if $\I_{00}(x^*) = \emptyset$, then all stationarity concepts are identical and collapse to the concept of S-stationarity.
From the definitions it follows that S-stationarity is the most restrictive since it has the tightest condition in the multipliers $\nu^*_i, \xi^*_i, i \in \I_{00}(x^*)$.
Weaker stationarity concepts, ordered by strength, are M, C, A, and W, with the stronger ones implying the weaker ones.
The relationship between multiplier-based stationarity (Def.~\ref{def:mpec_stationarity}) and B-stationarity (Def.~\ref{def:b_stationarity}) is established in the following theorem.
\begin{mytheorem}(\cite[Theorem 4]{Scheel2000})\label{th:b_and_s_stat}
	If $x^*$ is an S-stationary point of the \sloppy {\MPEC}~\eqref{eq:mpec}, then it is also B-stationary.
	If, in addition, the {\MPEC}-LICQ holds, then every B-stationary point is S-stationary.
\end{mytheorem}

\begin{myremark}
In summary, without assuming {\MPEC}-LICQ, B-stationarity cannot be verified using multiplier-based concepts.
This assumption is essential, as the next weaker condition, {\MPEC}-MFCQ, does not guarantee S-stationarity (cf.~\cite[Example 3]{Scheel2000}).  
Weaker termination criteria than B-stationarity may be considered unsatisfactory, since even M-stationarity allows first-order descent directions~\cite{Leyffer2007}.  
Even under {\MPEC}-LICQ, regularization~\cite{Ralph2004,Hoheisel2013,Kanzow2015} and some active set methods~\cite[Section 3]{Kirches2022} may converge to points weaker than S-stationary.  
Therefore, we propose a method using the {\LPCC}~\eqref{eq:lpec_reduced_theory}, to certify B-stationarity.  
\end{myremark}
 
\section{\MPECopt: A globally convergent method for computing B-stationary points of {\MPECs}}\label{sec:algorithm}

This section presents the MPEC Optimizer ({\MPECopt}), an easy-to-implement active set method that solves a finite sequence of {\LPCCs} {or LPs,} and BNLPs. 
It returns either a B-stationary point {or a stationary point of the problem of minimizing the constraint infeasibility.
The latter is a common outcome of nonlinear programming algorithms~\cite{Fletcher2002,Waechter2006}.}

Before providing all the details, we summarize the main ideas behind the new algorithm in Section~\ref{sec:algorithm_summary}. 
Thereafter, Section~\ref{sec:algorithm_lpecs} discusses different formulations of {\LPCCs}, their properties, and how {\LPCCs} can be solved numerically. 
{\MPECopt} consists of two phases. 
Details of Phase I and Phase II are discussed in Sections~\ref{sec:algorithm_phase_i} and ~\ref{sec:algorithm_phase_ii}, respectively. 
%The first phase, described in Section~\ref{sec:algorithm_phase_i}, identifies a first feasible BNLP or {returns a stationary point of the problem of minimizing the constraint infeasibility}. 
%The second phase, described in Section~\ref{sec:algorithm_phase_ii}, solves a sequence of BNLPs and {\LPCCs} and returns a B-stationary point.

\subsection{Summary of the approach}\label{sec:algorithm_summary}
\begin{algorithm}[h]\label{algo:mpec_opt_summary}
	\vspace{-0.1cm}
	\caption{\MPECopt: Computing a B-stationary point of {\MPECs}}
	\label{alg:mpec_opt_summary}
	\begin{algorithmic}[1]
		\State Call \texttt{Phase I} to find a feasible point $x^0$ of the \MPEC~\eqref{eq:mpec}; set $k = 0$
		\While{$d^{k,l} \neq 0$} \Comment{\textcolor{gray}{major iterations}}
		\State $x^{k,0} = x^k$
		\For{$l=0,\ldots,$} \Comment{\textcolor{gray}{minor iterations}} 
		\State Solve $\textrm{LPEC}(x^k,\rho^{k,l})$ in \eqref{eq:lpec_full} to obtain $d^{k,l}$
		\If{$d^{k,l} = 0$} \textbf{terminate} \Comment{\textcolor{gray}{B-stat. point found}}
		\EndIf
		\State Solve $\mathrm{BNLP}(\I_1(x^k+d^{k,l}),\I_2(x^k+d^{k,l}))$ to obtain $x^{k,l}$
		\If {BNLP unbounded} \textbf{terminate}	\Comment{\textcolor{gray}{MPEC is unbounded}}	
		\ElsIf{$f(x^{k,l}) < f(x^{k})$} 
		\State Set $x^{k+1} \!= \! x^{k,l}$, $k \! = \! k\!+\!1$, increase trust-region radius $\rho^{k+1,0}$, and \textbf{break} for-loop
		\Else
		\State Reduce trust-region radius $\rho^{k,l+1}$
		\EndIf 
		\EndFor
		\EndWhile
	\end{algorithmic}
	\vspace{-0.1cm}
\end{algorithm}
Algorithm~\ref{alg:mpec_opt_summary} provides simplified pseudo-code for {\MPECopt}, which consists of two phases.
Phase I either finds an MPEC-feasible point $x^0 \in \Omega$ or a stationary point of the problem of minimizing the constraint infeasibility.
Phase II is summarized in lines 2--11 of Algorithm~\ref{alg:mpec_opt_summary}, while the details of our actual implementation are given in Section~\ref{sec:algorithm_phase_ii}.
The outcome of Phase II is that it finds a B-stationary point of the MPEC {or detects that the problem is unbounded.}

If a feasible point is found in Phase I, Phase II creates a finite sequence of iterates $x^k \in \Omega$ by solving BNLPs and {\LPCCs}.
If Phase II finds that any BNLP is unbounded, then the MPEC is unbounded, and the algorithm terminates.
The {\LPCCs} are similar to~\eqref{eq:lpec_reduced_theory}, but they are defined at the linearization point $x^k$ and include an additional trust-region constraint $\|d \|_{\infty} \leq \rho$, $\rho >0$. 
These {\LPCCs} are compactly denoted by $\mathrm{LPEC}(x,\rho)$; see Section~\ref{sec:algorithm_lpecs} for details.

At each major iteration $k$, an inner loop with counter $l$ solves $\mathrm{LPEC}(x^k,\rho^{k,l})$ for $d^{k,l}$ and, possibly, a corresponding $\mathrm{BNLP}(\I_1(x^k+d^{k,l}),\I_2(x^k+d^{k,l}))$, until either a feasible point $x^{k+1} = x^{k,l}$ with an improved objective value is found (lines 7--10), or the algorithm terminates due to $d^{k,l} = 0$ (line 6) or the {\MPEC} is found to be unbounded (line 8).

If $d^{k,l} = 0$, the algorithm terminates as a B-stationary point has been found.
Note that the optimal value of the {\LPCC} may yield $\nabla f(x^k)^\top d^{k,l} = 0$ with $d^{k,l} \neq 0$, for example when the solution lies on the boundary of the feasible set.
Since $d^{k,l} = 0$ is always feasible for the {\LPCC}~\eqref{eq:lpec_reduced_theory}, we consistently select $d^{k,l}=0$ in such cases and terminate; cf. Section~\ref{sec:algorithm_lpecs}.
Otherwise, when $d^{k,l} \neq 0$, this solution is used to determine the active sets $\I_1(x^k+d^{k,l})$ and $\I_2(x^k+d^{k,l})$ for the $\BNLP(\I_1(x^k+d^{k,l}),\I_2(x^k+d^{k,l}))$.
There is some ambiguity in the selection of active sets $\I_1$ and $\I_2$, which is discussed in detail in Section~\ref{sec:algorithm_lpecs}. 
In line 7, we compute a stationary point of this BNLP, which is denoted by $x^{k,l}$.
If $f(x^{k,l}) < f(x^{k})$, we accept the step, increase or reset the trust-region radius, and set the major iteration counter to $k+1$. 
If the condition in line 9 is not satisfied, the trust-region radius is reduced, and the \LPCC$(x^k,\rho^{k,l})$ is resolved in the next minor iteration, with $l = l+1$.
In Section~\ref{sec:convergence_theory}, we prove that both the inner and outer loops terminate finitely.

Note that the LPEC steps $x^k+d^{k,l}$ are never used to update iterates directly.  
This avoids cases where, under nonlinear constraints, $x^k+d^{k,l} \notin \Omega$ while $f(x^k+d^{k,l})$ increases, which would require sophisticated globalization for the {\LPCC} steps.  
Instead, we accept only BNLP stationary points $x^{k,l}$, which remain feasible for the {\MPEC}~\eqref{eq:mpec}.

\subsection{The {\LPCC} subproblems}\label{sec:algorithm_lpecs}
In Section \ref{sec:formulating_lpecs} we discuss and compare two different formulations for the {\LPCC} used in our method, and in Section \ref{sec:solving_lpecs} we discuss how they are efficiently solved.

\subsubsection{On different {\LPCC} formulations}\label{sec:formulating_lpecs}
{The first formulation in inspired by Definition~\ref{def:b_stationarity}.
Given a feasible point $x^k \in \Omega$, Phase II in Algorithm~\ref{alg:mpec_opt_summary}} requires solving the following $\mathrm{LPEC}(x^k,\rho)$:
\begin{subequations}\label{eq:lpec_reduced}
	\begin{align}
		\underset{d\in \R^{n}}{\mathrm{min}} \;  \quad & \nabla f(x^k)^\top d \\
		\textnormal{s.t.} \quad 
		& c_i(x^k)+ \nabla c_i(x^k)^\top d\geq0, &{\forall i \in \A(x^k)}\\\
		& x_{1,i}^k + d_{1,i} = 0,\ x_{2,i}^k + d_{2,i} \geq 0 ,
		&\forall i \in  \I_{0+}(x^k), \label{eq:lpec_reduced_branch1}\\
		& x_{1,i}^k + d_{1,i} \geq 0 ,\  x_{2,i}^k + d_{2,i} = 0, &\forall i \in \I_{+0}(x^k), \label{eq:lpec_reduced_branch2}\\
		&0 \leq  x_{1,i}^k + d_{1,i}  \perp x_{2,i}^k + d_{2,i}  \geq0, &\forall i \in \I_{00}(x^k), \label{eq:lpec_reduced_comp}\\
		&\| d \|_{\infty} \leq \rho. \label{eq:lpec_reduced_tr}
	\end{align}
\end{subequations}
Compared to~\eqref{eq:lpec_reduced_theory} in Def.~\ref{def:b_stationarity}, the {\LPCC} above has the following modifications.
In addition to \eqref{eq:lpec_reduced_theory_branch1} and \eqref{eq:lpec_reduced_theory_branch2}, here we include the constraint $ x_{2,i}^k + d_{2,i} \geq 0 $ and $x_{1,i}^k + d_{1,i} \geq 0$ in \eqref{eq:lpec_reduced_branch1} and \eqref{eq:lpec_reduced_branch2}, respectively. 
This ensures that for any solution $d$ of \eqref{eq:lpec_reduced} the point $x^k+d$ satisfies the complementarity constraints, which is necessary for a meaningful definition of the active sets in ~\eqref{eq:active_set_partition}.
Moreover, we add a trust region constraint for the step $d$ in \eqref{eq:lpec_reduced_tr}, where $\rho>0$ is the trust region radius. 
This constraint guarantees that the {\LPCC} is bounded.  
We use an $\ell_{\infty}$ norm for the step, yielding simple lower and upper bounds on $d$.

If $d=0$ is a global minimizer of \eqref{eq:lpec_reduced}, then none of the additionally introduced constraints is active, and moreover, if the MPEC-ACQ holds, then by Theorem~\ref{th:b_stationarity}, $x^k$ is a B-stationary point of {\MPEC}~\eqref{eq:mpec}.
On the other hand, if $x^k$ is not B-stationary (i.e., $d=0$ is not a minimizer of \eqref{eq:lpec_reduced}), and if $\A(x^k+d) = \A(x^k)$, then $\nabla f(x^k)^\top d < 0$ and $d$ is a first-order descent direction. 
In this case, we take the sets $\I_1(x^k + d)$ and $\I_2(x^k+d)$ as the new active set guess for the next BNLP to be solved in Phase II; more details on the actual selection of these index sets are provided below.

Next, since algebraic B-stationarity is stronger than geometric B-stationarity, we show that the assumption of {\MPEC}-ACQ (i.e., that $\T_{\Omega}(x)=\mathcal{F}_{\Omega}(x)$), cannot be relaxed if one wishes to verify {geometric} B-stationarity ({cf. Theorem \ref{th:b_stationarity}}) by solving the {\LPCC}~\eqref{eq:lpec_reduced}, which verifies algebraic B-stationarity (cf. Def.~\ref{def:b_stationarity}).
%Recall that the {\MPEC}-Abadie CQ is implied by the {\MPEC}-MFCQ~\cite{Flegel2005a}.

%\vspace{-0.6cm}
\begin{figure}[h]
	\centering
	\subfloat[The  MPEC~\eqref{eq:mpec_cq_fail}.]{\includegraphics[width=0.33\textwidth]{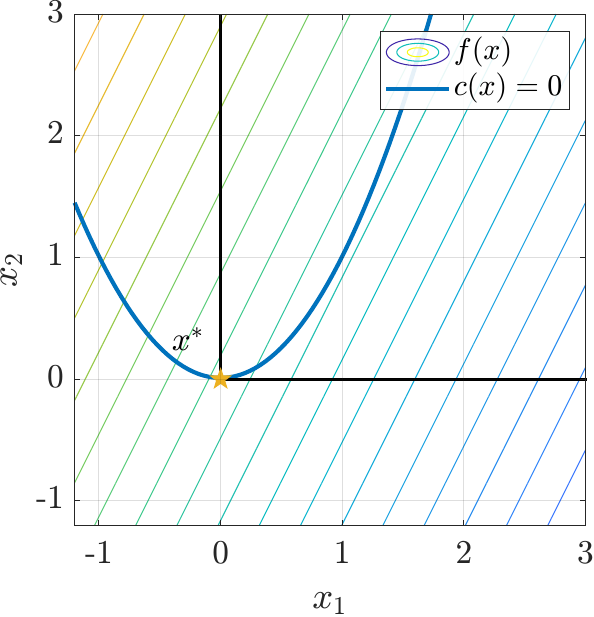}\label{fig:cq_mpec_1}}
	\subfloat[{\LPCC}~\eqref{eq:lpeq_cq_fail} $\rho = 1$.]{\includegraphics[width=0.33\textwidth]{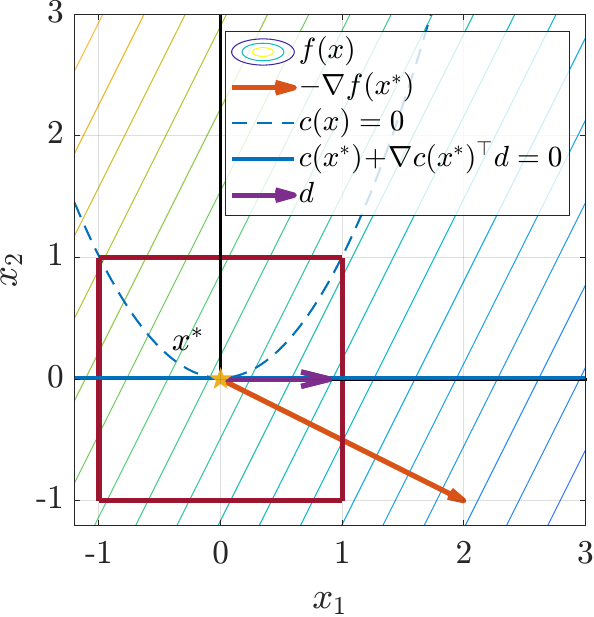}\label{fig:cq_lpec_1}}
	\subfloat[{\LPCC}~\eqref{eq:lpeq_cq_fail} with $\rho = 0.5$.]{\includegraphics[width=0.33\textwidth]{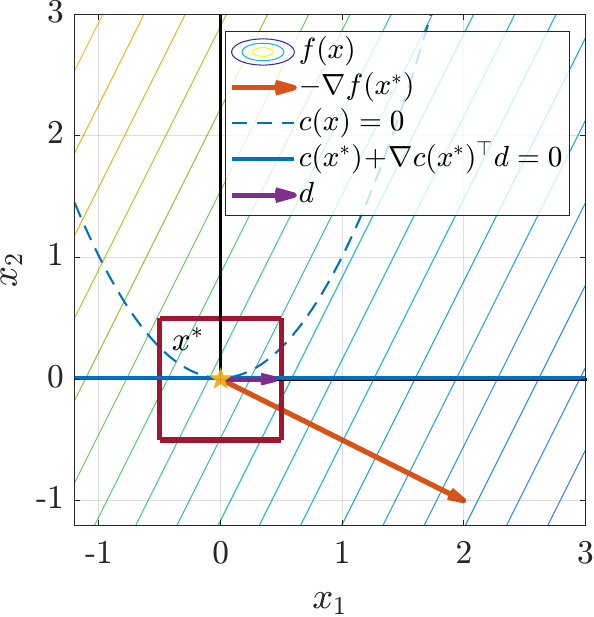}\label{fig:cq_lpec_2}}
	\caption{Illustration of Example~\ref{ex:mpec_cq_fail}. If a constraint qualification does not hold, the {\LPCC} cannot verify geometric B-stationarity.}
	\label{fig:mpec_cq_fail}
	\vspace{-0.5cm}
\end{figure}
\begin{myexample}\label{ex:mpec_cq_fail}
	Consider the following {\MPEC}:
	\begin{subequations}\label{eq:mpec_cq_fail}
		\begin{align*}
			\underset{x\in \R^{2}}{\mathrm{min}} \;  & -2x_1 + x_2\
			\textnormal{s.t.} \
			x_2-x_1^2\geq0,\  0\leq x_1 \perp x_2\geq0.
		\end{align*}
	\end{subequations}
	The point $x^* = (0,0)$ is the global minimizer of this problem and thus {geometric} B-stationary.
	It can be seen that the tangent cone of this problem is $\T_{\Omega}(x^*)  = \{d \in \R^2 \mid d_1 =0, d_2\geq0\}$.
	The corresponding {\LPCC} reads as: 
	\begin{subequations}\label{eq:lpeq_cq_fail}
		\begin{align}
			\underset{d\in \R^{2}}{\mathrm{min}} \;  \quad & -2d_1 + d_2\\
			\textnormal{s.t.} \quad 
			&d_2\geq0\\
			&0\leq d_1 \perp d_2\geq0\\
			& \|d \|_\infty \leq \rho.
		\end{align}
	\end{subequations}
	Clearly, $\T_{\Omega}(x^*) \subset \mathcal{F}_{\Omega}(x^*)$, which means that the {\MPEC}-ACQ does not hold at $x^*$.
	The only solution of \eqref{eq:lpeq_cq_fail} is $d = (\rho,0)$. 
	Therefore, this {\LPCC} cannot verify {geometric} B-stationarity for any $\rho>0$. 
	Fig.~\ref{fig:mpec_cq_fail} illustrates the feasible sets of the {\MPEC} and the {\LPCC} for different $\rho$.
\end{myexample}

The {\LPCC}~\eqref{eq:lpec_reduced} is well-defined at feasible points. 
However, in the Phase~I algorithm of Section~\ref{sec:algorithm_phase_i}, we must also consider {\LPCCs} at infeasible points, where index sets $\I_{+0}, \I_{0+}, \I_{00}$ are not defined. 
Existing strategies for identifying active sets from infeasible points often require overly restrictive assumptions~\cite{Lin2006}. 
Inspired by~\cite[Sec.~3.2]{Scholtes2004}, alternative approaches instead include all linearized constraints, yielding the {\LPCC}~\cite{Kirches2022,Leyffer2007}:

\begin{subequations}\label{eq:lpec_full}
	\begin{align}
		\underset{d\in \R^{n}}{\mathrm{min}} \;  \quad & \nabla f(x^k)^\top d \\
		\textnormal{s.t.} \quad 
		&c_i(x^k)+ \nabla c_i(x^k)^\top d\geq0,  &\forall i \in \{1,\ldots,n_c\},\\
		&0 \leq  x_{1,i}^k + d_{1,i}  \perp x_{2,i}^k + d_{2,i}  \geq0, &\forall i \in \{1,\ldots,m\},\label{eq:lpec_full_comp} \\
		&\| d \|_{\infty} \leq \rho. \label{eq:lpec_full_tr}
	\end{align}
\end{subequations}
We refer to the {\LPCC}~\eqref{eq:lpec_reduced} as the ``reduced'' {\LPCC} and \eqref{eq:lpec_full} as the ``full'' {\LPCC}.  
{For feasible points $x^k$, the full {\LPCC} contains the feasible set of the reduced one.  
Moreover, since $d = 0$ is always feasible for $x^k\in \Omega$, the optimal objective values of both~\eqref{eq:lpec_reduced} and~\eqref{eq:lpec_full} are nonpositive, regardless of the trust region radius $\rho>0$.  
At infeasible points $x^k \notin \Omega$, the optimal objective value of~\eqref{eq:lpec_full} can take any value in $\R$.
Moreover, at such points, the trust-region radius $\rho$ must be large enough to ensure that the {\LPCC} is feasible, cf. Theorem~\ref{th:phase_i_feasibility}.}

Let us introduce the following notation shorthand $\I_{+0}^k := \I_{+0}(x^k)$ and  $\I_{0+}^k := \I_{0+}(x^k)$.
Denote the set of inactive inequality constraints by $\bar{\A}^k := \bar{\A}(x^k) = \{i \in \{1,\ldots,n_c\} \mid c_i(x^k) > 0 \}$. 
For a feasible point $x^k \in \Omega$, the solution sets of the two {\LPCCs} are related as follows.
The proof can be found in Appendix~\ref{sec:lemma_proofs}.
\begin{mylemma}\label{lem:full_vs_reduced_lpec}
	Let $x^k\in \Omega$ be a feasible point of the {\MPEC}~\eqref{eq:mpec}. 
	For all trust region radii that satisfy
	\begin{align}\label{eq:full_vs_reduced_lpec}
		\begin{split}
		0 <\rho <  \bar{\rho}  = \min \Big\{&
		\{ {x^k_{1,i}} \mid i \in \I_{+0}^k \}  
		\cup 
		\{ {x^k_{2,i}} \mid i \in \I_{0+}^k \}\\
		&
		\cup 
		\{ {c_i(x^k)} \|\nabla c_i(x^k)\|^{-1} \mid i \in \bar{\A}^k , \|\nabla c_i(x^k)\|\neq 0\}
		\Big\}, 
		\end{split}
	\end{align}
	the sets of the local minimizers of ~\eqref{eq:lpec_reduced} and~\eqref{eq:lpec_full} are identical.
	In the special case of  $\I_{00}^k = \{1,\ldots,m\}$ {and $\bar{\A}^k =\emptyset$} the reduced and full {\LPCC} coincide and in that case, $\bar{\rho} = \infty$.
\end{mylemma}

In other words, at a feasible point $x^k$ and for a sufficiently small trust-region radius, all branches that are not biactive at $x^k$ are eliminated, and the linearizations of inactive constraints remain inactive. 
Thus, the feasible sets of the two {\LPCCs} coincide.
On the other hand, for $\rho \geq \bar{\rho}$, the global minimizer $d = 0$ of~\eqref{eq:lpec_reduced} may correspond only to a local minimizer of~\eqref{eq:lpec_full}, cf.\ Example~\ref{ex:global_local_lpec}. 
Hence, global optimality of $d = 0$ is sufficient but not necessary to verify B-stationarity via \eqref{eq:lpec_full}~\cite{Scholtes2004}.

In the next example, we show an additional advantage and disadvantage of the full LPEC relative to the reduced LPEC.
\begin{figure}[t]
	\centering
	\subfloat[$\rho = 0.6$.]{\includegraphics[width=0.33\textwidth]{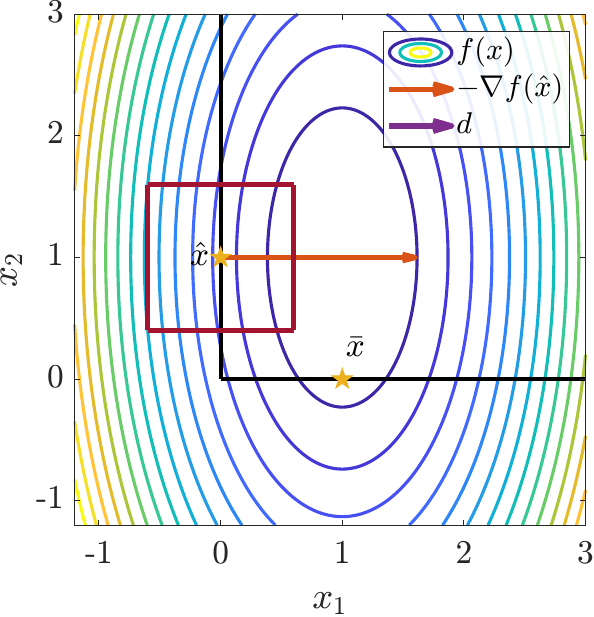}\label{fig:global_local_mpec1}}
	\subfloat[$\rho = 1.2$.]{\includegraphics[width=0.33\textwidth]{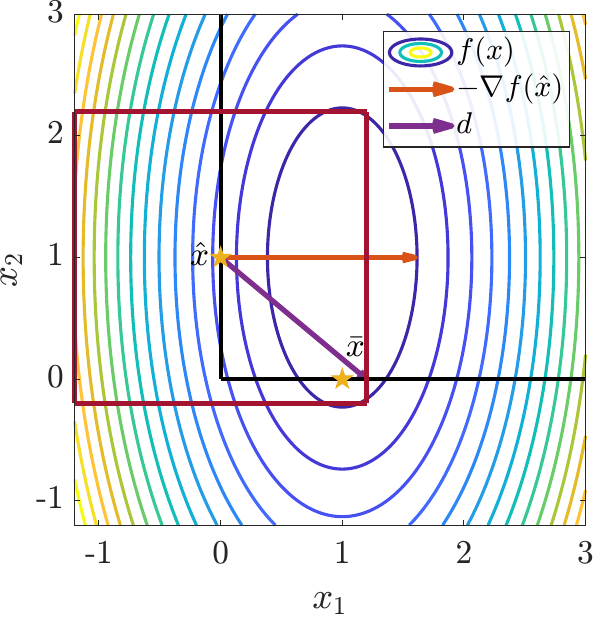}\label{fig:global_local_mpec2}}
	\subfloat[$\rho = 1.2$.]{\includegraphics[width=0.33\textwidth]{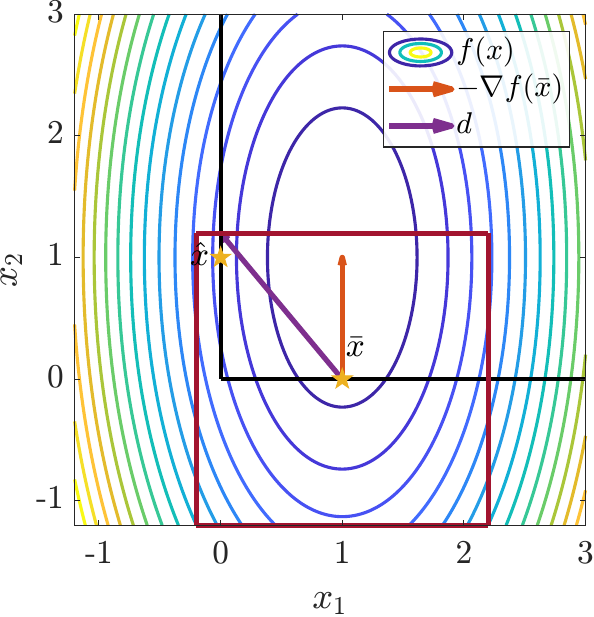}\label{fig:global_local_mpec3}}
	\caption{Illustration of the MPEC ~\eqref{eq:global_local_mpec}. 
		The yellow stars denote the two stationary points. 
		The figure illustrates the LPEC's~\eqref{eq:lpec_full} trust region, objective gradient, and optimal solution $d$ for different values of $\rho$ at different stationary points.}
	\label{fig:global_local_mpec}
\end{figure}
\begin{myexample}[Finding a better objective value]\label{ex:global_local_lpec}
	Consider the following {\MPEC}:
	\begin{align}\label{eq:global_local_mpec}
		\min_{x\in\R^2}\ 4(x_1-1)^2 + (x_2-1)^2 \quad \textnormal{s.t.} \quad 0 \leq x_1\perp x_2\geq 0,
	\end{align}
The example has two S-stationary points, $\bar{x} = (1,0)$ and $\hat{x} = (0,1)$, both with empty biactive sets, and respective objective values $f(\bar{x}) = 1$ and $f(\hat{x}) = 4$.
In Fig.~\ref{fig:global_local_mpec1}, if $\rho$ is sufficiently small or the reduced {\LPCC} is used, the point $\hat{x}$ with $f(\hat{x}) = 4$ is verified as a B-stationary point.
On the other hand, with a larger trust region radius, as depicted in Fig.~\ref{fig:global_local_mpec2}, the globally optimal solution of the full {\LPCC} identifies a BNLP with a lower objective value $f(\bar{x}) =1$.
Conversely, if we start at $\bar{x}$, the full {\LPCC} can identify a BNLP with a larger objective, as depicted in Fig.~\ref{fig:global_local_mpec3}. However, this step will be rejected, and the trust region radius reduced until $\bar{x}$ is verified as B-stationary.
\end{myexample}

We make several observations from this example.
First, solving the full {\LPCC} to global optimality can help us to find an objective with a lower value.
%Of course, as illustrated in Fig.~\ref{fig:global_local_mpec3}, there is no guarantee that this will always be the case.
Second, if the trust region is not small enough, there is no guarantee that the descent direction corresponding to the global optimum of the full {\LPCC} identifies a BNLP with a lower objective value.
{In Algorithm~\ref{alg:mpec_opt_summary} this leads to step rejection and a reduction of $\rho$.
Third, note that if we start from $\hat{x} \in \Omega$ and solve the BNLP predicted by the full {\LPCC} solution (Fig.~\ref{fig:global_local_mpec2}), then $\hat{x}\in \Omega$ is infeasible for the new BNLP. 
This may cause step rejection, and consequently reduction of $\rho$, if the NLP solver fails to solve the BNLP from this initial guess. 
However, for sufficiently small $\rho$, one can guarantee that any $x^k \in \Omega$ is feasible for $\BNLP(\I_1(x^k+d),\I_2(x^k+d))$, where $d$ feasible point of \LPCC$(x^k,\rho)$. 
The next lemma establishes this result.}

\begin{mylemma}\label{lem:reduced_lpec_feasible}
Let $x^k \in \Omega$ be a feasible point of the {\MPEC}~\eqref{eq:mpec} and suppose that the {\MPEC}-MFCQ holds at this point.
Let $d$ be any feasible point of the $\mathrm{LPEC}(x^k,\rho)$ in \eqref{eq:lpec_full} for $\rho \in (0,\bar{\rho})$, where $\bar{\rho}$ is defined in \eqref{eq:full_vs_reduced_lpec}.
Then $x^k$ is a feasible point of $\mathrm{BNLP}(\I_1(x^k+d),\I_2(x^k+d))$. 
\end{mylemma}
A proof is given in Appendix~\ref{sec:lemma_proofs}.

\subsubsection{On solving {\LPCCs}}\label{sec:solving_lpecs}
{\LPCCs} are nonconvex and combinatorial problems, and the effectiveness of our proposed method relies on solving them efficiently.
In practice, our computational results indicate that although the worst-case complexity of {\LPCCs} may be exponential, we typically only need to solve the root node of the MILP reformulation described below.

We summarize several favorable cases that frequently arise.
First, if the {\LPCC} contains only bound and complementarity constraints, it can be solved in linear time~\cite[Propositions 2.3 and 2.4]{Kirches2022}. This algorithm is incorporated into our implementation.
Second, at a feasible point $x^k$, if $\I_{00}(x^k) = \emptyset$, then the {\LPCC}~\eqref{eq:lpec_reduced} reduces to a standard LP. According to~\eqref{eq:mpec_feasible_sets}, the {\LPCC} can then be solved by enumerating all branch LPs, which is not computationally expensive when $|\I_{00}(x^k)|$ is small.
Third, if $x^k$ is not B-stationary, it is unnecessary to solve all branch LPs; it suffices to find a feasible direction $d$ with $\nabla f(x^k)^\top d < 0$ (cf. Corollary~\ref{cor:finite_termination}). 
Moreover, in the LPECs arising in Phase I of our method, it is sufficient to simply find a feasible point (cf. Theorem~\ref{th:phase_i_feasibility}). 
In a MILP formulation, the solver can therefore be terminated as soon as such a feasible point is obtained.
Section~\ref{sec:numerical_results} shows that this results in a significant reduction in the number of LPs solved in the MILP algorithm.
Finally, if $x^k$ is S-stationary, B-stationarity can be verified by solving a single LP (cf. Proposition~\ref{lem:milp_via_lp}).
In summary, except when verifying B-stationarity, it is rarely necessary to solve the LPEC to global optimality -- and even then, the task often reduces to solving a single LP.

In this paper, we focus on arguably the easiest-to-implement approach for solving \LPCCs, which does not suffer from the drawbacks of regularization methods.
We consider the equivalent MILP formulation of {\LPCC}~\eqref{eq:lpec_full}~\cite{Hu2012}:
\begin{subequations}\label{eq:lpec_milp}
	\begin{align}
		\underset{{d\in \R^{n},\ y\in\{0,1\}^m}}{\mathrm{min}} \;  \quad & \nabla f(x^k)^\top d \\
		\textnormal{s.t.} \quad 
		&c(x^k)+ \nabla c(x^k)^\top d\geq0, \\
		&0 \leq  x^k_{1,i} \!+\! d_{1,i} \!\leq\! {y_i} M,&\forall i \in \{1,\ldots,m\}, \label{eq:lpec_milp_1}\\
		&0 \leq  x^k_{2,i} \!+\! d_{2,i}\!\leq\! {(1-y_i)}M,&\forall i \in \{1,\ldots,m\},\label{eq:lpec_milp_2}\\
		& \| d \|_{\infty}  \leq \rho \label{eq:lpec_milp_tr}.
	\end{align}
\end{subequations}
The constraints \eqref{eq:lpec_milp_1} and \eqref{eq:lpec_milp_2} correspond to the big-M reformulation of the complementarity constraints in \eqref{eq:lpec_full_comp}, where $M$ is a sufficiently large positive constant.  
The binary variables $y\in \{0,1\}^m$ effectively select one of the branches of the L-shaped set defined by the complementarity constraints.

To reduce the problem complexity for a given $\rho$, if $d_i \leq \rho < x_{1,i}$ or $d_i \leq \rho < x_{2,i}$, we can fix the corresponding binary variable $y_i$ for every $i$ in a preprocessing step, as the other branch becomes unreachable.
In the special case of selecting $\rho$ as in~\eqref{eq:full_vs_reduced_lpec}, all binary variables $y_i$ for $i \notin \I_{00}(x^k)$ can be determined in advance, and it is observed that this contributes to the computational efficiency of our approach.

The MILP~\eqref{eq:lpec_milp} can be solved with state-of-the-art commercial and open-source solvers such as Gurobi~\cite{Gurobi} or HiGHS~\cite{Huangfu2018}.
We discuss some special properties of the MILP \eqref{eq:lpec_milp}, which are favorable for MILP solvers~\cite{Gurobi,Huangfu2018}.
We denote the optimal objective value of \eqref{eq:lpec_milp} by $V_{\mathrm{MILP}}(x^k)$, and of its LP relaxation, that is, a problem where $y\in\{0,1\}^m$ is replaced by $y\in [0,1]^m$, by $V_{\mathrm{LP}}(x^k)$.
Thus, it holds that $V_{\mathrm{LP}}(x^k) \leq V_{\mathrm{MILP}}(x^k)$. 
The tighter the relaxation, i.e., the smaller $|V_{\mathrm{MILP}}(x^k) - V_{\mathrm{LP}}(x^k)|$ is, the less branching may be needed in MILP solvers~\cite{Wolsey1998}.
Big-M formulations can be inefficient because they can lead to rather loose relaxations of the MILP~\cite{Hu2008}.
However, for $\rho < \frac{M}{2}$, the trust-region constraint~\eqref{eq:lpec_milp_tr} always dominates the bounds on $d$ from the big-M constraints~\eqref{eq:lpec_milp_1} and \eqref{eq:lpec_milp_2}.
Thus, $M$ only needs to be large enough to ensure the feasibility of~\eqref{eq:lpec_milp_1}--\eqref{eq:lpec_milp_2}.
A suitable choice is $M = \max(\max(x_1^k, x_2^k)) + \rho$.
In addition, observe that the constraint \eqref{eq:lpec_milp_tr} makes the feasible set of the MILP compact. 
Moreover, for feasible points $x^k$, $\rho$ can be arbitrarily small, leading to a smaller $V_{\mathrm{LP}}(x^k)$ and thus tighter relaxations of the MILP.

Furthermore, given a feasible complementarity pair, e.g. from the previous BNLP, a feasible binary vector $y$ can be constructed, {providing an upper bound for $V_{\mathrm{MILP}}(x^k)$. 
The availability of an upper bound can reduce branching in branch-and-bound MILP algorithms~\cite{Wolsey1998}, and our implementation includes the procedure for computing such a feasible binary vector.
Our numerical experiments indicate that these factors typically result in small number of LP that need to be solved when solving the MILP~\eqref{eq:lpec_milp}.
Interestingly, in our experiments the MILP formulations for {\LPCCs} are often substantially faster than solving the \LPCC~\eqref{eq:lpec_full} with noncombinatorial regularization-based methods, cf. Appendix~\ref{sec:numerics_macmpec_lpec}.

We discuss a special case that frequently arises in our numerical experiments, namely when $x^* \in \Omega$ is S-stationary.
In this case, the MILPs are particularly easy to solve since the relaxation is tight, i.e., $V_{\mathrm{MILP}}(x^*) = V_{\mathrm{LP}}(x^*)$.
Thus, verifying B-stationarity reduces to solving a single LP in the MILP solver.
We prove this for the reduced {\LPCC}~\eqref{eq:lpec_reduced}.  
However, if $\rho < \bar{\rho}$ (with $\bar{\rho}$ from \eqref{eq:full_vs_reduced_lpec}) and the preprocessing step for fixing integers is applied, then the feasible sets of the full MILP~\eqref{eq:lpec_milp} and reduced MILP~\eqref{eq:lpec_milp_reduced} coincide.

\begin{myproposition}\label{lem:milp_via_lp}
	Let $x^*\in \Omega$ be an S-stationary point of the {\MPEC}~\eqref{eq:mpec}
	at which the MPEC-MFCQ holds. 
	Consider the MILP reformulation of the reduced LPEC~\eqref{eq:lpec_reduced}:
	\allowdisplaybreaks{}
	\begin{subequations}\label{eq:lpec_milp_reduced}
		\begin{align}
			\underset{d\in \R^{n}\!, y\in\{ \! 0,1 \!\}^{|\I_{00}^*|}}{\mathrm{min}} \;  \quad & \nabla f(x^*)^\top d \\
			\textnormal{s.t.} \quad 
			&c_i(x^*)+ \nabla c_i(x^*)^\top d\geq0,  &\forall i\in \A(x^*),\\
			&x_{1,i}^* + d_{1,i}= 0,\ x_{2,i}^*+d_{2,i} \geq 0, &\forall i \in 	\I_{0+}(x^*), \\
			&x_{1,i}^* + d_{1,i} \geq 0,\ x_{2,i}^* +d_{2,i}= 0, &\forall  i \in 	\I_{+0}(x^*), \\
			&0 \leq  x^*_{1,i} \!+\! d_{1,i} \!\leq\! {y_i} M,  &\forall i \in \I_{00}(x^*), \label{eq:lpec_milp_reduced_comp1}\\
			& 0 \leq  x^*_{2,i} \!+\! d_{2,i}\!\leq\! {(1-y_i)}M,  &\forall i \in \I_{00}(x^*), \label{eq:lpec_milp_reduced_comp2}\\
			& \| d \|_{\infty}  \leq \rho. 
		\end{align}
		\end{subequations}
	For $\rho<\frac{M}{2}$, $d=0$ is a global minimizer of both the MILP~\eqref{eq:lpec_milp_reduced} and its relaxation with $y_i \in [0,1]$.
\end{myproposition}
\textit{Proof.} 
First we show that $d = 0$ is a global minimizer of the  MILP~\eqref{eq:lpec_milp_reduced}. 
It follows from Theorem~\ref{th:b_and_s_stat} that if $x^*\in \Omega$ is S-stationary it is also B-stationary.
Hence, $d=0$ is a global minimizer of the {\LPCC}~\eqref{eq:lpec_reduced_theory} and of $\LPEC(x^*,\rho)$ in~\eqref{eq:lpec_reduced}, because the additional trust-region constraint in the latter remains inactive. 
The MILP~\eqref{eq:lpec_milp_reduced} and LPEC~\eqref{eq:lpec_reduced} are equivalent, hence $d =0$ is a global minimizer of the MILP.
% Moreover, $d = 0$ is the global minimizer of the MILP~\eqref{eq:lpec_milp_reduced}, which is equivalent to~\eqref{eq:lpec_reduced}.

Second, we show that $d = 0$ is global minimizer of the  relaxed MILP~\eqref{eq:lpec_milp_reduced} with $y_i \in [0,1]$.
Since $x^*$ is S-stationary, it is also a stationary point of the RNLP$(x^*)$.
Consider the following LP:
\begin{subequations}\label{eq:relaxed_lp}
	\begin{align}
		\underset{{d\in \R^{n}}}{\mathrm{min}} \;  \quad & \nabla f(x^*)^\top d \\
		\textnormal{s.t.} \quad 
		&c_i(x^*)+ \nabla c_i(x^*)^\top d\geq0, \ \forall i\in \A(x^*), \label{eq:relaxed_lp_c}\\
		&x_{1,i}^* + d_{1,i}= 0,\ x_{2,i}^*+d_{2,i} \geq 0,\; &\forall i \in 	\mathcal{I}_{0+}(x^*), \\
		&x_{1,i}^* + d_{1,i} \geq 0,\ x_{2,i}^* +d_{2,i}= 0,\; &\forall  i \in 	\mathcal{I}_{+0}(x^*), \\
		&x_{1,i}^* +d_{1,i} \geq 0,\ x_{2,i}^*+ d_{2,i} \geq 0,\; &\forall  i \in 	\mathcal{I}_{00}(x^*), \label{eq:relaxed_lp_biactive}\\
		& \| d \|_{\infty}  \leq \rho 
	\end{align}
\end{subequations}
The MPEC-MFCQ implies that the standard MFCQ holds for the RNLP$(x^*)$.
Observe that constraints \eqref{eq:relaxed_lp_c}--\eqref{eq:relaxed_lp_biactive} define the tangent cone of RNLP$(x^*)$ at $x^*$.
Since $x^*$ is a stationary point of RNLP$(x^*)$, it follows that $d = 0$ is a global minimizer of \eqref{eq:relaxed_lp} (since this exactly recovers the KKT conditions of the RNLP at $x^*$ from the LP~\eqref{eq:relaxed_lp}).
The trust-region constraint remains inactive and does not affect the solution.
Note that this LP is a relaxation of the reduced LPEC~\eqref{eq:lpec_reduced} (and thus of the MILP~\eqref{eq:lpec_milp_reduced}).
%Observe also that this LP coincides with the relaxed LP of $\LPEC(x^*,\rho)$ in~\eqref{eq:lpec_reduced} at $d = 0$. 

In the relaxed MILP, $y_i \in \{0,1\}$ in \eqref{eq:lpec_milp_reduced} is replaced by $y_i \in [0,1]$.
Next, from the upper bounds in \eqref{eq:lpec_milp_reduced_comp1} and \eqref{eq:lpec_milp_reduced_comp2}, we have that  $ \frac{d_{1,i}}{M}  \leq y_i  \leq 1 - \frac{d_{2,i}}{M}$.
Since the MILP is feasible, it follows that this interval is nonempty if $ \frac{d_{1,i}}{M}   \leq 1 - \frac{d_{2,i}}{M})$ holds.
This implicitly defines the constraint:
$ d_{1,i}+ d_{2,i} \leq M$. 
If we choose $\rho < \tfrac{M}{2}$, then the upper bounds in \eqref{eq:lpec_milp_reduced_comp1} and \eqref{eq:lpec_milp_reduced_comp2} are never active.
For such a choice, the feasible sets of the relaxed MILP~\eqref{eq:lpec_milp_reduced} and the LP~\eqref{eq:relaxed_lp} coincide for all components $d \in \R^n$, for all $y \in [0,1]^{|\I_{00}|}$.
Since the objective of \eqref{eq:lpec_milp_reduced} does not depend on $y$, it follows that $d = 0$ is also a global minimizer of the relaxation of \eqref{eq:lpec_milp_reduced}.
In conclusion, $d =0$ is a global solution of the MILP~\eqref{eq:lpec_milp_reduced} and {\LPCC}~\eqref{eq:lpec_reduced} (verifies B-stationarity), relaxed LP (verifies S-stationarity) and relaxed MILP.
\qed

\begin{myremark}[On Algorithm~\ref{alg:mpec_opt_summary} and mixed-integer algorithms]
	Any {\MPEC}~\eqref{eq:mpec} can be reformulated into a mixed-integer NLP (MINLP) using the big-M reformulation:
	\begin{subequations}\label{eq:mpec_minlp}
		\begin{align}
			\underset{x\in \R^{n}, y\in\{0,1\}^m}{\mathrm{min}} \;  \quad &f(x)\\
			\textnormal{s.t.} \quad 
			&c(x)\geq0, \label{eq:minlp_ineq}\\
			&0 \leq x_{1,i} \leq y_i M , & \forall i \in \{1,\ldots,m\},\\
			&0 \leq x_{2,i} \leq (1-y_i)M, & \forall i \in \{1,\ldots,m\}.
		\end{align}
	\end{subequations}
	However, the convergence properties of our method and MINLP methods differ in several aspects. 
	Most MINLP methods converge to the global optimum when the relaxed integer problem is convex~\cite{Belotti2013}, but Algorithm~\ref{alg:mpec_opt_summary} does not.
	This distinction is illustrated by the {\MPEC} example in Example~\ref{ex:global_local_lpec}, whose equivalent MINLP (after relaxing $y\in\{0,1\}^m$ to $y\in [0,1]^m$) yields a convex quadratic program. 
	There, Algorithm~\ref{alg:mpec_opt_summary} converges only to a B-stationary point, not necessarily the global optimum.

	While Algorithm~\ref{alg:mpec_opt_summary} solves NLPs and MILPs similar to many MINLP methods~\cite{Belotti2013}, there are significant differences in termination criteria and the particular subproblems that are solved.
	Nonlinear branch-and-bound algorithms terminate when the optimality gap falls below a threshold~\cite{Belotti2013}, which does not directly relate to B-stationarity unless the global optimum is found. 
	Outer approximation algorithms also solve sequences of NLPs and MILPs~\cite{Belotti2013}, but their MILPs include linearized constraints from multiple different linearization points, making them different from the MILP~\eqref{eq:lpec_milp}.
	Consequently, it is not immediately clear when these methods terminate at B-stationary points.
	In summary, Algorithm~\ref{alg:mpec_opt_summary} is not a MINLP algorithm; rather, it employs MILP solvers as a tool for solving {\LPCCs}. 
	It aims to find a local optimum of an MPEC, which—unlike in MINLPs -- is a well-defined notion.
\end{myremark}

\subsection{Phase I: Finding a feasible branch NLP}\label{sec:algorithm_phase_i}
In this section, we introduce an algorithm to find an initial feasible BNLP and its corresponding solution $x^0$, or {return a stationary point of an auxiliary problem that minimizes a measure of constraint infeasibility.}

Regularization-based methods may not always find B-stationary points, but they can be efficient in finding (almost) feasible points.  
We outline a Phase I method that combines solving {\LPCCs} with regularization-based methods.
While any regularization-based method can be used in this framework, we focus in our exposition on Scholtes' global relaxation method, which performs well in numerical experiments~\cite{Hoheisel2013,Nurkanovic2024b}.
It solves a sequence of NLPs of the form:
\begin{subequations}
	\label{eq:mpec_scholtes}
	\begin{align}
		\underset{x\in \R^{n}}{\mathrm{min}} \;  \quad &f(x)\\
		\textnormal{s.t.} \quad 
		&c(x)\geq0, \\
		&x_1, x_2 \geq 0, \\
		&x_{1,i} x_{2,i} \leq \tau, &\forall i \in \{1,\ldots,m\}. 
	\end{align}
\end{subequations}
We often denote the problem~\eqref{eq:mpec_scholtes} by $\mathrm{Reg}(\tau)$ to highlight its dependency on the parameter $\tau$.
We compute $x^*(\tau)$, a stationarity point of \eqref{eq:mpec_scholtes}, which is not necessarily feasible for the {\MPEC}~\eqref{eq:mpec}. 
An {\LPCC} is then constructed at $x^*(\tau)$ (for some $\tau>0$), and its solution determines an active set corresponding to a feasible BNLP.

It is well-known that certifying infeasibility of a general nonconvex optimization problem is difficult.  
In many NLP solver implementations, if the globalization strategy cannot make sufficient progress toward a feasible point, a feasibility restoration phase is triggered, cf. \cite[Sec. 3.3.]{Fletcher2002}, \cite[Sec. 3.3.]{Waechter2006}, and \cite[Sec. 6]{Byrd2006}.  
In this step, an auxiliary NLP that minimizes the constraint violation is solved~\cite{Byrd2006,Fletcher2002,Waechter2006}.
If its stationary point is infeasible for the original problem, the solver terminates.
Infeasibility of the original problem is sufficient but not necessary for this outcome.
This outcome is often called local infeasibility~\cite{Waechter2006}, and it is not a global certificate of infeasibility.

\begin{mylemma}\label{lem:phase_i_infeasibility}
Regard the {\MPEC}~\eqref{eq:mpec} and the corresponding relaxed problem $\mathrm{Reg}(\tau)$ in~\eqref{eq:mpec_scholtes} for some $\tau>0$. 
If $\mathrm{Reg}(\tau)$ is infeasible then the same holds for the {\MPEC}~\eqref{eq:mpec}.
\end{mylemma}
\textit{Proof.}
This follows from the fact that the feasible set of $\mathrm{Reg}(\tau)$ is a superset of~$\Omega$. 
\qed

On the other hand, if \ref{eq:mpec_scholtes} is feasible for all $\tau$, then the {\MPEC}~\eqref{eq:mpec} is also feasible.
We expect that if $\tau^k$ is sufficiently small, then the feasible set of $\mathrm{LPEC}(x^*(\tau^k),\rho)$ is locally a good approximation of the original feasible set, and can help identify a feasible BNLP.
Algorithm~\ref{alg:phase_i} details such a procedure.
Before going into details, we illustrate how an {\LPCC} solution may fail or succeed in identifying a feasible BNLP for different values of $\tau$.
\begin{myexample}\label{ex:feasibility_lpec}
	Consider the following {\MPEC}:
	\begin{align*}
		\underset{x\in \R^{2}}{\mathrm{min}} \;  \quad & -2x_1 + x_2 \quad
		\textnormal{s.t.} \
		-x_1-(x_2-a)^2+1\geq0,\
		0\leq x_1 \perp x_2\geq0,
	\end{align*}
	with $a>1$. 
	In the example, we set $a=1.1$. 
	For larger values of $\tau$ both branches are feasible for the {\LPCC}~\eqref{eq:lpec_full}, but only the branch $x_1 = 0, x_2 \geq 0$ is feasible for the {\MPEC}.
	On the other hand, for smaller values of $\tau$, also in the LPEC, only the branch $d_1 = 0,\ d_2 \geq 0$ is feasible, thus the LPEC identifies only active sets $(\I_1(x^*(\tau)+d),\I_2(x^*(\tau)+d))$, which results in a feasible BNLP.
	Fig.~\ref{fig:feasibility_lpec} illustrates two sample solutions for its $\mathrm{Reg}(\tau)$ relaxation (left plots) and the corresponding $\mathrm{LPEC}(x^*(\tau),\rho)$ (right plots).
	In this example, the smaller the parameter $a$ is the smaller $\tau$ must become.
\end{myexample}
\begin{figure}[t]
	\centering
	\subfloat[$\tau = 1.0$.]{\includegraphics[width=0.35\textwidth]{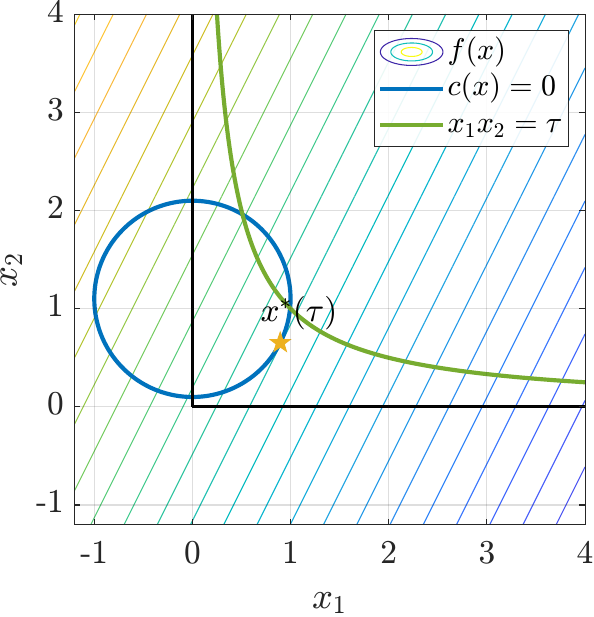}\label{fig:phase1_mpec_1}}
	\subfloat[$\tau = 1.0$.]{\includegraphics[width=0.35\textwidth]{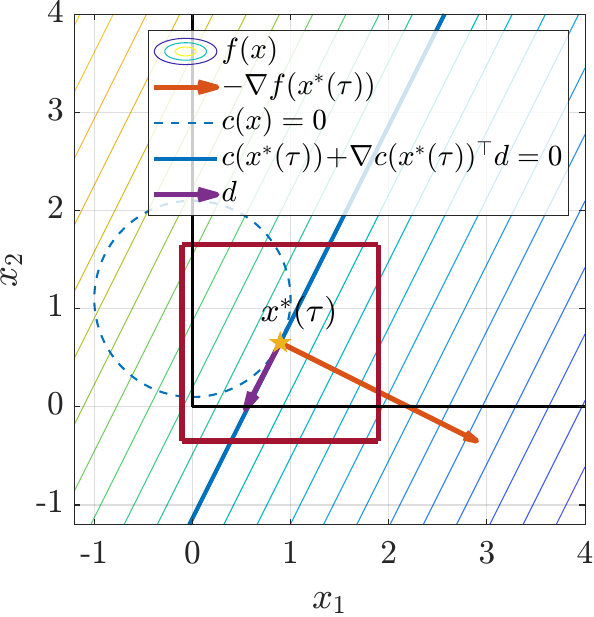}\label{fig:phase1_lpec_1}}\\
	\vspace{-0.34cm}
	\subfloat[$\tau = 0.04$.]{\includegraphics[width=0.35\textwidth]{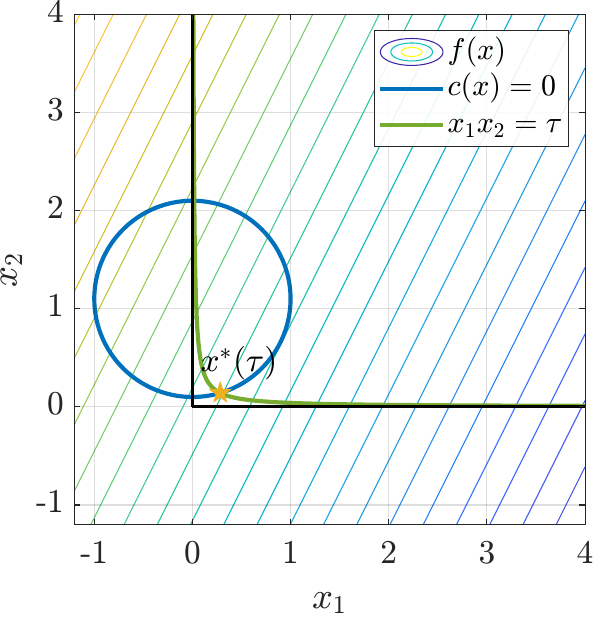}\label{fig:phase1_mpec_2}}
	\subfloat[$\tau  =0.04$.]{\includegraphics[width=0.35\textwidth]{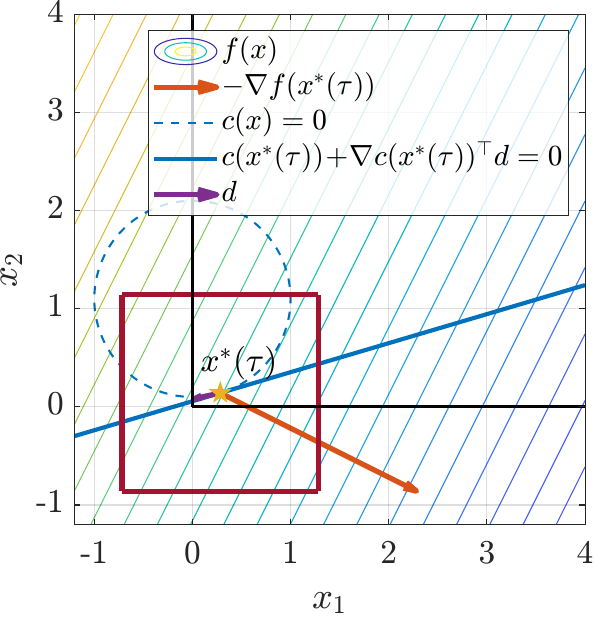}\label{fig:phase1_lpec_2}}
	\vspace{-0.2cm}
	\caption{Illustration of Example~\ref{ex:feasibility_lpec}, where the {\LPCC}~\eqref{eq:lpec_full} can correctly and wrongly identify a feasible BNLP depending on the linearization point $x^*(\tau)$.
	{Top: the linearization point $x^*(\tau)$ is not close enough to the feasible set, so the LPEC admits an MPEC-infeasible branch.
	Bottom: $x^*(\tau)$ is sufficiently close, so all LPEC branches are MPEC-feasible.}
	}
	\label{fig:feasibility_lpec}
\end{figure}

\begin{algorithm}[h!]
	\caption{Regularization-based Phase I of {\MPECopt}}
	\label{alg:phase_i}
	\begin{algorithmic}[1]
		\Statex \textbf{Input:} $\rho_I >0$, $\kappa \in (0,1)$, $\tau^0>0$, $N^\tau>1$  
%		 $N^\rho  >1$
%		$\gamma^{\mathrm{U}} >1$
		%		\Statex \textbf{Phase I (finding a feasible point)}
		\While{}
%		\For{ $k=0,\ldots, N^\tau$}
		\State Solve $\mathrm{Reg}(\tau^k)$ for 
		$x^*(\tau^k)$
		%		$x^*(\tau^k)$
		\If{$\mathrm{Reg}(\tau^k)$ locally infeasible}
		\State \textbf{break} \Comment{\textcolor{gray}{{\MPEC} locally infeasible}} \label{alg:reg_infeasible}
		\ElsIf{$x^*(\tau^k) \in \Omega$ }
		\State \textbf{break}; return $x^0 = x^*(\tau^k)$ \label{alg:reg_feasible} \Comment{\textcolor{gray}{$\mathrm{Reg}(\tau^k)$ solution feasible for {\MPEC}}}
		\ElsIf { $h_{\perp}(x^*(\tau^k)) < \rho_I$ } \label{alg:phase_i_start_lpec}
		\State Solve $\mathrm{LPEC}(x^*(\tau^k),\rho_I)$ for $d^{k}$ \label{alg:phase_i_solve_lpec}
		\If{ $d^{k}$ exists} \label{alg:phase_i_end_lpec}
			\State Solve $\mathrm{BNLP}(\I_1(x^*(\tau^k)+d^k), \I_2(x^*(\tau^k)+d^k))$ for $x^0$.\label{alg:bnlp_solve} 
				\If{Solution exists}
				\State \textbf{break}, return $x^0$ \label{alg:bnlp_solved} 
				\Comment{\textcolor{gray}{feasible point found}}
			\EndIf
		\EndIf
		\Else 
		\State $\tau^{k+1} = \kappa \tau^k$.  \label{alg:tau_update} 
%		, reset $\rho^{k,0} \in [\rho^{\mathrm{lb}},\rho^{\mathrm{ub}}]$
		\Comment{\textcolor{gray}{no feasible LPEC/BNLP found, reduce homotopy parameter}}
		\EndIf
				\EndWhile
%		\EndFor
	\end{algorithmic}
\end{algorithm}
We show below in Corollary~\ref{col:phase_i_feasibility} that, if the MPEC admits a feasible point, then for sufficiently small $\tau$, a feasible point of the LPEC yields a complementarity active set corresponding to a feasible BNLP.
Consequently, the while loop in Algorithm~\ref{alg:phase_i} has finite termination.
In a practical implementation, a point is declared feasible if the constraint residual is below a certain tolerance~{\texttt{tol}$_h$}.
Formally, we require $h(x)\leq \texttt{tol}_h$, where $h(x)$ is the total infeasibility defined as:
%\begin{align}\label{eq:total_infeasiblity}
%	%	\label{eq:infeasiblity}
%	&h(x) \!=\! \|(h_c(x),h_{\!\perp\!}(x))\|_{\infty},\ \!\!
%	h_c(x) \!=\!\|c^-(x)\|_{\infty},\ \!
%	h_{\!\perp\!}(x) \!=\! \max(|\! \min(x_1,x_2) \! |),
%\end{align} 
\begin{align}\label{eq:total_infeasiblity}
	%	\label{eq:infeasiblity}
	\begin{split}
		&h(x) =  \max(h_c(x),h_{\!\perp\!}(x)),\\
		&h_c(x) = \max(c^-(x)),\ \!
		h_{\!\perp\!}(x) \!=\! \max(|\min(x_1,x_2)|),
	\end{split}	
\end{align} 
with $c^-(x) =\min(0,c(x))$.

If $\mathrm{Reg}(\tau^k)$ is infeasible (line~\ref{alg:reg_infeasible}), or its solution is feasible for the {\MPEC} (line~\ref{alg:reg_feasible}), the algorithm terminates.  
Otherwise, the solution $x^*(\tau^k)$ is used to construct $\LPEC(x^*(\tau^k),\rho_I)$.  
The fixed trust-region radius $\rho_I>0$ must be sufficiently large to ensure feasibility of this {\LPCC}.  
Since a solution of $\mathrm{Reg}(\tau^k)$ already satisfies the general inequality constraints, i.e., $c(x^*(\tau^k)) \geq 0$, the smallest value of $\rho_I$ required for feasibility is $\rho_I > h_{\perp}(x^*(\tau^k))$.  
The {\LPCC} is solved only if this inequality holds (line~\ref{alg:phase_i_start_lpec}).
If it is not yet satisfied or the {\LPCC} is infeasible, the homotopy parameter is reduced to obtain a better linearization point $x^*(\tau^k)$ (line~\ref{alg:tau_update}).  
If a solution to the {\LPCC} is found, Algorithm~\ref{alg:phase_i} attempts to solve the corresponding $\mathrm{BNLP}(\I_1(x^*(\tau^k)+d^{k}),\I_2(x^*(\tau^k)+d^{k}))$ (line~\ref{alg:bnlp_solve}).  
If this is successful, a feasible point is obtained.  
Otherwise, the homotopy parameter is reduced to compute a better linearization point $x^*(\tau^{k+1})$.  

Theorem~\ref{th:phase_i_feasibility} shows that if $x^*(\tau^k)$ is sufficiently close to a feasible point $x\in \Omega$, then $\mathrm{BNLP}(\I_1(x^*(\tau^k)+d^{k}),\I_2(x^*(\tau^k)+d^{k}))$ is feasible, where $d^{k}$ is a feasible point of $\mathrm{LPEC}(x^*(\tau^k),\rho^{I})$.  
Moreover, under MPEC-MFCQ, the sequence $x^*(\tau^k)$ converges to a stationary and thus feasible point of the {\MPEC}~\cite{Hoheisel2013,Kanzow2015}.  

\begin{myremark}
	Algorithm~\ref{alg:phase_i} with the outcomes in lines~\ref{alg:reg_feasible} or~\ref{alg:bnlp_solved} can also be viewed as a standalone method that computes a stationary point of the {\MPEC}~\eqref{eq:mpec}, though not necessarily a B-stationary point. 
	However, if the computed point is B-stationary, this will be immediately verified in Phase II. 
	If the trust-region radius is sufficiently small (cf. Lemma~\ref{lem:full_vs_reduced_lpec}), the first {\LPCC} to be solved will have a solution $d=0$. 
	This behavior was occasionally observed in our numerical experiments in Section~\ref{sec:numerical_results}.
\end{myremark}

\subsection{Phase II: Computing a B-stationary point}\label{sec:algorithm_phase_ii}
Algorithm~\ref{alg:phase_ii} provides a detailed version of the Phase II algorithm outlined in Algorithm~\ref{alg:mpec_opt_summary}.
It consists of two nested loops: an outer loop with counter $k$, which makes progress toward B-stationary points, and an inner loop with counter $l$, which controls the trust-region radius $\rho^{k,l}$ and solves {\LPCCs} and BNLPs until a point $x^{k,l}\in \Omega$ with improved objective value is found, i.e., $f(x^{k,l}) < f(x^{k})$.
Under reasonable assumptions, the inner loop is guaranteed to terminate in finitely many steps; cf. Proposition~\ref{prop:finite_termination}.
With additional conditions, the outer loop may also terminate finitely; cf. Theorem~\ref{th:main_convergence}.

\begin{algorithm}
	\caption{Phase II of {\MPECopt}}
	\label{alg:phase_ii}
	\begin{algorithmic}[1]
		\Statex \textbf{Input:} $x^0 \in \Omega$, $\gamma^\mathrm{L} \in (0,1)$, $N^{\mathrm{out}}$, {$N^{\mathrm{in}}$,
%		$\rho^{0,0} \in [\rho^{\mathrm{lb}},\rho^{\mathrm{ub}} ]$
		$\rho^{\mathrm{lb}}$, $\rho^{\mathrm{ub}}$.
		}
		%		$\gamma^\mathrm{U} >1$
		\For{$k=0,\ldots,N^{\mathrm{out}}$} \Comment{\textcolor{gray}{major/outer iterations}}
		\State Reset trust region radius $\rho^{k,0} \in [\rho^{\mathrm{lb}},\rho^{\mathrm{ub}}]$;
		\For{$l=0,\ldots, N^{\mathrm{in}}$} \Comment{\textcolor{gray}{minor/inner iterations}}
		\State Solve $\mathrm{LPEC}(x^k,\rho^{k,l})$ for $d^{k,l}$ \label{alg:phase_ii_solve_lpec}
		%		\If{$\|d_\mathrm{L}^{k,l}\|_{\infty} \leq \texttt{tol}_{\mathrm{B}}$} \label{alg:phase_ii_check_B}
		\If{$d^{k,l}=0$} \label{alg:phase_ii_check_B}
		\State \textbf{terminate} \Comment{\textcolor{gray}{B-stationary point found}} \label{alg:phase_ii_terminate}
		\Else
		\State Select partition $(\I_1^*(x^k+d^{k,l}),\I_2^*(x^k+d^{k,l}))$ via Eq. \eqref{eq:lpec_select_partition} \label{alg:phase_ii_select_partition}
		\If{$l=0$ or $\I_1^*	(x^k+d^{k,l})\neq \I_1^*(x^k+d^{k,l-1})$} \label{alg:phase_ii_new_partition_selected}
		\State {Solve $\mathrm{BNLP}(\I_1^*(x^k+d^{k,l}), \I_2^*(x^k+d^{k,l}))$ to obtain $x^{k,l}$} \label{alg:phase_ii_solve_bnlp}
		\If{BNLP solver failed or $f(x^{k,l}) \geq f(x^{k})$ } \label{alg:phase_ii_check_objective_decerase}
		\State ${\rho}^{k,l+1}={\gamma^\mathrm{L} \rho^{k,l}}$; 
		\Comment{\textcolor{gray}{reduce TR radius}} \label{alg:phase_ii_reduce_tr_bnlp}
		\Else
		%		\State ${\rho}^{k,l+1}={\gamma^\mathrm{U} \rho^{k,l} }$		\Comment{\textcolor{gray}{increase TR radius}}
		\State $x^{k+1} = x^{k,l}$; \textbf{break} inner for-loop;   \Comment{\textcolor{gray}{step accepted}}	\label{alg:phase_ii_step_accept}
		\EndIf 			
		\Else 
		\State ${\rho}^{k,l+1}={\gamma^\mathrm{L} \rho^{k,l}}$ \label{alg:phase_ii_} \Comment{\textcolor{gray}{reduce TR radius}} \label{alg:phase_ii_reduce_tr_lpec}
		\EndIf 
		\EndIf 
		\EndFor
		\EndFor
	\end{algorithmic}
\end{algorithm}

Algorithm~\ref{alg:phase_ii} starts from a feasible point $x^0 \in \Omega$. 
At the beginning of each outer iteration, the trust-region radius is reset to $\rho^{k,0} \in [\rho^{\mathrm{lb}},\rho^{\mathrm{ub}} ]$ with $\rho^{\mathrm{ub}} > \rho^{\mathrm{lb}} > 0$. 
In theory any initial $\rho^{k,0}$ from the interval can be picked. 
In our implementation, we set $\rho^{k,0} = 10^{-3}$, with $\rho^{\mathrm{ub}} = 10^{3}$, and $\rho^{\mathrm{lb}} = 10^{-7}$.

Given the current iterate $x^k$ and trust-region radius $\rho^{k,l}$, we solve the {\LPCC}$(x^{k}, \rho^{k,l})$ for a step $d^{k,l}$, using its equivalent MILP reformulation~\eqref{eq:lpec_milp} (line~\ref{alg:phase_ii_solve_lpec}). 
If $d^{k,l} = 0$, then a B-stationary point is found (line~\ref{alg:phase_ii_terminate}). 
As discussed in Sec.~\ref{sec:algorithm_summary}, if $\nabla f(x^k)^\top d^{k,l} = 0$, then $d^{k,l} = 0$ is a global minimizer that can be selected.

If $d^{k,l} \neq 0$, the algorithm proceeds by selecting and solving the corresponding BNLP (lines~\ref{alg:phase_ii_select_partition}--\ref{alg:phase_ii_solve_bnlp}). 
If $\I_{00}(x^k+d^{k,l}) \neq \emptyset$, the choice of the next partition is not unique.
Proposition~\ref{prop:finite_termination} ensures that there exists a threshold $\bar{\rho}$ such that for all $\rho^{k,l} <\bar{\rho}$, every possible partition from  $\mathcal{P}(x^k+d^{k,l})$ yields a BNLP whose local minimizer improves the objective. 
Therefore, we may pick any partition from this set.
We denote the selected partition by $(\I_1^*(x^k+d^{k,l}),\I_2^*(x^k+d^{k,l})) \in \mathcal{P}(x^k+d^{k,l})$.
In our implementation, we use the MILP solution to pick the partition:
\begin{align}\label{eq:lpec_select_partition}
	\I^*_1(x^k+d^{k,l}) = \{ i \mid y_i^{k,l} = 0\}, \ \I^*_2(x^k+d^{k,l}) = \{ i \mid y_i^{k,l} = 1\}.
\end{align}
Here the integer values $y_i^{k,l}$ implicitly select the corresponding branches (cf. Eq.~\eqref{eq:lpec_milp_1}--\eqref{eq:lpec_milp_2}).
Alternative partitioning rules are also valid; for instance, we also implemented a rule based on the objective gradient components corresponding to complementarity variables $x_{1,i}$ and $x_{2,i}$~\cite{Kirches2022}:
\begin{align}\label{eq:lpec_select_rule2}
	\begin{split}
	&\mathcal{D}_1^*(x^k \!+\! d^{k,l}) =\{ i \in \I_{00}(x^k \!+\! d^{k,l}) \! \mid \!  \nabla f_{1,i}(x^k \!+\! d^{k,l}) \geq  \nabla f_{2,i}(x^k+d^{k,l})   \}, \\
	&\mathcal{D}_2^*(x^k+d^{k,l}) =	\I_{00}(x^k+d^{k,l}) \setminus \mathcal{D}_1^*(x^k+d^{k,l}), \\
	&\I^*_1(x^k+d^{k,l}) = \I_{0+}(x^k+d^{k,l}) \cup \mathcal{D}_1^*(x^k+d^{k,l}) \\
	&\I^*_2(x^k+d^{k,l}) = \I_{+0}(x^k + d^{k,l}) \cup \mathcal{D}_2^*(x^k + d^{k,l}).
	\end{split}
\end{align}
Since the BNLPs are in general nonconvex problems, starting from $x^k$ can give different objective improvements $f(x^k)-f(x^{k,l}) > 0$ for different BNLPs.
For convergence, it is not necessary to pick the best or a particular one, but only sufficient decrease must be guaranteed by the NLP solver, cf.~Theorem~\ref{th:main_convergence}.
In our numerical experiments, we did not observe performance differences between using rule~\eqref{eq:lpec_select_partition} and~~\eqref{eq:lpec_select_rule2}.

If $\rho^{k,l} > \bar{\rho}$, the partition may not change between iterations, which is why line~\ref{alg:phase_ii_new_partition_selected} checks whether a new BNLP needs to be solved.
If no new partition is selected, the trust-region radius is reduced (line~\ref{alg:phase_ii_reduce_tr_lpec}).
Moreover, for larger $\rho^{k,l}$, if the current BNLP does not improve the objective (cf. Example~\ref{ex:global_local_lpec}), the trust-region radius is reduced as well (line~\ref{alg:phase_ii_reduce_tr_bnlp}).
In numerical experiments, we observed that allowing larger initial trust-region radii often accelerates convergence, though they are not necessary for convergence.

\section{Convergence analysis}\label{sec:convergence_theory}
This section provides the convergence proof of the proposed method.
Before stating the main results, we discuss the main assumption and some auxiliary results.

\subsection{Assumptions}\label{sec:convergence_theory_assumptions}
Since the feasible set of any BNLP is a subset of the MPEC’s feasible set, an unbounded MPEC can be detected via an unbounded BNLP. 
Most NLP solvers readily detect this. 
Hence, in our convergence analysis we restrict our attention to bounded problems.
\begin{myassumption}\label{ass:compactness}
	There exists a compact set $X$ such that for the feasible set $\Omega$ of the {\MPEC}~\eqref{eq:mpec} it holds that $\Omega\subset X$.
\end{myassumption}

Some results in the previous section required only the {\MPEC}-ACQ.
For the convergence results in this section, we require a stronger assumption.
\begin{myassumption}\label{ass:cq}
	The {\MPEC}-MFCQ holds at any feasible point of the {\MPEC}~\eqref{eq:mpec}. 
\end{myassumption}
{\MPEC}-MFCQ implies {\MPEC}-ACQ~\cite[Theorem 3.1]{Flegel2005a}, and therefore the degenerate case of Example~\ref{ex:mpec_cq_fail} is avoided.
Other useful consequences are as follows. 
\begin{myproposition}\label{prop:mfcq_consequences}
	Suppose that Assumption~\ref{ass:cq} holds at $x^*\in \Omega$ for~\eqref{eq:mpec}. Then:
	\begin{enumerate}
		\item the usual MFCQ is satisfied at $x^*$ by the corresponding $\BNLP(\I_1(x^*),\I_2(x^*))$ for any partition $(\I_1(x^*),\I_2(x^*)) \in \mathcal{P}(x^*)$,
		\item there exists a neighborhood $\mathcal{N}$ of $x^*$ and a constant $\bar{\tau}>0$ such that the usual MFCQ for $\mathrm{Reg}(\tau) $ in~\eqref{eq:mpec_scholtes} is satisfied at all $x \in \mathcal{N} \cap \Omega_{\mathrm{Reg}(\tau)}$ for all $\tau \leq \bar{\tau}$, where $\Omega_{\mathrm{Reg}(\tau)}$ is the feasible set of $\mathrm{Reg}(\tau)$ in~\ref{eq:mpec_scholtes}; and,
		\item independently of Assumption~\ref{ass:cq}, the tangent cones of the MPEC and the BNLPs satisfy
		\begin{align}
			\mathcal{T}_{\Omega}(x^*) 
			= \bigcup_{(\I_1^*,\I_2^*) \in \mathcal{P}^*} 
			\mathcal{T}_{\Omega_{\BNLP(\I_1^*,\I_2^*)}}(x^*).
		\end{align}
		Moreover, under Assumption~\ref{ass:cq}, the left-hand side can be replaced by the MPEC-linearized feasible cone $\mathcal{F}_{\Omega}(x^*)$, and the terms on the right-hand side by the usual linearized feasible cones of the BNLPs; cf. Eq.~\eqref{eq:mpec_lin_cones}.
	\end{enumerate}
\end{myproposition}
\textit{Proof}. The first fact is~\cite[Lemma 3.3]{Flegel2005a}, the second is~\cite[Theorem 3.2]{Hoheisel2013}, and the third follows from~\cite[Lemma 3.1]{Flegel2005a} and \cite[Theorem 3.1]{Flegel2005a}. \qed

Most \MPEC\ methods require solving a sequence of NLPs, e.g., \eqref{eq:mpec_scholtes} with a decreasing sequence $\{\tau^k\}$.
In their convergence analysis, it is typically assumed that a constraint qualification holds and that a stationary point of such an NLP is available, without further specifying what assumptions on the NLP solver are required to compute such a stationary point~\cite{Hoheisel2013,Kanzow2015,Ralph2004,Steffensen2010}.
If a minimizer, and not only a stationary point, is required, then a suitable MPEC second-order sufficient condition may be assumed.
Here we proceed similarly by making assumptions on the convergence of the NLP solver.

\begin{myassumption}[NLP solver]\label{ass:nlp_solver}
	Consider a nonlinear program 
	\begin{align}\label{eq:nlp}
		\min_x f(x) \quad \textnormal{s.t.} \quad g(x) \geq 0,
	\end{align} 
	where $f:\R^n \to \R$ and $g:\R^n \to \R^{n_g}$ are twice continuously differentiable functions.
	Starting from an initial guess $x^0 \in \R^n$, an NLP solver returns either:
	\begin{enumerate}
%		\item a stationary point of \eqref{eq:nlp}; moreover, if $x^0$ is a feasible but not stationary point of \eqref{eq:nlp}, then it returns a point $x^*$ such that for some $\alpha >0$, it holds that:
%		\begin{align}\label{eq:nlp_sufficient_decrease2}
%			f(x^0)- f(x^*)  \geq  \alpha >0,
%		\end{align}
		\item a stationary point $x^*$ of \eqref{eq:nlp}; moreover, for every feasible but
		nonstationary point $\bar{x}$ of \eqref{eq:nlp}, there exist a neighborhood
		$\mathcal{N}_{\bar{x}}$ of $\bar{x}$ and a constant $\alpha_{\bar{x}}>0$ such that, whenever the NLP solver is initialized at a
		feasible but nonstationary point $x^0\in\mathcal{N}_{\bar{x}}$, it returns a stationary point $x^*$ satisfying
		\begin{align}\label{eq:nlp_sufficient_decrease}
			f(x^0)-f(x^*) \geq \alpha_{\bar{x}}>0.
		\end{align}
		\item a stationary point of the feasibility problem $\min_x \|g^-(x)\|$.
	\end{enumerate} 
\end{myassumption}
This assumption is justified as follows.
The NLP~\eqref{eq:nlp} covers all NLPs considered in this paper, including BNLPs~\eqref{eq:bnlp} and $\mathrm{Reg}(\tau)$ in \eqref{eq:mpec_scholtes}, for which MFCQ holds by Assumption~\ref{ass:cq} and Proposition~\ref{prop:mfcq_consequences}.
Under the MFCQ and additional assumptions (e.g., boundedness of the Hessian, second-order sufficient conditions), the two outcomes of  Assumption~\ref{ass:nlp_solver} are outcomes of several NLP algorithms, e.g.~\cite[Theorem 1]{Chin2003}, \cite[Theorem 7]{Fletcher2002b}) or~\cite{Waechter2006a}.
In practice, they are met by several globalized NLP solver implementations, e.g., \ipopt~\cite{Waechter2006}, \texttt{filterSQP}~\cite{Fletcher2002b}, \texttt{SNOPT}~\cite{Gill2005}, \texttt{KNITRO}~\cite{Byrd2006}, \texttt{MadNLP}~\cite{Shin2024}, \text{uno}~\cite{Vanaret2026}, to name a few.

%convergence of stabilized sequential quadratic programming~\cite{Wright2002,Wright1998} and interior-point methods~\cite{Vicente2002} to a local minimizer is guaranteed.
%The outcomes of Assumption~\ref{ass:nlp_solver} are met in practice by most globalized NLP solver implementations, e.g., \ipopt~\cite{Waechter2006}, \texttt{filterSQP}~\cite[Theorem 7]{Fletcher2002b}, \texttt{SNOPT}~\cite{Gill2005}, \texttt{KNITRO}~\cite{Byrd2006}, to name a few.

Note that the globalization strategies of all listed NLP solvers enforce sufficient decrease conditions on every step to guarantee convergence to a stationary point.
When starting from a feasible point, we assume that the step towards a stationary point satisfies a sufficient decrease condition, which we summarize in Assumption~\ref{ass:nlp_solver} by inequality~\eqref{eq:nlp_sufficient_decrease} (cf. \cite[Theorem 1]{Chin2003}, \cite[Theorem 7]{Fletcher2002b}, \cite{Waechter2006a}).

Even when starting from a feasible point, the NLP solver may fail to produce a point satisfying the globalization strategy (i.e., sufficiently feasible with improved objective). 
In such cases, a feasibility restoration phase is typically triggered~\cite{Waechter2006}, possibly returning a point with higher objective. 
We did not observe this in any of our numerical experiments and therefore do not treat it further. 
A possible remedy is to add the constraint {$f(x) \leq f(x^k) + \epsilon$}, with $\epsilon$ small, to the current BNLP.

Lastly, we comment on the convergence of an {\LPCC} solver.
We solve the \LPCC~\eqref{eq:lpec_full} via its MILP reformulation~\eqref{eq:lpec_milp}.
MILP solvers solve LPs or certify their infeasibility, e.g., within a branch-and-bound algorithm, which always converges~\cite{Wolsey1998}.
Both the interior-point method and the dual-simplex method (with suitable anti-cycling rules) either solve an LP or certify its infeasibility; see, e.g., \cite{Nocedal2006}.
Thus, in our setting, the LPECs can always be solved or certified infeasible.
We implement this using either the commercial solver {\gurobi}~\cite{Gurobi} or the open-source solver {\highs}~\cite{Huangfu2018}, both based on branch-and-bound methods.

%% SL START ==============================================================
%%SL If MFCQ holds at a point $x$, then any LP relaxation of~\eqref{eq:lpec_milp} will satisfy it as well. The MFCQ is sufficient for both the interior-point~\cite{what} and the dual-simplex method for LPs~\cite{What} to either solve the LP or certify its infeasibility. Both the interior-point~\cite{what} and the dual-simplex method for LPs~\cite{What} to either solve the LP or certify its infeasibility.

%% SL E N D ==============================================================

\subsection{Main convergence results}
We start with a useful proposition. 
A similar result is available in~\cite{Kirches2022,Leyffer2007}.
\begin{myproposition}\label{prop:non_b_stationarity}
	Let $x\in \Omega$ be a feasible point for the {\MPEC}~\eqref{eq:mpec}, {satisfying Assumption~\ref{ass:cq}}, and assume that $x$ is not B-stationary in the sense of Definition~\ref{def:b_stationarity}.
	Then there exist $\epsilon>0$, a step $s\in \R^n$ with $\|s\|_{\infty}=1$, and a partition $(\D_1(x),\D_2(x))$ of $\I_{00}(x)$ such that:
	\allowdisplaybreaks
	\begin{subequations}\label{eq:non_b_stationarity}
		\begin{align}
			&\nabla f(x)^\top s < -\epsilon,\\
			& \nabla c_i(x)^\top s \geq \epsilon, & i \in \A(x), \label{eq:non_b_stationarity_ineq}\\
			&s_{1,i} = 0, & i\in \I_{0+}(x),\\
			&s_{2,i} = 0, & i\in \I_{+0}(x),\\
			&s_{1,i} = 0,\ s_{2,i} \geq 0, & i\in \D_1(x),\\
			&s_{1,i} \geq 0,\ s_{2,i} = 0, & i\in \D_2(x).
		\end{align}
	\end{subequations}
\end{myproposition}

\textit{Proof.} 
If $x$ is not B-stationary, then there exists at least one BNLP for which $x$ is not stationary. 
Denote by $(\D_1(x),\D_2(x))$ the partition of $\I_{00}(x)$ associated with this BNLP.
By Proposition~\ref{prop:mfcq_consequences}, the usual MFCQ holds for this BNLP. 
From this we draw two conclusions.

First, since $x$ is not B-stationary for this BNLP, there exists a $d \in \mathcal{T}_{\Omega_{\mathrm{BNLP}}}(x)$ such that $\nabla f(x)^\top d< 0$. 
Moreover, because MFCQ holds, the tangent cone $\mathcal{T}_{\Omega_{\mathrm{BNLP}}}(x)$ coincides with the usual linearized feasible cone of the BNLP, rendering~\eqref{eq:b_stationariry} into an LP.

Second, considering the active sets for the BNLP at $x$, MFCQ implies the existence of a nonzero direction $s$ that points into the strict interior of the LP's~\eqref{eq:b_stationariry} feasible set, i.e., the following set is nonempty:
\begin{align*}
	S = \{\, &s \in \R^n \mid \nabla c_i(x)^\top s > 0, \ i \in \A(x); \\ 
	&s_{1,i} = 0, \ i \in \I_{0+}(x);\ 
	s_{2,i} = 0, \ i \in \I_{+0}(x); \\
	& s_{1,i} = 0,\ s_{2,i} > 0, \ i \in \D_1(x);\ 
	s_{1,i} > 0,\ s_{2,i} = 0, \ i \in \D_2(x)\, \}.
\end{align*}

Next, let $\bar s\in S$.
If $\nabla f(x)^\top \bar s<0$, set $s:=\bar s$. 
Otherwise, pick any $d\in\mathcal{T}_{\Omega_{\mathrm{BNLP}}}(x)$ with $\nabla f(x)^\top d<0$.
Because $\mathcal{T}_{\Omega_{\mathrm{BNLP}}}(x)$ is a cone then for every $\lambda>0$ we can define a new descent direction
$d_\lambda:=\lambda d \in \mathcal{T}_{\Omega_{\mathrm{BNLP}}}(x)$.
Next, define the vector $s(t):=\bar s+ t d_\lambda$.
Since $\bar s\in S$, there exists $t_0>0$ such that $s(t)\in S$ for all $t\in[0,t_0]$.
Moreover,
\[
\nabla f(x)^\top s(t)=\nabla f(x)^\top \bar s + t\,\nabla f(x)^\top d_\lambda
= \nabla f(x)^\top \bar s + t\lambda\,\nabla f(x)^\top d.
\]
Choose $\lambda$ so large that $0 < \frac{\nabla f(x)^\top \bar s}{-\lambda \nabla f(x)^\top d} < t_0$, and then take any
$t\in\big(\frac{\nabla f(x)^\top \bar s}{-\lambda \nabla f(x)^\top d},\,t_0\big]$. This yields $s(t)\in S$ and $\nabla f(x)^\top s(t)<0$.
Fix such $(\lambda,t)$ and set $s:=s(t)$.

Because inequalities are strict for $s(t)$, define
\(
\epsilon_1 := -\nabla f(x)^\top s, \
\epsilon_2 := \min_{i\in\A(x)} \nabla c_i(x)^\top s, \
\epsilon_3 := \min(
\min_{i\in\D_1(x)} s_{2,i},\!
\min_{i\in\D_2(x)} s_{1,i}
),
\)
which are all positive, and set $\bar\epsilon := \min\{\epsilon_1,\epsilon_2,\epsilon_3\}>0$. With $\hat s := s/\|s\|_\infty$ and $\epsilon := \bar\epsilon/\|s\|_\infty$, the vector $\hat s$ satisfies \eqref{eq:non_b_stationarity}.

%Combining these two facts yields the existence of $\epsilon > 0$ and $s$ with $\|s\|_{\infty} = 1$ such that \eqref{eq:non_b_stationarity} holds.
\qed

Next, we show that an \LPCC{} constructed at a point sufficiently close to a feasible point can ``predict'' a feasible BNLP.
In our implementation, such a point, and the corresponding {\LPCC} are found in Algorithm~\ref{alg:phase_i}.

\begin{mytheorem}\label{th:phase_i_feasibility}
	Let Assumptions~\ref{ass:compactness}-\ref{ass:nlp_solver} hold. 
	Consider $x^*\in \Omega$ a feasible point of the {\MPEC}~\eqref{eq:mpec}, and $\tilde{x} \in X$ (not necessarily an element of $\Omega$) such that $\|\tilde{x}-x^*\|\leq \delta$.
	Then, there exist positive constants $\kappa$, $\epsilon$ such that for a sufficiently small $\delta$ the following interval is nonempty
	\begin{align}\label{eq:feasiblity_lpec_rho}
		\max\Big\{ \frac{h_c(\tilde{x})}{\epsilon}, h_{\perp}(\tilde{x})\Big\} \leq \rho \leq \kappa,
	\end{align}
	and for any $\rho$ from this interval the $\mathrm{LPEC}(\tilde{x},\rho)$ in \eqref{eq:lpec_full} is feasible. 
	Moreover, for any feasible point $d$ of this {\LPCC}, it follows that $\mathrm{BNLP}(\I_1(\tilde{x}+d),\I_2(\tilde{x}+d))$ is feasible as well.
\end{mytheorem}

\textit{Proof.} 
Consider the case when $\tilde{x}\in \Omega$, i.e., $\delta =0$ and $h_c(\tilde{x})  = h_\perp(\tilde{x}) = 0$.
The $\mathrm{LPEC}(\tilde{x},\rho)$ is always feasible for this point and let $d$ be an element of its feasible set. 
Now set $\kappa = \bar{\rho}$, where $\bar{\rho}$ is defined in \eqref{eq:full_vs_reduced_lpec}.
Then according to Lemma~\ref{lem:reduced_lpec_feasible}, for all $\rho \in (0,\kappa)$, every $\BNLP(\I_1(\tilde{x}+d),\I_2(\tilde{x}+d))$ is feasible.

Next, we prove the assertion for $\tilde{x} \notin \Omega$.
The main idea of the proof is to consider a feasible $x^* \in \Omega$, and to show that Lemma~\ref{lem:reduced_lpec_feasible} is also true for points $\tilde{x}$ close enough to $x^*$. 
Denote by $\mathcal{B}(x,\rho) = \{ d \in \R^n \mid \|d-x\|_{\infty} \leq \rho\}$ the $\ell_{\infty}$-ball centered at $x \in \R^n$ with the radius $\rho>0$.

First, we show that {\LPCCs} defined at points close to $x^*$ are also feasible under suitable conditions.
Let $x^*\in \Omega$ be feasible but not B-stationary.
It follows from Proposition~\ref{prop:non_b_stationarity} and the continuity of the functions $f$, and $c$, that there exists a neighborhood $\mathcal{N}^* = \mathcal{B}(x^*,\kappa_1) \cap X$, such that for all $\tilde{x}\in \mathcal{N}^*$ it holds that: 
\begin{align}\label{eq:lpec_continity_f_c}
		\nabla f(\tilde{x})^\top s < -\epsilon, \nabla c_i(\tilde{x})^\top s \geq \epsilon,\  i \in \A(x^*).
\end{align}
Next, we prove that every $\LPEC(\tilde{x},\rho)$ for all $\tilde{x} \in \mathcal{N}^*$ and a suitable $\rho$ is feasible.
We start first with the general inequality constraints and use an argument similar to \cite[Lemma 4]{Chin2003}.
Take the step $d=\rho s$, which satisfies $\|d\|\leq \rho$, and consider the linearization of the inequality constraints active at $x^*$:
\begin{align}\label{eq:feasible_lp_ineq_1}
	&c_i(\tilde{x}) + \rho \nabla c_i(\tilde{x})^\top s \geq -h_c(\tilde{x}) + \rho \epsilon, &i \in \A(x^*).
\end{align}
The first term is lower bounded by the infeasibility of the inequality constraints, the second uses the fact that $\|s\|_{\infty} = 1$ and \eqref{eq:lpec_continity_f_c}.
For the inequality constraints inactive at $x^*$, it follows from the continuity of $c$ that there exists positive constants $\gamma_1$ and $\gamma_2$ independent of $\tilde{x}$ such that 
$c_i(\tilde{x})  \geq \gamma_1$ and $\nabla c_i(\tilde{x})^\top s \geq -\gamma_2$. 
Thus for the inactive constraints, we have
\begin{align}\label{eq:feasible_lp_ineq_2}
	&c_i(\tilde{x}) + \rho \nabla c_i(\tilde{x})^\top s \geq \gamma_1 - \rho\gamma_2, &i \notin \A(x^*).
\end{align}
By requiring that the right-hand-side of~\eqref{eq:feasible_lp_ineq_1} and \eqref{eq:feasible_lp_ineq_2} are nonnegative, we find that the linearized inequality constraints are feasible if $\rho$ satisfies:
\begin{align}\label{eq:rho_bound_general_ineq}
	\frac{h_c(\tilde{x})}{\epsilon} \leq \rho \leq \frac{\gamma_1}{\gamma_2}.
\end{align}

We now consider the complementarity constraints:
%\begin{align*}
	 $0 \leq \tilde{x}_1  + d_1 \perp \tilde{x_2} + d_2 \geq 0$.
%\end{align*}
It can be seen that to make these constraints always feasible the trust region radius $\rho$ must be within the following bounds:
\begin{align*}
	\max(|\min(\tilde{x}_1,\tilde{x}_2))|) \leq \rho \leq \max(\max(|\tilde{x}_1|,|\tilde{x}_2|)).
\end{align*}
The lower bound is obtained if the complementarity constraint is satisfied by always setting the smaller component of $(\tilde{x}_{1,i},\tilde{x}_{2,i})$ to zero via the step $d$.
Similarly, the upper bound is obtained by setting the larger component of $(\tilde{x}_{1,i},\tilde{x}_{2,i})$ to zero via the step $d$.
The lower bound is equal to the complementarity residual $h_{\perp}(\tilde{x})$. 
Moreover, from Lipschitz continuity it follows that $h_{\perp}(\tilde{x}) \leq \delta$ and 
$h_{c}(\tilde{x}) \leq \delta$.
Because, $\tilde{x} \in X$, there exists a constant $\kappa_2 >0$ independent of $\tilde{x}$ such that $\|\tilde{x}\| \leq \kappa_2$. 
Therefore, for the feasibility of the complementarity constraints we require that:
\begin{align}\label{eq:rho_bound_comp}
%	{h_\perp(\tilde{x})}\leq \rho \leq \kappa_2
	 \delta \leq \rho \leq \kappa_2
\end{align}

Therefore, by combining \eqref{eq:rho_bound_general_ineq} and \eqref{eq:rho_bound_comp} we have that 
for all $\tilde{x} \in \mathcal{N}^*$ and $\rho$ for which holds:
\begin{align}\label{eq:lpec_rho_fesible}
	\max\Big\{\frac{\delta}{\epsilon},\delta\Big\} \leq \rho \leq  \min\Big\{\kappa_2,\frac{\gamma_1}{\gamma_2}\Big\}.
\end{align}
Since there always exists a $\varepsilon  > 0$ for every $x^*\in \Omega$ (see proof of Proposition~\ref{prop:non_b_stationarity}), the term $\lim\limits_{\delta \to 0 }\frac{\varepsilon}{\delta} =0 $ is well-defined.
Therefore, the $\LPEC(\tilde{x},\rho)$ is feasible, since for a sufficiently small $\delta$ the interval \eqref{eq:lpec_rho_fesible} is nonempty. 

In the final step of the proof, we make sure that: 
\begin{align}\label{eq:lpec_tr_relation}
	\mathcal{B}(\tilde{x},\rho) \subset \mathcal{B}(x^*,\bar{\rho}).
\end{align}
This can be achieved by computing:
\begin{align*}
\| x-x^*\|_{\infty} \leq \| x-\tilde{x} \|_{\infty}+\| \tilde{x}-x^* \|_{\infty} \leq \delta  + \rho \leq \bar{\rho}, 
\end{align*}
yielding the upper bound on $\rho \leq \bar{\rho} - \delta$. 
Combining this with \eqref{eq:lpec_rho_fesible} we require that 
\begin{align}\label{eq:lpec_rho_final}
	\max\Big\{\frac{\delta}{\epsilon},\delta\Big\} \leq \rho  \leq  \min\Big\{\kappa_2,\frac{\gamma_1}{\gamma_2}, \bar{\rho} - \delta\Big\}.
\end{align}
For sufficiently small $\delta$ this interval is nonempty and we define
\sloppy $\kappa = \min\Big\{\kappa_2,\frac{\gamma_1}{\gamma_2}, \bar{\rho} - \delta\Big\}$.

Furthermore, we conclude the following about the feasible $\LPEC(\tilde{x},\rho)$.
Let $\tilde{d}$ be any feasible point of $\LPEC(\tilde{x},\rho)$.
Due to \eqref{eq:lpec_tr_relation}, all points $\tilde{x}+\tilde{d}$ that can be reached from $\tilde{x}$ and that satisfy the complementarity constraints, can also be reached via $x^* + d$, where $d$ is a feasible point of $\LPEC(x^*,\bar{\rho})$.
Because Lemma \ref{lem:reduced_lpec_feasible} holds for $\LPEC(x^*,\bar{\rho})$, we conclude that also all $\BNLP(\I_1(\tilde{x}+\tilde{d}),\I_2(\tilde{x}+\tilde{d}))$ are feasible.
Fig.~\ref{fig:feasibility_lpec_proof} illustrates the arguments.
This concludes the proof. \qed

\begin{figure}[t]
	\centering
	\subfloat[Illustration for some $i\in \I_{+0}(x^*)$ ]{\includegraphics[width=0.43\textwidth]{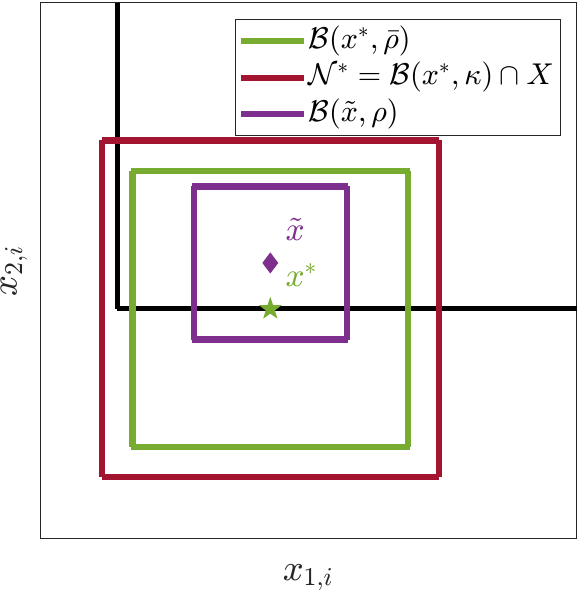}\label{fig:feasibility_lpec_proof_1}}
	\hspace{1.2cm}
	\subfloat[Illustration for some $j\in \I_{00}(x^*)$ ]{\includegraphics[width=0.43\textwidth]{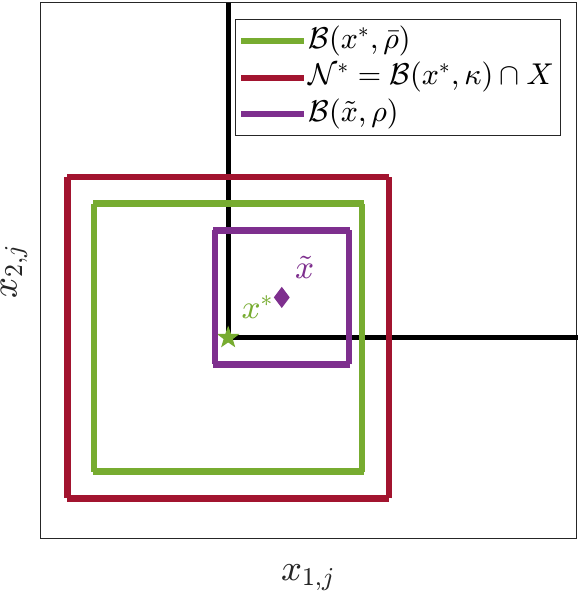}\label{fig:feasibility_lpec_proof_2}}
	\caption{Illustration of the argument used in the proof of Theorem~\ref{th:phase_i_feasibility}.}
	\label{fig:feasibility_lpec_proof}
\end{figure}
\begin{myremark}
Theorem~\ref{th:phase_i_feasibility} ensures that any algorithm that produces a point sufficiently close to the MPEC feasible set can be used in Phase~I.
Subsequently, finding a feasible point of the LPEC constructed at this point helps identify a feasible BNLP.
In Section~\ref{sec:algorithm_phase_i}, we introduced a specific Phase~I algorithm and, in the following corollary, we specialize Theorem~\ref{th:phase_i_feasibility} to this algorithm.
\end{myremark}

\begin{mycorollary}\label{col:phase_i_feasibility}
	Let Assumptions~\ref{ass:compactness}-\ref{ass:nlp_solver} hold. 
	Let $x^*(\tau)$ be a stationary point of $\mathrm{Reg}(\tau)$~in~\eqref{eq:mpec_scholtes} for some $\tau>0$. 
	There exist positive constants $\bar{\tau}$, $\kappa$ such that for any	$\tau \in \left(0,\bar{\tau}\right]$ and $\rho \in [\sqrt{\tau},\kappa]$ the $\mathrm{LPEC}(x^*(\tau),\rho)$ in \eqref{eq:lpec_full} is feasible. 
	Moreover, for any feasible point $d$ of this {\LPCC}, it follows that $\mathrm{BNLP}(\I_1(x^*(\tau)+d),\I_2(x^*(\tau)+d))$ is feasible as well.
\end{mycorollary}
\textit{Proof.}
First, regard the case when $x^*(\tau) \in \Omega$, for some $\bar{\tau} \geq \tau > 0$. 
The $\mathrm{LPEC}(x^*(\tau),\rho)$ is always feasible for this point and let $d$ be an element of its feasible set. 
Now for $\kappa = \bar{\rho}$ and all $\rho \in (0,\kappa)$, according to Lemma~\ref{lem:reduced_lpec_feasible} all $\mathrm{BNLP}(\I_1(x^*(\tau)+d),\I_2(x^*(\tau)+d))$ are feasible.

Second, for $x^*(\tau) \notin \Omega$ set $\tilde{x} = x^*(\bar{\tau})$ and apply Theorem~\ref{th:phase_i_feasibility}. 
Since $x^*(\bar{\tau})$ is a stationary point of $\mathrm{Reg}(\bar{\tau})$,
it holds that $c(x^*(\bar{\tau}))\geq 0$, i.e., $h_c(\tilde{x})  = 0$.
Moreover, from $x_{1,i} x_{2,i} \leq \tau, \
\forall i \in \{1,\ldots,m\}$, it follows that $h(\tilde{x}) \leq \sqrt{\bar{\tau}}=\delta$.
Hence, for a sufficiently small $\delta=\bar{\tau}$ and every $0 <\tau \leq \bar{\tau}$ the interval in \eqref{eq:lpec_rho_final} is nonempty and we can apply Theorem~\ref{th:phase_i_feasibility}. 
\qed 

Remarkably, it follows from the last two results, that the {\LPCC} in Phase I does not even need to be solved to optimality. 
It is sufficient to find a feasible point of the {\LPCC} to identify a feasible BNLP.
Before we present the main convergence result, we prove that the inner loop in Algorithm~\ref{alg:phase_ii} terminates in a finite number of iterations.

\begin{myproposition}\label{prop:finite_termination}
	Let Assumptions~\ref{ass:compactness}-\ref{ass:nlp_solver}  hold and let $x^k$ be feasible for the	MPEC~\eqref{eq:mpec}. 
	Then the inner loop of Algorithm~\ref{alg:phase_ii} terminates after finitely many iterations. 
	Furthermore, if $d^{k,l}\neq 0$ solves $\operatorname{LPEC}(x^k,\rho^{k,l})$ with	$\rho^{k,l}<\bar\rho$, then $x^k$ is not stationary for any
	$\BNLP(\I_1(x^k+d^{k,l}),\I_2(x^k+d^{k,l}))$ with $(\I_1(x^k+d^{k,l}),\I_2(x^k+d^{k,l})) \in \mathcal{P}(x^k+d^{k,l})$.
\end{myproposition}
\textit{Proof.}
If $d^{k,l}=0$ at some inner iteration, Algorithm~\ref{alg:phase_ii} terminates. 
Otherwise, each rejected step multiplies the trust-region radius by $\gamma_L\in(0,1)$. 
Hence, after finitely many rejections,	there exists an index $N\geq l$ such that $\rho^{k,N}=(\gamma_L)^{N-l}\rho^{k,l}<\bar\rho$ and $d^{k,N}\neq0$.

For $\rho^{k,N}<\bar\rho$, due to Lemma~\ref{lem:full_vs_reduced_lpec} we can focus on the reduced LPEC~\eqref{eq:lpec_reduced}.
In this case, because $\| d^{k,N} \|_{\infty} < \bar  \rho$  is sufficiently small, from the proof of Lemma~\ref{lem:reduced_lpec_feasible} we have
\[
\mathcal{P}(x^k+d^{k,N}) \subseteq \mathcal{P}(x^k),
\]
so every BNLP that can be selected at $x^k+d^{k,N}$ can also be selected at $x^k$, hence $x^k$ is feasible for all these BNLPs.

Since no termination occurred before index $N$, it holds that $\nabla f(x^k)^\top d^{k,N}<0$ and $d^{k,N}\neq 0$.
By \eqref{eq:mpec_feasible_sets} the vector $d^{k,N}$ solves all branch LPs of the reduced {\LPCC},
and by \eqref{eq:mpec_lin_cones} each such branch LP corresponds to an LP of the form~\eqref{eq:b_stationariry} for the BNLPs.
Therefore $d^{k,N}\neq 0$ is a common minimizer of all these LPs with $\nabla f(x^k)^\top d^{k,N}<0$, which implies that $x^k$ is feasible but not stationary for any BNLP indexed by $(\I_1(x^k+d^{k,N}),\I_2(x^k+d^{k,N})) \in \mathcal{P}(x^k+d^{k,N})$.

By Assumption~\ref{ass:nlp_solver}, the NLP solver applied to any such BNLP produces a point $x^{k,N}$ with
$f(x^{k,N})<f(x^k)$. 
Hence an acceptance occurs and the inner loop terminates.
Since $N$ is finite, the inner loop terminates after finitely many iterations.
\qed
An immediate consequence of this proposition is that we can terminate an LPEC solver early at a nonzero feasible point $d \neq 0$.
\begin{mycorollary}\label{cor:finite_termination}
Under the assumptions of Proposition~\ref{prop:finite_termination}, 
if $x^k$ is not B-stationary, then for sufficiently small $\rho^{k,N} < \bar{\rho}$, 
any feasible point $d \neq 0$ of LPEC$(x^k,\rho^{k,N})$ satisfies 
$\nabla f(x^k)^\top d < 0$, and $x^k$ is feasible but not stationary for any BNLP 
indexed by $(\mathcal{I}_1(x^k+d),\mathcal{I}_2(x^k+d)) \in \mathcal{P}(x^k+d)$.
\end{mycorollary}
\textit{Proof.} 
It follows from the proof of the previous proposition that, for sufficiently large $N$, the solution of the LPEC, $d^{k,N}\neq0$, also solves every branch LP. 
By item (3) of Proposition~\ref{prop:mfcq_consequences}, and for $\rho^{k,N}<\bar{\rho}$, the feasible set of a branch LP is a convex cone intersected with the trust-region constraints. 
Because it is a cone, from the feasibility of $\hat{d}=0$ and $d^{k,N}$ it follows that any $d= t d^{k,N}$ with $t>0$ is also feasible (within the trust region), and $\nabla f(x^k)^\top (t d^{k,N}) = \nabla f(x^k)^\top d < 0$. 
Moreover, because $d = t\ d^{k,N}$ it holds that  $(\mathcal I_1(x^k+d),\mathcal I_2(x^k+d))$ coincides with $(\mathcal I_1(x^k+d^{k,N}),\mathcal I_2(x^k+d^{k,N}))$.
\qed

Finally, we prove the main convergence result.
The first outcome is inherited from the NLP solver failing to find a feasible point (cf. Assumption~\ref{ass:nlp_solver}~(2)).
The second outcome is the most common when the algorithm succeeds, and it is guaranteed to occur if all BNLP stationary points are isolated.
The third outcome occurs if the algorithm generates an infinite sequence, and its proof follows similar lines to those of \cite[Theorem 1, f-monotonic case]{Chin2003} and \cite[Theorem 7, case (C)]{Fletcher2002b}.

\begin{mytheorem}\label{th:main_convergence}
{Regard the {\MPEC}~\eqref{eq:mpec} and suppose that Assumptions~\ref{ass:compactness}-\ref{ass:nlp_solver} hold.}
Then {\MPECopt} given in Algorithm~\ref{alg:phase_ii}, with a Phase I from Section~\ref{sec:algorithm_phase_i} has one of the following mutually exclusive outcomes:
\begin{enumerate}
	\item {In Phase I, a stationary point of the problem of minimizing a measure of constraint infeasibility is returned.}
	\item It computes a point $x^*$ that is B-stationary for the MPEC~\eqref{eq:mpec} after solving finitely many BNLPs and LPECs.
	\item {It generates an infinite sequence of iterates $\{x^k\}$ with decreasing objective values, whose accumulation point is B-stationary.}
\end{enumerate}
{If, in addition, every stationary point $\hat{x}$ of each BNLP is isolated, i.e., there exists a neighborhood $\mathcal{N}_{\hat{x}}$ of $\hat{x}$ containing no other stationary points, then outcome 3 cannot occur.}
\end{mytheorem}

\textit{Proof.}
First, we discuss Phase I and outcome one. 
Consider Algorithm~\ref{alg:phase_i}. 
It follows from Lemma~\ref{lem:phase_i_infeasibility} that if \eqref{eq:mpec_scholtes} is infeasible, then so is the {\MPEC} ~\eqref{eq:mpec}.
Infeasibility is sufficient for outcome two of Assumption~\ref{ass:nlp_solver}.
Moreover, if the NLP solver fails to solve \eqref{eq:mpec_scholtes} for some $k$, then we have again outcome two of Assumption~\ref{ass:nlp_solver}.
Now consider the case when a stationary point of  all \eqref{eq:mpec_scholtes} is found. 
If $x^*(\tau^k)$ is feasible, Phase I terminates successfully with a feasible point.
If not, then it follows from Corollary~\ref{col:phase_i_feasibility} that there exists a finite $k$ such that the LPEC$(x^*(\tau^k),\rho)$ is feasible, and that the corresponding BNLP is feasible. 
If the NLP solver finds a stationary point of this BNLP, then Phase I terminates successfully with a feasible point.
In conclusion, Algorithm~\ref{alg:phase_i} either finds a feasible point $x^0$ or terminates with outcome one in a finite number of iterations.
This completes the proof for the first outcome.

Given a feasible point $x^0$, we discuss the outcomes of Algorithm~\ref{alg:phase_ii}.
It follows from Proposition~\ref{prop:finite_termination} that the inner loop terminates in a finite number of iterations.
Therefore, the outer loop produces a sequence of points $\{x^k\} \subset \Omega$ with decreasing objective values $\{f^k\}$.
Consider some iterate $x^k$, which is a stationary point of some $\mathrm{BNLP}(\I_1(x^k),\I_2(x^k))$, and a $\rho^{k,l}>0$. 
If $d^{k,l}=0$ is a minimizer of $\mathrm{LPEC}(x^k,\rho^{k,l})$, then a B-stationary point is found and the algorithm terminates in a finite number of iterations, which is outcome two.

We prove the occurrence of outcome three.
Since $\{x^k\}\subset\Omega\subset X$, $f$ is continuous, and $X$ is compact, the sequence $\{f(x^k)\}$ is bounded from below.
By monotonicity, $\{f(x^k)\}$ converges to some $f^\infty\in\R$. 
Consequently, $f(x^k)-f(x^{k+1})\to0$.
Moreover, compactness of $X$ implies that $\{x^k\}$ has an accumulation point $x^\infty\in X$. Let $\{x^{k_j}\}$ be a subsequence such that
$x^{k_j}\to x^\infty$.

Suppose, for contradiction, that $x^\infty$ is not B-stationary. 
By Assumption~\ref{ass:nlp_solver}, applied locally at $x^\infty$, there exist a neighborhood $\mathcal{N}_{x^\infty}$ and a constant
$\alpha_{x^\infty}>0$ such that $f(x^k)-f(x^{k+1})\geq\alpha_{x^\infty} > 0$ whenever $x^k\in\mathcal{N}_{x^\infty}$, which already
contradicts the fact that $\{f(x^k)\}$ is bounded from below.
Since $x^{k_j}\to x^\infty$, there exists $J\in\mathbb{N}$ such that $x^{k_j}\in\mathcal{N}_{x^\infty}$ for all $j\geq J$.
Hence, the contradiction holds also for the subsequence $\{x^{k_j}\}$.
Therefore, $x^\infty$ is B-stationary.

Finally, we consider the case where all stationary points of the BNLPs are isolated.  
Since $X$ is compact and there is a finite number of possible BNLPs, it follows that the total number of stationary points is finite.  
Consequently, the sequence $\{f(x^k)\}$ is finite, and rules out the third outcome.
Algorithm~\ref{alg:phase_ii} generates a sequence of strict decreases 
\(
\Delta f^k = f(x^k) - f(x^{k+1}) > 0 .
\)  
Let $x^K$ denote the last point in this sequence, attained after a finite number of steps $K$.  
Suppose that for the optimal solutions of the LPEC$(x^k,\rho^{k,l})$ it holds that $d^{k,l} \neq 0$ for all $l$ and all $k < K$.  
Since no problem $\mathrm{BNLP}(\I_1(x^K),\I_2(x^K))$ can produce a better objective value, it follows from \eqref{eq:mpec_feasible_sets} that $x^K$ is stationary for all these BNLPs.  
Hence, $x^K$ is a B-stationary point of the {\MPEC}~\eqref{eq:mpec} and we have outcome two.
This concludes the proof.
\qed

{We note that the assumption of isolated stationary points for finite termination is not overly restrictive. 
Under MPEC-MFCQ in Ass.~\ref{ass:cq}, the usual MFCQ holds for the associated BNLPs. 
Moreover, if we additionally assume generalized second-order sufficient conditions at the BNLP stationary points, these points are isolated local minimizers~\cite[Theorem 6.9]{Forsgren2002}.}

\section{Numerical results}\label{sec:numerical_results}
In this section, we extensively compare several {\MPECopt} {variants with regularization-based methods, as well as with solving the MINLP and NLP reformulations of {\MPECs}} on two distinct nonlinear {\MPEC} benchmarks.

\subsection{Implementation}
In our implementation and experiments, we treat {\MPECs} in a more general form than \eqref{eq:mpec}:
\begin{subequations}\label{eq:mpec_implementation}
	\begin{align}
		\underset{x\in \R^{n}}{\mathrm{min}} \;  \quad &f(x)\\
		\textnormal{s.t.} \quad 
		&c^{\mathrm{lb}} \leq c(x) \leq c^{\mathrm{ub}}, \\
		&x^{\mathrm{lb}} \leq x \leq x^{\mathrm{ub}},\\
		&0 \leq g(x)\perp h(x) \geq 0,\label{eq:mpec_implementation_comp}
	\end{align}
\end{subequations}
where $c^{\mathrm{lb}}, c^{\mathrm{ub}} \in \R^{n_c}$ are lower and upper bounds of $c(x)$ with $ c^{\mathrm{ub}}_i \geq c^{\mathrm{lb}}_i, \ \forall i \in \{1,\ldots,n_c\}$. 
Equality constraints are imposed by setting $ c^{\mathrm{ub}}_i = c^{\mathrm{lb}}_i$.
The vectors $x^{\mathrm{lb}}, x^{\mathrm{ub}} \in \R^{n}$ are lower and upper bounds on $x$. 
If there is no lower or upper bound on some $x_i$ or $c_i(x)$, then we set the respective value on $-\texttt{inf}$ or $\texttt{inf}$.
The functions $g,h:\R^{n} \to \R^{m}$ are assumed twice continuously differentiable.

To solve the NLPs arising in {\MPECopt} we use the \texttt{IPOPT} solver~\cite{Waechter2006}. 
For better performance, we change some of the default options by setting the tolerance to $\texttt{tol} = 10^{-12}$ (default is $10^{-8}$), 
$\texttt{acceptable\_tol} = 10^{-9}$ (default is $10^{-6}$), and
$\texttt{warm\_start\_init\_point} = \texttt{'yes'}$.
Moreover, we change the barrier parameter strategy by setting:
\texttt{mu\_strategy} = \texttt{'adaptive'} and
\texttt{mu\_oracle} = \texttt{'quality-function'}, cf.~\cite{Nocedal2009}.
In all experiments, we use \texttt{MA27}~\cite{Duff1982} from the \texttt{HSL} library~\cite{HSL} as the linear solver in \ipopt.
The same settings are used for all methods and NLP solutions in the following experiments.

All problems are formulated in the symbolic framework of \texttt{CasADi}~\cite{Andersson2019} via its \textsc{Matlab} interface.
{\casadi} provides the first and second-order derivatives for all problem functions via its automatic differentiation and has an interface for calling \ipopt.

We formulate the LPECs as MILPs using~\eqref{eq:lpec_milp} and solve them with the commercial solver {\gurobi}~\cite{Gurobi} or the open-source solver {\highs}~\cite{Huangfu2018}.
{In our experiments, the MILP reformulation was significantly faster than regularization- and penalty-based methods, cf. Appendix~\ref{sec:numerics_macmpec_lpec}.}
For the MILPs, we set the maximum number of nodes in the branch and bound algorithm to 5000 and the maximum computation time to 5 minutes. 
The tolerance for the LP solves in both MILP solvers is set to $10^{-9}$ (default is $10^{-6}$).

For comparisons, we also implement the regularization and penalty-based methods.
Thereby, an adaptation of \ref{eq:mpec_scholtes}, $\ell_{1}$ and $\ell_{\infty}$~\cite{Kim2020,Ralph2004,Nurkanovic2024b} to the form of ~\eqref{eq:mpec_implementation} are solved in a homotopy loop with a decreasing parameter $\tau$.
We initialize $\tau^0 = 1$ and update it using the rule $\tau^k = 0.1 \tau^{k-1}$.
The loop terminates if either: 
the NLP solution is feasible and the complementarity tolerance $h_{\perp}(x)$ is below $10^{-9}$; 
NLP is locally infeasible; 
or a maximum number of $15$ homotopy iterations is reached.
A homotopy loop is declared successful if it finds a stationary point.
Additionally, we compare {\MPECopt} to solving the {\MPEC} directly as an NLP in~\eqref{eq:mpec_nlp} with {\ipopt}, and solving its MINLP reformulation~\eqref{eq:mpec_minlp} with a nonlinear branch-and-bound method implemented in \texttt{Bonmin}~\cite{Bonami2008}.

We mention the remaining parameter choices in our implementation.
The parameters in Algorithm~\ref{alg:phase_i} are $\kappa = 0.1$, $\tau^0=1$, ${\rho}^I = 10^{-1}$, and $\texttt{tol}_h = 10^{-10}$.
In the implementation of Algorithm~\ref{alg:phase_ii}, in line~\ref{alg:phase_ii_check_B}, a point is declared to be B-stationary if $\|d^{k,l}\|_{\infty} \leq \texttt{tol}_{\mathrm{B}}$.
However, since there is also the trust region constraint $\|d\| \leq \rho^{k,l}$ we should avoid having $\rho^{k,l}< \texttt{tol}_{\mathrm{B}}$. 
Since the inner loop terminates within a finite number of iterations (cf. Proposition~\ref{prop:finite_termination}), this situation can be avoided by choosing appropriate parameters of the algorithm; we pick $\texttt{tol}_{\mathrm{B}} = 10^{-8}$, $\rho^{k,0} = 10^{-3}$, $N^{\mathrm{in}} = 6$ and $\gamma^{\mathrm{L}} = 0.1$.
We also set a maximum number of outer iterations to $N^{\mathrm{out}} = 25$.
If no B-stationary point is found, the algorithm terminates after reaching the maximum number of iterations.

Our implementation of {\MPECopt} is open-source and available at \url{github.com/nosnoc/mpecopt}. 
It contains code to reproduce the two benchmark results below, and several tutorial examples.

\subsection{Results on MacMPEC benchmark}\label{sec:numerics_macmpec}
The {\MacMPEC} problem set is available in the form of a package of AMPL  \texttt{.mod} and \texttt{.dat} files from \url{wiki.mcs.anl.gov/leyffer/index.php/MacMPEC}.
We consider in total 191 problems. 
Results are reported using Dolan and Mor\'e performance profiles~\cite{Dolan2002}, here referred to as relative performance profiles, and analogous absolute timing profiles, for which no problem-wise scaling is applied and the horizontal axis shows wall time in seconds.

\begin{figure}[h]
	\centering
	\subfloat[Absolute performance profiles.\label{fig:macmpec_general_absolute}]{\includegraphics[width=0.5\textwidth]{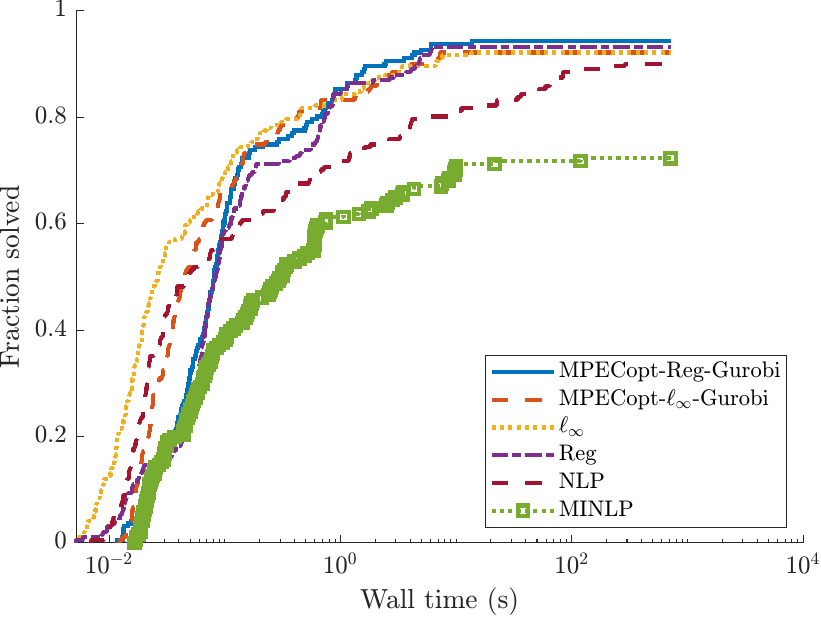}}
	\subfloat[Relative performance profiles.\label{fig:macmpec_general_relative}]{\includegraphics[width=0.5\textwidth]{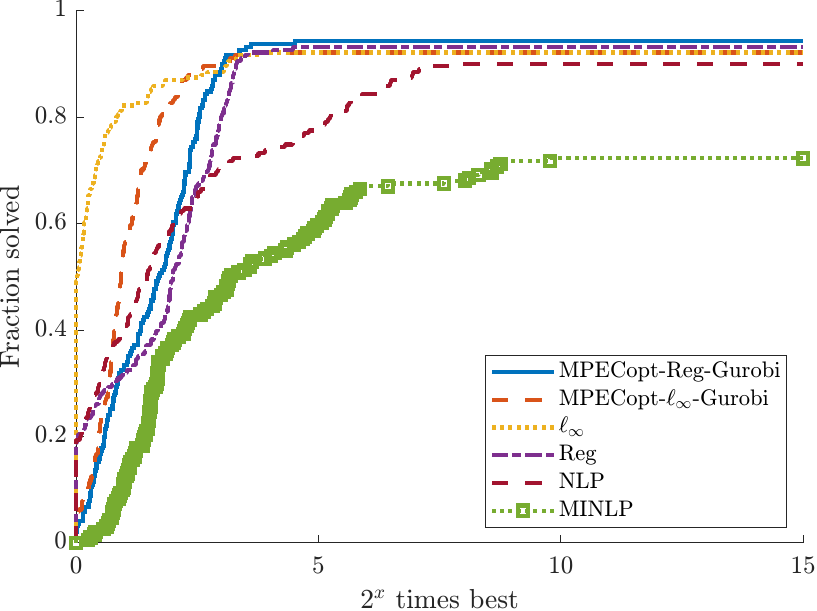}}\\
	\caption{Evaluating different {\MPEC} solution methods on the {\MacMPEC} test set in terms of finding a B-stationary point.}
	\label{fig:macmpec_general}
\end{figure}

We compare {\MPECopt} to other mentioned approaches.
For Phase I of {\MPECopt}, we use Algorithm~\ref{alg:phase_i} with {\gurobi} as {\LPCC} solvers, as it had the highest success compared to other variants we tested.
In the experiments, we compare in total six algorithms:
\begin{enumerate}
	\item \MPECopt, with regularization-based Phase I, and \texttt{Gurobi} as LPEC solver -- \texttt{MPECopt-Reg-Gurobi}.
	\item \MPECopt, with $\ell_\infty$-penalty-based Phase I, and \texttt{Gurobi} as LPEC solver -- \texttt{MPECopt-$\ell_\infty$-Gurobi}.
%	\item \MPECopt, with regularization-based Phase I, and HiGHS as LPEC solver -- \texttt{MPECopt-HiGHS}.
	\item Regularization method that solves \eqref{eq:mpec_scholtes} in a homotopy loop -- \texttt{Reg}.
	\item $\ell_\infty$-penalty method that solve an $\ell_{\infty}$ reformulation of \eqref{eq:mpec}~\cite{Kim2020} in a homotopy loop -- \texttt{$\ell_\infty$}.
	\item Directly solving the MPEC as an NLP~\eqref{eq:mpec_nlp} with \ipopt -- \texttt{NLP}.
	\item Solving the MINLP reformulation \eqref{eq:mpec_minlp} with a nonlinear branch and bound in \texttt{Bonmin} -- \texttt{MINLP}.
\end{enumerate}

We also tried the $\ell_{1}$ penalty method, which has a similar computation time as $\ell_{\infty}$, but a slightly lower success rate.
Moreover, using {\highs} instead of {\gurobi} led to slightly slower method and marginally lower success rate. 
Therefore, for better readability, we omit these variants from the plots.
Since all methods use \texttt{IPOPT} as the NLP solver, we can directly compare the computation times.
Fig.~\ref{fig:macmpec_general} reports the absolute and relative performance profiles for this experiment.

Note that, if {\MPECopt} is successful, it always returns a B-stationary point, whereas for the other methods, we can only know that it is a B-stationary point if it is also S-stationary.
To verify B-stationarity of the solutions computed by all other methods, we proceed as described in Appendix~\ref{sec:stat_points}.
We do not add the computational cost of verifying B-stationarity to the overall computation time of the other methods.

We observe that \texttt{MPECopt-Reg-Gurobi} achieves the highest overall success rate of 94.24\%.
It is both faster and more successful than its counterpart \texttt{Reg} (93.2\%).
The $\ell_{\infty}$-penalty method is the fastest overall and faster than \texttt{MPECopt-Reg-Gurobi}, but it attains a lower success rate of 92.14\%.
Meanwhile, \texttt{MPECopt-$\ell_\infty$-Gurobi} is only slightly slower than \texttt{$\ell_\infty$} and achieves the same success rate.
As expected, the plain NLP and MINLP reformulations are less competitive than the MPEC-tailored methods, with success rates of 90.00\% and 72.52\%, respectively.
The main reasons for failure are timeouts and unsuccessful restoration phases.
Additionally, the MINLP method frequently terminates at feasible but non–B-stationary points (39 cases).

The {\MPECopt} method inherits the efficacy of the standard methods used in Phase I but can often be even faster, as it requires solving fewer NLPs.
This is because, already for moderate values of $\tau$, the LPEC can identify a feasible BNLP in Phase I.
Consequently, {\MPECopt} needs to solve fewer NLPs than \texttt{Reg} and \texttt{$\ell_\infty$} in most cases.
The branch-and-bound MINLP method solves the largest number of NLPs, which also explains its slow computation times.

{In practice, only a few {\LPCCs} need to be solved.
For instance, on the {\MacMPEC} benchmark, the number of {\LPCCs} solved ranged from 1 to 6.
In many cases, only two are required: one in Phase I to identify a feasible BNLP, and one in Phase II to verify B-stationarity.}

We have shown in Section~\ref{sec:convergence_theory} that, unless we need to verify B-stationarity, the {\LPCCs} do not need to be solved to global optimality; finding a feasible point is sufficient.
Moreover, Proposition~\ref{lem:milp_via_lp} shows that if the solution is S-stationary, it can be verified by solving a single LP.
We examine these theoretical results by plotting histograms of the number of nodes in a MILP call (corresponding to the number of LPs solved) in both Phase I and Phase II.
Fig.~\ref{fig:macmpec_lpec_full} shows the histogram when the MILP is solved to global optimality, while Fig.~\ref{fig:macmpec_lpec_early} shows the case when it is terminated early, at a feasible point.
In both cases, Phase II only needed to solve LPs, and in Phase I we observe a significant reduction in the number of LPs to be solved in several cases.
In fact, with early termination, on {\MacMPEC} we only need to solve LPs.
This illustrates our theoretical results: early termination of the MILP method reduces the number of LPs that must be solved.
This did not significantly affect the overall computation time, as the largest computational burden on the {\MacMPEC} problem set came from solving the NLP subproblems.

\begin{figure}[t]
	\centering
	\subfloat[Solving to optimality.\label{fig:macmpec_lpec_full}]{\includegraphics[width=0.49\textwidth]{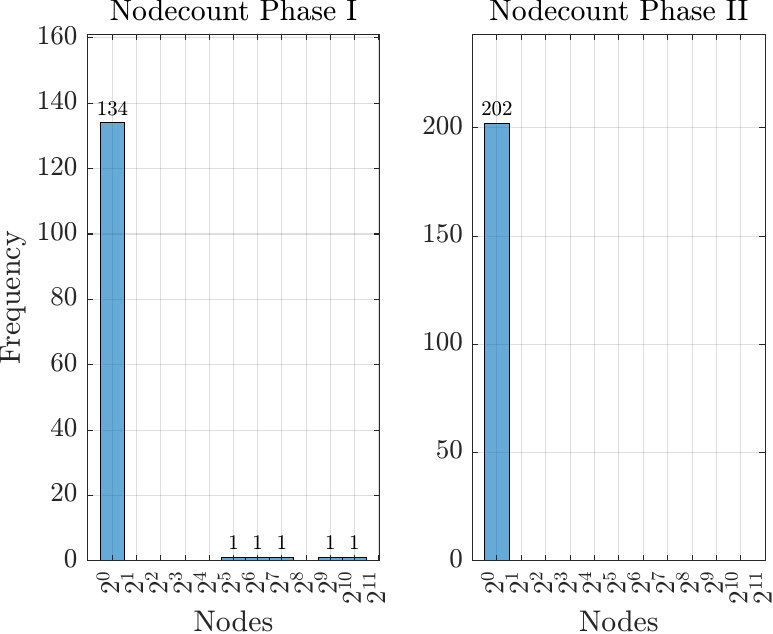}}
	\subfloat[Early termination.\label{fig:macmpec_lpec_early}]{\includegraphics[width=0.49\textwidth]{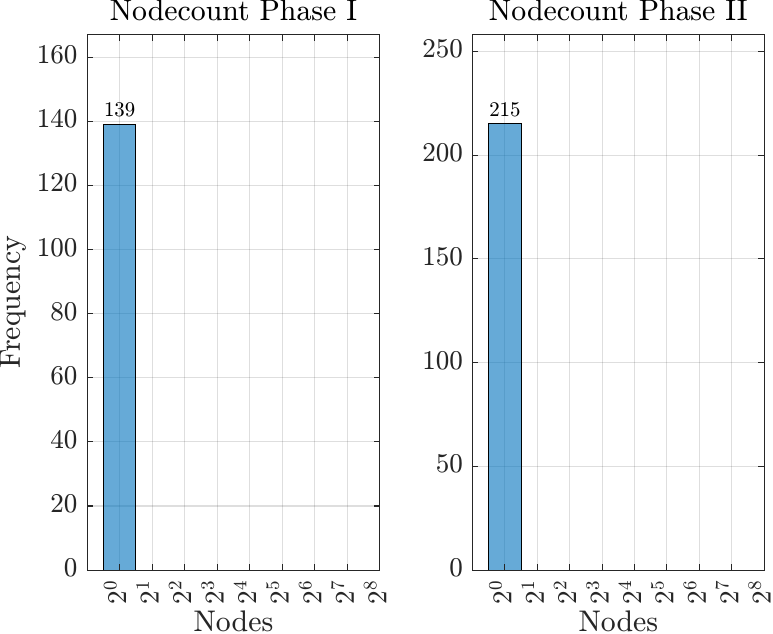}}\\
	\caption{Number of nodes in the MILP solver in Phase I and II. The node counter on the horizontal axis is on logarithmic scale.}
	\label{fig:macmpec_lpec}
\end{figure}

In summary, all methods -- sometimes with minor tuning -- can successfully compute B-stationary points on most {\MacMPEC} problems, although only {\MPECopt} directly certifies this.
All methods are quite fast on easy problems without biactive constraints, where often a single NLP solve is sufficient.
In fact, more than two-thirds of the solutions found by any algorithm correspond to problems with an empty biactive set.
To provide a broader perspective, we present another benchmark with larger problems in the next section.

\begin{figure}[t]
	\centering
	\subfloat[Absolute performance profiles.\label{fig:nonlinear_general_absolute}]{\includegraphics[width=0.5\textwidth]{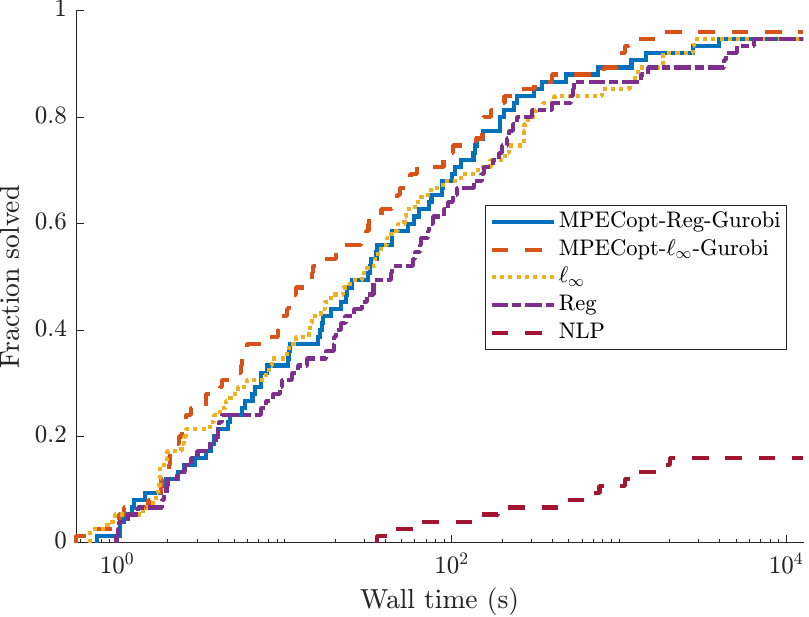}}
	\subfloat[Relative performance profiles.\label{fig:nonlinear_general_relative}]{\includegraphics[width=0.5\textwidth]{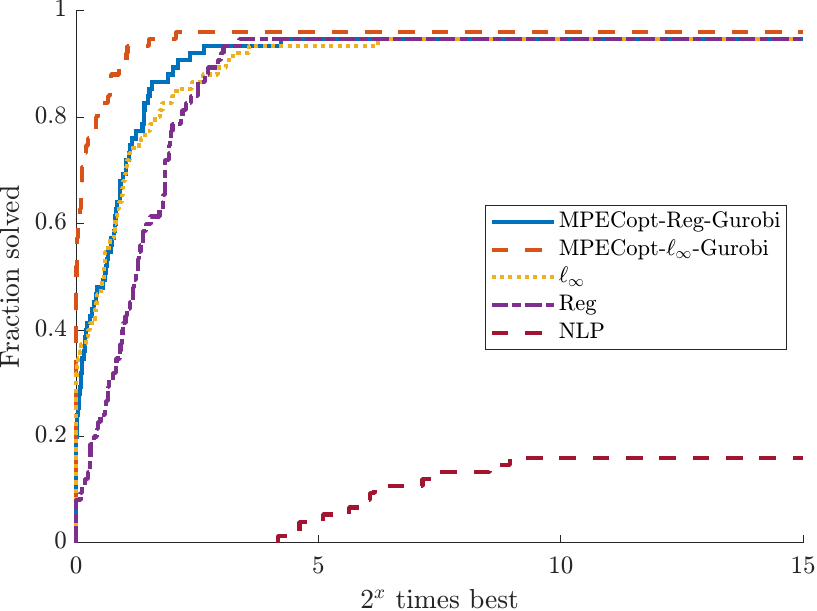}}\\
	\caption{Evaluating different {\MPEC} solution methods on the synthetic test set in terms of finding a stationary point.}
	\label{fig:nonlinear_general}
\end{figure}

\subsection{Synthetic nonlinear {\MPEC} benchmark}\label{sec:numerics_synthetic}
Many {\MacMPEC} problems are small or easily solved by all examined algorithms.
To further compare \MPECopt, we use a synthetic nonlinear MPEC benchmark described in Appendix~\ref{sec:synthetic_benchmark}.
The benchmark consists of 75 significantly larger problems, ranging from 200 complementarity constraints and 500 variables in the smallest instances to 4000 complementarity constraints and 10000 variables in the largest.
These problems are specifically designed to violate the MPEC-LICQ.
In this experiment, we exclude the MINLP from the comparison, as it was already not competitive on the easier problem set.
Fig.~\ref{fig:nonlinear_general} shows the relative and absolute performance profiles of the methods we compared.

\begin{figure}[t]
	\centering
	\subfloat[Solving to optimality.\label{fig:nonlinear_lpec_full}]{\includegraphics[width=0.49\textwidth]{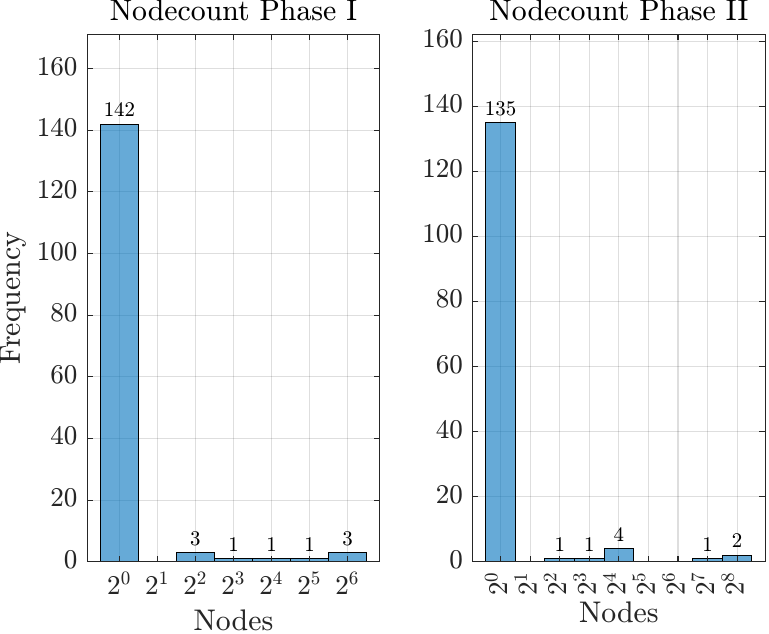}}
	\subfloat[Early termination.\label{fig:nonlinear_lpec_early}]{\includegraphics[width=0.49\textwidth]{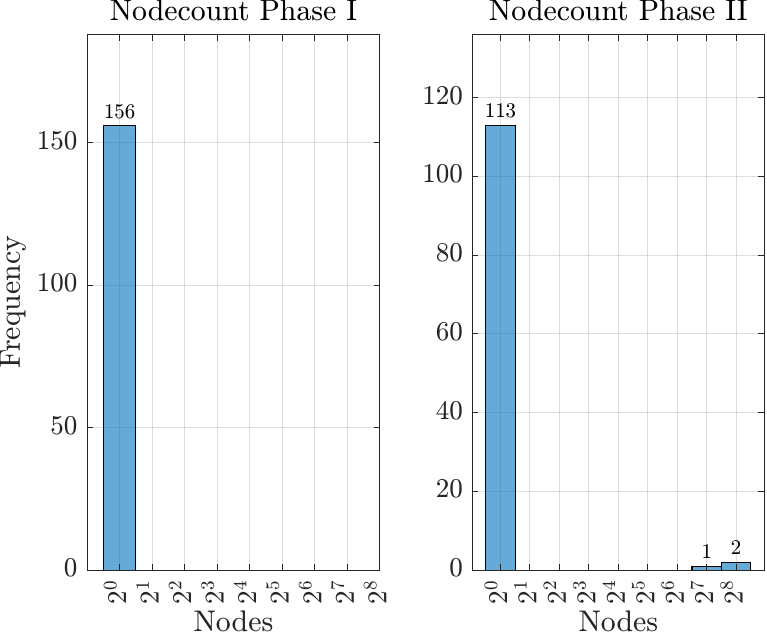}}\\
	\caption{Number of nodes in the MILP solver in Phase I and II. The node counter on the horizontal axis is on logarithmic scale.}
	\label{fig:nonlinear_lpec}
\end{figure}

Similar to the results on {\MacMPEC} (Fig.~\ref{fig:macmpec_general}), we see that the {\MPECopt} variants are the fastest and most successful algorithms even on the large problem instances.
For the same reasons as in the previous experiments, they require fewer total NLP solves.
Moreover, on this set, the differences in computation times between \texttt{MPECopt-Reg-Gurobi} (96.16\% success rate), \texttt{MPECopt-$\ell_{\infty}$-Gurobi} (94.67\% success rate), and \texttt{Pen}-$\ell_{\infty}$ (94.67\% success rate) are more pronounced, whereas \texttt{Reg} (94.67\% success rate) needs to solve more NLPs and is therefore slower.
The lower success rates of the regularization- and penalty-based methods compared to their {\MPECopt} counterparts are due to either reaching the maximum number of iterations or hitting the time limit.
In other cases, the methods failed by converging to a point of local infeasibility, i.e., when the restoration phase failed.
Directly solving the NLP is by far the least successful approach (16.00\% success rate), and in most cases it timed out after reaching the 45-minute limit.
This further justifies the development of tailored methods for solving \MPECs.
As shown in Fig.~\ref{fig:nonlinear_general_absolute}, the MPEC-tailored methods compute solutions considerably faster and rarely run into timeouts.

Lastly, we examine the effects of early termination in the LPEC solvers in Fig.~\ref{fig:nonlinear_lpec}.
In Fig.~\ref{fig:nonlinear_lpec_full}, where the MILP (and thus the LPEC) is solved to global optimality, we see that in both phases there are instances requiring a significant number of LPs.
By contrast, Fig.~\ref{fig:nonlinear_lpec_early} shows the effects of early termination in Phase I (due to Corollary~\ref{col:phase_i_feasibility}) and Phase II (due to Corollary~\ref{cor:finite_termination}), where the number of LPs is drastically reduced.
In particular, in Phase I a single LP was always sufficient.
In Phase II, we also observe a reduction.
We terminate the MILP solver if a feasible point with an objective below $-0.1 \| \nabla f(x^k)\| \rho^{k,l}$ is found.
With more careful tuning of this factor, the number of LPs in Phase II can often be reduced further.

In conclusion, we can see that {\MPECopt} is competitive with other existing approaches. 
Unlike all other methods, it always provides a certificate of B-stationarity.
Moreover, the experiments demonstrate that solving an {\LPCC} was never a computational bottleneck, even for the larger problem instances.
Significant savings can be achieved through early termination, in line with our theoretical predictions.

\section{Conclusions and summary}\label{sec:conclusion} 
This paper introduces {\MPECopt}, a two-phase method for computing B-stationary points of mathematical programs with complementarity constraints (MPECs).
The first phase identifies a feasible point of the MPEC or certifies local infeasibility.
The second phase iterates over feasible points to find a B-stationary point by solving a finite sequence of branch nonlinear programs (BNLPs) and linear programs with complementarity constraints (LPECs).
We perform a detailed convergence analysis of the algorithm.

We also provide an open-source implementation of {\MPECopt}, available at \url{github.com/nosnoc/mpecopt}.
The implementation leverages robust and mature NLP and MILP solvers.
We tested {\MPECopt} on two numerical benchmarks: the {\MacMPEC} set and a new synthetic nonlinear MPEC test set containing problems with up to four thousand complementarity constraints.
The numerical results confirm the efficacy of the proposed method for finding B-stationary points.
{\MPECopt} is generally faster and more reliable than existing methods, requiring fewer NLP solves than homotopy-based approaches.

In general, LPECs are challenging combinatorial optimization problems that can be formulated as mixed-integer linear programs (MILPs).
Interestingly, a combinatorial number of LPs only needs to be solved if the limit point is B-stationary but not S-stationary.
We prove that in all other cases, the LPEC solver can be terminated early or that solving a single LP is sufficient.
In our numerical experiments on the two benchmarks, even for large-scale problem instances, it was often sufficient to solve only the root node of the MILP reformulation of the LPEC, which corresponds to a single LP solve.

An intriguing observation from the benchmarks is that the limit point of regularization- or penalty-based methods is in most cases B-stationary, even when it is not S-stationary.
To our knowledge, current convergence results for these methods do not explain this phenomenon.
However, verifying B-stationarity typically requires solving additional nonlinear programs and LPECs, whereas {\MPECopt} inherently certifies it, since solving LPECs is fully integrated into the algorithm.

Finally, in this paper we assume that the BNLPs are solved to full convergence.
A natural direction for future work is to study inexact BNLP solutions while maintaining robust and fast convergence.

\section*{Acknowledgments}
A.N. acknowledges the generous support of this research by the DFG via Research Unit FOR 2401, project 424107692, 504452366 (SPP 2364), and 525018088, by BMWK via 03EI4057A and 03EN3054B, and by the EU via ELO-X 953348.
Moreover, A.N. thanks Anton Pozharskiy for translating \texttt{MacMPEC} to \texttt{CasADi}, and for many useful remarks on the {\texttt{\MPECopt}} and benchmark implementations.
This work was also supported by the Applied Mathematics activity within the U.S. Department of Energy, Office of Science, Advanced Scientific Computing Research, under contract number DE-AC02-06CH11357.
\section*{Data availability}
All data necessary to reproduce the results in this paper are publicly available in the GitHub repository \url{https://github.com/nosnoc/mpecopt}. 
This repository contains the code, experiment scripts, and data used to generate the numerical results reported in the manuscript.
\section*{Competing interests}
The authors declare no competing interests.

\appendix
\section{Proofs of lemmas from Section~\ref{sec:algorithm}}\label{sec:lemma_proofs}
\textit{Proof of Lemma~\ref{lem:full_vs_reduced_lpec}}:
%\textit{Proof.}
Denote by $\Omega^{\mathrm{R}}$ and $\Omega^{\mathrm{F}}$ the feasible sets of the reduced {\LPCC}~\eqref{eq:lpec_reduced} and full {\LPCC}~\eqref{eq:lpec_full}, excluding the constraints $\|d  \| \leq \rho$, respectively.
Compared to $\Omega^{\mathrm{R}}$, the set $\Omega^{\mathrm{F}}$ is additionally restricted by the linear inequalities for all inactive constraints, i.e.,  $c_i(x^k) + \nabla c_i(x^k)^\top d \geq 0, i \in \bar{\A}^k$.
Recall that the constraints $0 \leq x_{1,i}^k + d_{1,i} \perp x_{2,i}^k + d_{2,i} \geq 0$ define the set $\{ d \mid x_{1,i}^k + d_{1,i} \geq 0, \ x_{2,i}^k + d_{2,i} =0  \} \cup  
\{ d \mid x_{1,i}^k + d_{1,i}   = 0,\ x_{2,i}^k + d_{2,i} \geq 0 \}$.
Moreover, $\Omega^{\mathrm{F}}$ contains all branches excluded for the strictly active complementarity indices $i \in \I_{+0}^k \cup \I_{0+}^k$ in $\Omega^{\mathrm{R}}$.
Therefore, the two sets are related as follows:

\begin{align}\label{eq:lpec_sets}
	\begin{split}
		\Omega^{\mathrm{F}} = \Omega^{\mathrm{R}} 
		& \bigcap_{i\in \bar{\A}^k} \{ d \mid  c_i(x^k) + \nabla c_i(x^k)^\top d \geq 0 \}\\
		&\bigcup_{i\in \I_{0+}^k} \{ d \mid x_{1,i}^k + d_{1,i} \geq 0, \ x_{2,i}^k + d_{2,i} =0  \}\\ 
		&\bigcup_{i\in \I_{+0}^k} \{ d \mid x_{1,i}^k + d_{1,i}   = 0,\ x_{2,i}^k + d_{2,i} \geq 0 \}.
	\end{split}
\end{align}

To ensure that the two {\LPCCs} have the same set of solutions, we need to restrict $\|d\|_{\infty} \leq \rho $ such that $\Omega^{\mathrm{F}} = \Omega^{\mathrm{R}}$.

Let us first consider the inequality constraints. 
For all $i \in \bar{\A}^k$, we have $c_i(x^k) > 0$, and we want 
\(
c_i(x^k) + \nabla c_i(x^k)^\top d > 0.
\) 
If $\nabla c_i(x^k) = 0$, this inequality holds trivially for any $d$ (and thus any $\rho$). 
Otherwise, if $\nabla c_i(x^k) \neq 0$, we can apply the Cauchy--Schwarz inequality
\(
| \nabla c_i(x^k)^\top d | \leq \| \nabla c_i(x^k) \| \, \| d \|.
\)
Hence, the linearized inequality remains inactive as long as 
\(
\| d \| < \frac{c_i(x^k)}{\| \nabla c_i(x^k) \|}.
\) 
To ensure this holds for all relevant $i$, we define
\(
\bar{\rho}_1 = \min_{i \in \bar{\A}^k, \, \nabla c_i(x^k) \neq 0} \frac{c_i(x^k)}{\| \nabla c_i(x^k) \|}.
\) 
Next, we compute upper bounds for $\rho$ to exclude the union terms in \eqref{eq:lpec_sets}.

Now observe that, for example, for any $i\in\I_{+0}^k$ it holds that $x_{1,i}^k > 0$ and $x_{2,i}=0$.
If $\rho$ is chosen such that $d_{1,i} <\rho < x_{1,i}^k$ then  $x_{1,i}^k+d_{1,i}>0$ and 
$0 \leq  x_{1,i}^k + d_{1,i} \perp d_{2,i} \geq 0 $ can be only satisfied for $d_{2,{i}} = 0$. 
Thus, the branch $\{ d \mid x_{1,i}^k + d_{1,i} = 0,\ x_{2,i}^k + d_{2,i} \geq 0 \}$ cannot be reached, i.e., it is infeasible. 
Similar reasoning holds for any $i\in\I_{0,+}^k$.
Therefore, by picking 
$\rho < \bar{\rho}_2 = \min \{
\{ x_{1,i} \mid i \in \I_{+0}^k \}  
\cup 
\{ x_{2,i} \mid i \in \I_{0+}^k \}\}$, 
none of the terms in the unions in \eqref{eq:lpec_sets} is feasible.
The equation \eqref{eq:full_vs_reduced_lpec} follows from writing $\bar{\rho} = \min(\bar{\rho}_1,\bar{\rho}_2)$ explicitly.
Thus, for $\rho < \bar{\rho}$, we have $\Omega^{\mathrm{F}} \cap \{ \|d\|_\infty \leq \rho \} = \Omega^{R} \cap \{ \|d\|_\infty \leq \rho \}$.
As the two {\LPCCs} have the same objective, their sets of local minimizers coincide.
\qed

\noindent\textit{Proof of Lemma~\ref{lem:reduced_lpec_feasible}}:
It follows from Lemma~\ref{lem:full_vs_reduced_lpec} that for $\rho \in (0,\bar{\rho})$ both the full and reduced {\LPCC} have the same feasible set.  
Thus, we focus solely on the reduced {\LPCC}.  
Denote by $\mathcal{P}^k$ the set of all possible partitions $(\mathcal{I}_1^k,\mathcal{I}_2^k)$ at $x^k \in \Omega$ (see Eq.~\eqref{eq:active_set_partition} for definitions).   
The MPEC-MFCQ implies the MPEC-ACQ, hence the MPEC-linearized feasible cone and the standard linearized feasible cones of the BNLPs satisfy the following~\cite[Lemma 3.1]{Flegel2005a}:
\begin{align}\label{eq:mpec_lin_cones}
	\begin{split}
		\mathcal{F}_{\Omega}(x^k) = \bigcup_{(\mathcal{I}_1^k,\mathcal{I}_2^k)\in \mathcal{P}^k} \{ d \in \R^{n} \mid   
		&c_i(x^k)+ \nabla c_i(x^k)^\top d \geq 0, \ \forall i \in \A^k,\\  
		&x_{1,i}^k + d_{1,i} = 0,\ x_{2,i}^k + d_{2,i} \geq 0, \ \forall i \in  \I_{1}^k, \\  
		&x_{1,i}^k + d_{1,i} \geq 0,\ x_{2,i}^k + d_{2,i} = 0, \ \forall i \in \I_{2}^k \}.  
	\end{split}
\end{align}  
Moreover, the MPEC--MFCQ implies the usual MFCQ for each BNLP, which in turn implies the usual ACQ for each BNLP. 
Therefore, each set on the right-hand side of the above relation is equal to the tangent cone of the respective BNLP$(\I_1^k,\I_2^k)$~\cite[Corollary 3.1]{Flegel2005a}.  
Note that $x^k$ is feasible for all these BNLPs.  
The intersection of $\mathcal{F}_{\Omega}(x^k)$ with $\| d \|_\infty \leq \rho < \bar{\rho}$ recovers the feasible set of \eqref{eq:lpec_reduced}.  
Furthermore, it holds that $\mathcal{P}(x^k+d) \subseteq \mathcal{P}(x^k)$, because biactive indices can become strict (i.e., $\I_{00}(x^k+d)\subseteq \I_{00}(x^k)$), but strict indices must remain strict for $\| d \|_\infty \leq \rho < \bar{\rho}$ (i.e., $\I_{0+}(x^k) \subseteq \I_{0+}(x^k+d)$ and $\I_{+0}(x^k) \subseteq \I_{+0}(x^k+d)$). 
Any feasible point $d$ of this {\LPCC} determines a particular branch of the {\LPCC}, defined by a partition $(\I_1^*(x^k+d), \I_2^*(x^k+d)) \in \mathcal{P}(x^k+d) \subseteq \mathcal{P}^k$.  
This leads to the BNLP$(\I_1^*(x^k+d), \I_2^*(x^k+d))$, for which $x^k$ is a feasible point.  
This completes the proof.  
\qed

\section{Identifying stationary points}\label{sec:stat_points}
To determine the type of stationary point from Definition~\ref{def:mpec_stationarity} we must solve the TNLP~\eqref{eq:tnlp}, which requires a correct identification of the active sets $\I_{+0},\I_{0+}$ and $\I_{00}$.
{\MPECopt} is an active set method where the complementarity constraints are usually satisfied exactly or close to machine precision, making it easy to find the active sets for the TNLP.

On the other hand, \texttt{Reg}, and \texttt{Pen-$\ell_1$} rarely have a complementarity residual $h_{\perp}(x)$ much smaller than the prescribed tolerance of $10^{-9}$, making an active set identification more difficult.
We take an iterative procedure, and identify the active set via:
\begin{subequations}\label{eq:epsilon_active_set}
	\begin{align}
		&\I_{0+}^{\epsilon}(x) = \{ i \mid x_{1,i} \leq \epsilon, x_{2,i} > \epsilon\},\\
		&\I_{+0}^{\epsilon}(x) = \{ i \mid x_{1,i} > \epsilon, x_{2,i} \leq \epsilon\},\\
		&\I_{00}^{\epsilon}(x) = \{ i \mid x_{1,i} \leq \epsilon, x_{2,i} \leq \epsilon\}, 
	\end{align}
\end{subequations}
with $\epsilon = 10^{-6}$.
To minimize the possibility of wrong conclusions, we implement several safeguards.
We denote the solution computed via \texttt{Reg} or \texttt{Pen-$\ell_1$} by $x^*$, and the corresponding TNLP solution by $x^{\mathrm{tnlp}}$. 
We declare the active set to be correctly identified if $\|x^*-x^{\mathrm{tnlp}}\| \leq 10^{-7}$.
If this is not satisfied we reduce $\epsilon$ by a factor of $10$ and repeat the procedure.
Once $\|x^*-x^{\mathrm{tnlp}}\| \leq 10^{-7}$ is true or a maximum number of 5 attempts is exceeded, we compute the multipliers in the TNLP and classify the stationary point.
To verify B-stationarity, we take several steps.
First, we use the point $x^{\mathrm{tnlp}}$ and construct an $\LPEC(x^{\mathrm{tnlp}},\rho)$ with $\rho = 10^{-3}$.
If the objective of this LPEC is below $\texttt{tol}_{\mathrm{B}} =10^{-8}$, then we declare $x^*$ to be B-stationary.
If this is not the case, we take up to three more LPEC solve attempts while reducing $\rho$ by a factor of 10. 
The minimal value is $\rho = 10^{-6}$, which is the minimal allowed trust region radius in \MPECopt.
This delivers a valid certificate in most cases.
To have more robustness to wrong active set identification, if this does not verify B-stationarity, we repeat the same but now with $x^*$ as linearization point $\LPEC(x^*,\rho)$.
Finally, if this also fails, we look at the objective $f(x^*)$ and the objective computed by {\MPECopt} for the same problem, which is denoted by $f(\hat{x})$. 
If the relative difference satisfies $\frac{|f(x^*)-f(\hat{x})|}{|f(\hat{x})|+10^{-16}} < 10^{-3}$, then we declare $x^*$ to be B-stationary, if $\hat{x}$ is B-stationary.
If the difference between the objective values is larger, and the LPECs find a descent direction, then the point $x^*$ is declared to be not B-stationary.

\section{Comparison of different {\LPCC} solvers}\label{sec:numerics_macmpec_lpec}
We compare the cumulative NLP and LPEC computation times in every {\MPECopt} solver call on every {\MacMPEC} problem. 
In all experiments, we use the regularization-based Phase I described in Algorithm~\ref{alg:phase_i}, but with different LPEC solvers in both phases.
We compare solving the LPEC as MILP~\eqref{eq:lpec_milp} and via regularization-based methods.
The {\highs} solver performs significantly better if the trust region in the LPEC is smaller, thus we set  $\rho_I=10^{-3}$ in Phase I, while {\gurobi} is less sensitive to $\rho_I$ values.
We consider four variants: 
\texttt{Gurobi-MILP} ($\rho_I=10^{-1}$), 
\texttt{HiGHS-MILP} ($\rho_I=10^{-3}$), 
\texttt{Reg-MPEC}, 
and \texttt{$\ell_{\infty}$-MPEC}.

We solve all {\MacMPEC} problems with four different LPEC solution strategies, comparing total cumulative computation times for solving all LPECs in each instance.
In 34 of 191 problems, the MPEC is already solved or locally infeasible after the first NLP, so we exclude these (effective LPEC time is zero).
The remaining 157 problems are solved by all approaches except \texttt{HiGHS-MILP}, which fails on \texttt{siouxfls} and \texttt{siouxfls1} due to unsolved {\LPCCs}.

Fig.~\ref{fig:lpec_vs_nlp} shows the results, where the diagonal represents equal NLP and LPEC solving times.
For MILP methods, most points fall below the diagonal, indicating that BNLPs dominate computation time rather than LPECs.
For regularization methods, points cluster around/above the diagonal, showing computational burden is split almost equally between NLPs and LPECs.
Specifically, cumulative NLP time exceeds LPEC time in 0\% of cases for \texttt{Gurobi-MILP}, 20\% for \texttt{HiGHS-MILP}, 87.26\% for \texttt{Reg-MPEC}, and 85.99\% for \texttt{$\ell_{\infty}$-MPEC}.
MILP solvers benefit from smaller initial trust region radii in Phase I, but smaller $\rho_I$ requires solving more Phase I NLPs to find feasible LPECs, shifting computational load toward NLPs.

\begin{figure}[t]
	\centering
	{\includegraphics[width=0.99\textwidth]{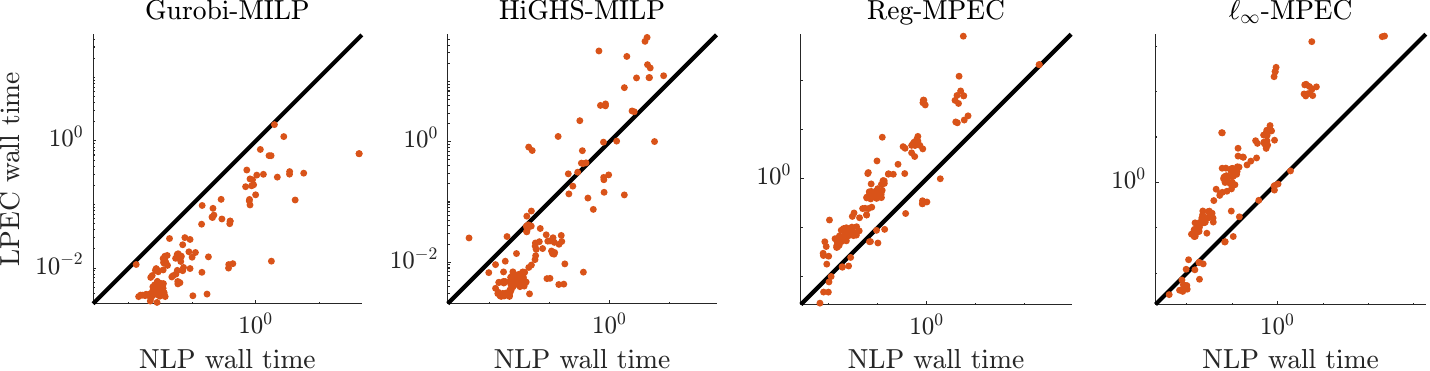}}
	\caption{Comparison of total NLP and LPEC computation times on the {\MacMPEC} benchmark.}
	\label{fig:lpec_vs_nlp}
\end{figure}

{Somewhat surprisingly, the global MILP algorithms significantly outperform the local regularization and penalty-based methods.}
The MILP-based methods are up to two orders of magnitude faster than regularization-based methods for LPECs on the {\MacMPEC} test set.
This difference can be explained by several observations.
The MILP methods solve the LPEC to global optimality (sufficient for verifying B-stationarity), whereas the others compute only a stationary point (which may not be a local minimum or B-stationary point). 
Consequently, LPECs solved by regularization-based methods are more often locally infeasible in Phase I or require solving more LPECs until a B-stationary point is verified in Phase II. 
This confirms the arguments in Section~\ref{sec:solving_lpecs} {on efficiently solving LPECs using MILP reformulations.}
Moreover, when the solution is S-stationary, solving the MILP requires solving a single LP, complying with Proposition~\ref{lem:milp_via_lp}.
This experiment suggests that solving {\LPCCs} with MILP methods as in {\MacMPEC} can be done very efficiently in practice, both with commercial (\gurobi) and open-source (\highs) solvers.

\section{Generating a synthetic MPEC benchmark}\label{sec:synthetic_benchmark}
In this section, we describe how we generate a random set of MPECs that are used in the numerical experiments of Section~\ref{sec:numerics_synthetic}.
The problem-generating function and code to reproduce the benchmark results is available at \url{https://github.com/nosnoc/benchmarks/nonlinear_mpec_benchmark}.

A summary of how the random MPECs are generated is given in Algorithm~\ref{alg:random_mpecs}.
Additionally, we describe how the problems are generated.
First, we generate random LPECs as described in~\cite{Hu2008}. 
They have the following form:
\begin{subequations}\label{eq:lpec_random}
	\begin{align}
		\underset{x_0 \in \R^{p}, x_1 \in \R^{m} } \min \;  \quad & g_0^\top x_0 + g_1^\top x_1 \label{eq:lpec_random_objective}\\
		\textnormal{s.t.} \quad 
		&A x_0 + B x_1 \geq a\\
		&0 \leq x_1 \perp Mx_0 + Nx_1 + q  \geq 0 \label{eq:lpec_random_cc},
	\end{align}
\end{subequations}
where $g_0\in \R^{p}$, $g_1 \in \R^m$ are the cost vectors. 
The inequality constraint data has the following dimensions
$A\in \R^{n_\mathrm{ineq} \times p}$, $B\in \R^{n_\mathrm{ineq} \times m}$, $a \in \R^{n_{\mathrm{ineq}}}$. 
The dimensions of the data in left hand side of \eqref{eq:lpec_random_cc} are $M \in \R^{m \times p}$, $N\in \R^{m \times m}$ and $q\in \R^{m}$.
A pseudocode for the generation of these LPECs can be found in~\cite[p. 22]{JaraMoroni2020}. 
We use the same description as a starting point but allow more flexibility in changing the problem data and dimensions.

To generate a nonlinear MPEC, we make the following modification to the LPEC~\eqref{eq:lpec_random}. 
We replace the objective in \eqref{eq:lpec_random_objective}  by $f(x)=f^{\mathrm{nl}}(x_0,x_1) + \alpha  (g_0^\top x_0 + g_1^\top x_1)$, where $f^{\mathrm{nl}}(x_0,x_1) $ is a nonlinear objective test function from unconstrained problems from the \texttt{CUTEst} test set~\cite{Gould2013}.
Let $c^{\mathrm{ineq}}(x) = A x_0 + B x_1 - a \geq 0$ define the inequality constraints.
To introduce more nonlinearity, we consider a subset of inequality constraint indices, i.e., $\I^\mathrm{ineq} \subseteq \{1,\ldots, n_\mathrm{ineq} \}$. 
Next, all inequality constraints $c_i^{\mathrm{ineq}}(x)\geq 0$ with $i\in \I^\mathrm{ineq}$ are replaced by $c_i^{\mathrm{ineq}}(x) + c_i^{\mathrm{ineq}}(x)^2 + c_i^{\mathrm{ineq}}(x)^4 \geq 0$.
Moreover, we duplicated a subset of the inequality constraints (both linear and nonlinear), with indices $i \in \I^\mathrm{d} \subseteq \{1,\ldots, n_\mathrm{ineq}\}$.
This makes the MPEC more degenerate as it increases the probability of violating MPEC-LICQ if inequality constraints are active. Next, we define equality constraints by introducing a slack variable $x_2\in \R^m$ for the left hand side in ~\eqref{eq:lpec_random_cc}, which reads as $c^{\mathrm{eq}}(x) = x_2 - Mx_0 + Nx_1 + q$. 
Analogously to the inequality constraints, we replace some constraints $c^{\mathrm{eq}}_i(x) =0$ with
$c^{\mathrm{eq}}_i(x) + c^{\mathrm{eq}}_i(x)^2 + c^{\mathrm{eq}}_i(x)^4 = 0, i \in \I^\mathrm{eq} \subseteq \{1,\ldots, m\}$. 
Moreover, we add upper bounds $x^{\mathrm{ub}}_0, x^{\mathrm{ub}}_1$ on $x_0$ and $x_1$. 

The full pseudocode for generating a single problem instance is given in Algorithm~\ref{alg:random_mpecs}.
Thereby, $\mathcal{N}(\mu,\sigma)$ denotes the normal distribution with $\mu$ mean and $\sigma$ standard deviation. 
The elements of $\mathcal{U}_n(a,b)$ are $n$ dimensional vectors whose elements are drawn from a uniform distribution on the interval $[a,b]$.
Similarly, the elements of $\mathcal{U}_{n\times k}(a,b,s)$ are $n \times k$ matrices whose elements take values from a uniform distribution on $[a,b]$. 
The scalar $s\in \left(0,1\right]$ is the density of this random matrix, which is defined as the ratio between the number of nonzero elements and the total number of elements of this matrix.
We chose the density of the random matrices to be up to 10\%.
For larger problems with more than 3000 variables, this results in quite dense matrices, and generating the computational graphs for the first and second-order derivatives can take several hours.
To ensure that the benchmark can be solved in under a week, we set the density of a random matrix $A \in \R^{n\times m}$ to $s = \frac{n}{n m}$ for problems with more than 3000 variables.

\begin{algorithm}
	\caption{Random MPEC generator}
	\label{alg:random_mpecs}
	\begin{algorithmic}[1]
		\Statex \textbf{Input:} 
		Dimensions: $p,m,n_{\mathrm{ineq}}$; nonlinear objective function: $f^{\mathrm{nl}}(x_0,x_1)$;
		cost weight: $\alpha \geq  0$; proportion of nonlinear equality and inequality constraints: $s^{\mathrm{eq}}, s^{\mathrm{ineq}} \in [0,1]$; proportion of duplicated inequality constraints $s^d \in [0,1]$.
		% initial guess and bounds
		\Statex \textcolor{gray}{Initial guess and upper bounds.}
		\State Generate initial guess $x_0^0 \sim  \mathcal{N}_p(0,1)$ and $x_0^0 = |x_0^0|$.
		\State Generate initial guess $x_1^0 \sim  \mathcal{N}_m(0,1)$ and $x_1^0 = \max(0,x_1^0)$.
		\State Generate upper bounds $x_0^\mathrm{ub} \sim  \mathcal{U}_p(100,1000)$ and $x_1^\mathrm{ub} \sim  \mathcal{U}_m(100,1000)$.	
		\Statex \textcolor{gray}{Equality constraints.}
		\State Choose density for the random matrices $s \sim  \mathcal{U}_1(0.05,0.1)$.
		\State Generate random integer number $r \sim  \mathrm{round}(\mathcal{U}_1(1,m))$.
		\State Generate matrix $E \sim  \mathcal{U}_{r \times (n-r)}(2,3,s)$, vectors $d_1 \sim  \mathcal{U}_r(0,2)$, $d_2 \sim  \mathcal{U}_{m-r}(0,2)$, and let $D_1 = \mathrm{diag}(d_1)$, $D_2 = \mathrm{diag}(d_2)$,
		 $M = \begin{bmatrix}
			D_1 & E\\
			-E^\top  & D_2
		\end{bmatrix}$.
		\State Generate matrix $N \sim  \mathcal{U}_{m \times p}(-1,1,s)$ and vector $q \sim  \mathcal{U}_r(-20,-10)$,
		\State Define $c^{\mathrm{eq}}(x) = M x_1 + N x_0 + q - x_2 = 0$ and $x = (x_0,x_1,x_2) \in \R^n, n = p+2m$.
		\State Pick a random subset $\I^{\mathrm{eq}} \subseteq \{1,\ldots,m\}$ with cardinality $|\I^{\mathrm{eq}}| = \mathrm{round}(s^{\mathrm{eq}} m)$.
		\State Replace $c_i^{\mathrm{eq}}(x) = 0$ by $c_i^{\mathrm{eq}}(x)+c_i^{\mathrm{eq}}(x)^2+c_i^{\mathrm{eq}}(x)^4 = 0$ for all $i \in \I^{\mathrm{eq}}$.
		\Statex \textcolor{gray}{Inequality constraints}
		\State Generate $A \sim \mathcal{U}_{n_{\mathrm{ineq}} \times p }(0,1,s)$, $B = \mathcal{U}_{n_{\mathrm{ineq}} \times m }(0,1,s)$
		\State Generate $a = A x_0^0 + B x_1^0 - |\epsilon|$, $\epsilon\in \mathcal{N}_{n_{\mathrm{ineq}}}(0,1)$.
		\State Define inequality constraints $c^{\mathrm{ineq}}(x) = A x_0 + B x_1 - a \geq 0$. 
		\State Pick a random subset $\I^{\mathrm{ineq}} \subseteq \{1,\ldots,n_{\mathrm{ineq}}\}$ with cardinality $|\I^{\mathrm{ineq}}| = \mathrm{round}(s^{\mathrm{ineq}} n_{\mathrm{ineq}})$.
		\State Replace $c_i^{\mathrm{ineq}}(x) \geq 0$ by $c_i^{\mathrm{ineq}}(x)+c_i^{\mathrm{ineq}}(x)^2+c_i^{\mathrm{ineq}}(x)^4 \geq 0$ for all $i \in \I^{\mathrm{ineq}}$.
		\State Pick a random subset $\I^{\mathrm{d}} \subseteq \{1,\ldots,n_{\mathrm{ineq}}\}$ with cardinality $|\I^{\mathrm{d}}| = \mathrm{round}(s^{\mathrm{d}} n_{\mathrm{ineq}})$ and duplicate inequality constraints $c_i^{\mathrm{ineq}}(x) \geq 0$ for all $ i \in  \I^{\mathrm{d}}$.
		\Statex \textcolor{gray}{Objective function.}
		\State Generate $g_0 \in \mathcal{U}_{p}(0,1)$, $g_1\in \mathcal{U}_m(1,3)$.
		\State Define objective function $f(x) = f^\mathrm{nl}(x_0,x_1)  + \alpha  g_0^\top x_0 + \alpha g_1^\top x_1$.
	\end{algorithmic}
\end{algorithm}

Our implementation allows the replacement of any random scalar variable (e.g. density $s$) by fixed values. 
In the implementation, the ranges of values from where the random variables are drawn can also be adapted, but for ease of notation, we write in Algorithm~\ref{alg:random_mpecs} the values used in our experiments.

To generate a test problem set, we consider 25 objective functions and create three instances of different sizes for each problem, resulting in a total of 75 problems.
For the nonlinear objective functions , we use the following nonlinear test functions from~\cite{Gould2013}:
\texttt{Fletcher,
	McCormick,
	Powell,
	Rosenbrock,
	Raydan1,
	Raydan2,
	Diagonal4,
	Diagonal5,
	Extended\_Tridiagonal,
	Three\_Exponential\_Terms,
	Generalized\_PSC1,
	Fletcvb3,
	Bdqrtic,
	Tridia,
	EG2,
	Edensch,
	Indef,
	Cube,
	Bdexp,
	Genhumps,
	Arwhead,
	Quartc,
	Cosine,
	NCVXQP6,
	DIXCHLNV.
}
The dimensions of $p$ are unique random integers from the set $\{100,\ldots,2100\}$, $p = m$ and $n_{\mathrm{ineq}} = \mathrm{round}(\mathcal{U}_1(0.5 p, 2p))$.
The remaining parameters we use in Algorithm~\ref{alg:random_mpecs} are 
 $\alpha = 10$,
 $s^{\mathrm{eq}} = 0.01$,
 $s^{\mathrm{ineq}} =0.25$, and 
 $s^d =0.25$.

%\includebibliography
\clearpage
\bibliographystyle{plain}
\bibliography{syscop}

\end{document}